\documentclass[oneside]{amsart}

\usepackage[latin1]{inputenc}
\usepackage{amsmath}
\usepackage{amssymb}
\usepackage{graphicx}
\usepackage{color}
\usepackage[all]{xy}
\usepackage{epsfig}
\usepackage{float}
\usepackage{psfrag}
\usepackage{enumerate}
\usepackage[T1]{fontenc}
\usepackage[sc]{mathpazo}
\usepackage{tabularx}
\usepackage{colortbl}

\usepackage{colordvi}
\usepackage{array}

\newcommand{\N}{\mathbb{N}}
\newcommand{\Z}{\mathbb{Z}}
\newcommand{\Q}{\mathbb{Q}}
\newcommand{\R}{\mathbb{R}}
\newcommand{\C}{\mathbb{C}}
\newcommand{\F}{\mathbb{F}}
\newcommand{\X}{\mathbb{X}}

\newcommand{\Aa}{\mathcal{A}}

\newcommand{\Ss}{\mathcal{S}}

\newcommand{\Mm}{\mathcal{M}}

\newcommand{\Oo}{\mathcal{O}}
\newcommand{\Xx}{\mathcal{X}}

\newcommand{\jj}{\mathcal{J}}

\newcommand{\Uu}{\mathcal{U}}

\newcommand{\om}{\omega}

\newcommand{\omucc}{\omega_{\mu,c_1,c_2}}

\newcommand{\image}{\mathrm{Im}}
\newcommand{\codim}{\mathrm{codim}}
\newcommand{\coker}{\mathrm{coker}}
\newcommand{\delbar}{\bar\partial}

\newcommand{\tT}{\widetilde T}
\newcommand{\tDelta}{\widetilde\Delta}

\newcommand{\tjj}{\widetilde{\jj}}

\newcommand{\tjjucc}{\tjj_{\mu,c_1,c_2}}

\newcommand{\ccrit}{c_{\text{crit}}}

\newcommand{\Ssucc}{\Ss_{\mu,c_1,c_2}}
\newcommand{\Aaucc}{\Aa_{\mu,c_1,c_2}}
\newcommand{\Xxucc}{\Xx_{\mu,c_1,c_2}}
\newcommand{\uccp}{\mu',c_1',c_2'}
\newcommand{\ucc}{\mu,c_1,c_2}

\newcommand{\PD}{\mathrm{PD}}
\newcommand{\id}{\mathrm{id}}

\newcommand{\rk}{\operatorname{rk}}

\newcommand{\GL}{\mathrm{GL}}

\newcommand{\U}{\mathrm{U}}
\newcommand{\SO}{\mathrm{SO}}

\newcommand{\PSL}{\mathrm{PSL}}

\newcommand{\Sym}{\mathrm{Sym}}
\newcommand{\psl}{\mathrm{sl}}

\newcommand{\Symp}{\mathrm{Symp}}
\newcommand{\Diff}{\mathrm{Diff}}

\newcommand{\Aut}{\mathrm{Aut}}
\newcommand{\Iso}{\mathrm{Iso}}
\newcommand{\Hol}{\mathrm{Hol}}
\newcommand{\Map}{\mathrm{Map}}

\newcommand{\FDiff}{\mathrm{FDiff}}

\newcommand{\Gucc}{{G}_{\mu,c_1,c_2}}

\newcommand{\CP}{\mathbb{C}P}

\newcommand{\STS}{S^2\times S^2}
\newcommand{\PbP}{\CP^2\#\,\overline{\CP}\,\!^2}
\newcommand{\NTB}{S^{2}\widetilde{\times} S^{2}}
\newcommand{\tJ}{\widetilde{J}}

\newcommand{\tM}{\widetilde{M}}
\newcommand{\tMuc}{\tM_{\mu,c_1}}
\newcommand{\Muo}{M^0_{\mu}}
\newcommand{\Mul}{M^1_{\mu}}

\newcommand{\tMuco}{\widetilde{M}^0_{\mu,c_1}}

\newcommand{\tMucc}{\tM_{\mu,c_1,c_2}}
\newcommand{\CTC}{\CP^2\#\,3\overline{\CP}\,\!^2}

\newcommand{\into}{\hookrightarrow}

\theoremstyle{plain}
\newtheorem{thm}{Theorem}[section]
\newtheorem*{thm*}{Theorem}
\newtheorem{prop}[thm]{Proposition}
\newtheorem*{prop*}{Proposition}
\newtheorem{lemma}[thm]{Lemma}
\newtheorem*{lemma*}{Lemma}
\newtheorem{cor}[thm]{Corollary}
\newtheorem*{cor*}{Corollary}

\newtheorem*{conj*}{Conjecture}
\newtheorem{claim}[thm]{Claim}
\newtheorem*{claim*}{Claim}

\theoremstyle{definition}
\newtheorem{defn}[thm]{Definition}
\newtheorem*{defn*}{Definition}

\theoremstyle{remark}
\newtheorem{remark}[thm]{Remark}
\newtheorem*{remark*}{Remark}

\newtheorem*{remarks*}{Remarks}

\begin{document}

\title[Homotopy algebra of symplectomorphism groups]{The homotopy Lie algebra of symplectomorphism groups of 3--fold blow-ups of the projective plane}

\author[S. Anjos]{S\'ilvia Anjos}
\address{Center for Mathematical Analysis, Geometry and Dynamical Systems \\ Departamento de Matemática \\  Instituto Superior T\'ecnico \\  Av. Rovisco Pais \\ 1049-001 Lisboa \\ Portugal}
\email{sanjos@math.ist.utl.pt}
\thanks{Partially supported by FCT through program POCTI/FEDER and grant PTDC/MAT/098936/2008}

\author[M. Pinsonnault]{Martin Pinsonnault}
\address{The University of Western Ontario \\ London, Ontario, Canada N6A 3K7}
\email{mpinso@uwo.ca}
\thanks{Partially supported by NSERC grant RGPIN 371999.}

\begin{abstract}
By a result of Kedra and Pinsonnault, we know that the topology of groups of symplectomorphisms of symplectic 4-manifolds is complicated in general. However, in all known (very specific) examples, the rational cohomology rings of symplectomorphism groups are finitely generated. In this paper, we compute the rational homotopy Lie algebra of symplectomorphism groups of the  3--point blow-up of the projective plane (with an arbitrary symplectic form) and show that in some cases, depending on the sizes of the blow-ups, it is infinite dimensional. Moreover, we explain how the topology is generated by the toric structures one can put on the manifold. Our method involve the study of the space of almost complex structures compatible with the symplectic structure and it depends on the inflation technique of Lalonde--McDuff.\end{abstract}

\maketitle

\section{Introduction}
The study of the topology of the group of symplectomorphisms of symplectic 4--manifolds is by now a classic topic in symplectic topology, but still very incomplete. Indeed, the only closed 4--manifolds for which the symplectomorphism group is well understood are $\CP^2$, rational ruled symplectic 4--manifolds, their one point blow--ups and the 3-, 4- and 5-point blow--ups of the projective plane, considered as {\em monotone}\footnote{A symplectic manifold $(M, \omega)$ is monotone if $c_1(M)=k[\omega] \in H^2(M,\R)$ for some $k>0$.} symplectic manifolds. The most recent developments on this subject are given by Abreu, Granja and Kitchloo in ~\cite{AGK}, Anjos, Lalonde and Pinsonnault in ~\cite{ALP} and Evans in ~\cite{Ev} , respectively, where they give a complete description of the homotopy type of the symplectomorphism group of these manifolds.  In particular Evans showed that the symplectomorphism group of the monotone symplectic 3--point blow--up of $\CP^2$ is homotopy equivalent to a torus $T^2$.  Consider  the  symplectic  manifold $\X_n$ obtained from $\CP^2$ by performing $n \geq 1$ blow--ups of capacities $\delta_i$ where $i=1, \dots, n.$ The manifold $\X_n$ carries a family of symplectic forms $\omega_\nu$ where $\nu >0$ determines the cohomology class $[\omega_\nu]$. It is well known that for $n=1,2$ the rational cohomology algebra of $\Symp(\X_n, \omega_\nu)$ is finitely generated (see ~\cite{AGK}, ~\cite{AM}, ~\cite{ALP}). However, in certain cases, depending on the sizes of the blow-ups, that no longer holds if $n\geq 3$, as suggested by Kedra in \cite{Ke}. More precisely, Kedra and Pinsonnault  showed that  for a large class of symplectic 4--manifolds, namely, if $(M,\omega)$ is a simply connected symplectic 4--manifold   such that $b_2=\dim H^2(M) \geq 3$, then the rational cohomology algebra  $H^*(\Symp({\widetilde M}_\epsilon); \Q)$, of the symplectomorphism group of the symplectic blow--up of size $\epsilon$ of $(M,\omega)$, is infinitely generated provided $\epsilon > 0$ is sufficiently small. This does not mean that the topology of the symplectomorphism group cannot be understood. It indicates a more appropriate homotopy-theoretic invariant to investigate might be its homotopy algebra (with the Samelson product) or its homology ring (with the Pontryagin product).  Another motivation to study this problem comes from an earlier work of Seidel  \cite{Se} in which he proves that for the  monotone symplectic surface $\X_5$, the group $\Symp \cap \Diff_0 $ is not connected. 
In this paper we show that in the case of the 3--point blow--up the subgroup $\Symp_h(\X_3,\omega_\nu)$ of the group of symplectomorphisms of $(\X_3, \omega_\nu)$ acting trivially in homology is always connected (see Section \ref{Connexity}). We also compute  the homotopy Lie algebra $\pi_*(\Symp_h(\X_3, \omega_\nu))\otimes \Q$ and the Pontryagin ring of $\Symp_h(\X_3,\omega_\nu)$ with rational coefficients for all $\nu >0$.  Although the latter is always finitely generated, it follows from the computation of the rational homotopy algebra that in all cases except in some very particular ones (see (R3) in Section \ref{RemainingCases}) which include the monotone one, the rational cohomology algebra $H^*(\Symp(\X_3, \omega_\mu);\Q)$ is indeed infinitely generated.

Let $M_\mu$ be the symplectic manifold $S^2 \times S^2$ endowed with the split symplectic form with area $\mu \geq 1$ for the first $S^2$--factor, and with area 1 for the second factor,  and $\tMucc$ be obtained from $M_\mu$ by performing two successive symplectic blow-ups of capacities $c_1$ and $c_2$, with $\mu\geq 1> c_1+c_2>c_1>c_2$. In Section \ref{SymplecticCone} we explain why $\tMucc$ is symplectomorphic to $\X_3$ and why it is sufficient to consider values of $c_1$ and $c_2$ in this range.

Recall that, if $G$ is a connected topological group, the Samelson product $[\cdot,\cdot]: \pi_p(G) \otimes \pi_q(G) \to \pi_{p+q}(G)$ is defined by the commutator map

\begin{equation}\label{Samelsonproduct}
S^{p+q}=S^p \times S^q/ S^p \vee S^q \to G: (s,t) \mapsto a(s)b(t) a^{-1}(s) b^{-1}(t). 
\end{equation} 
Moreover, this product gives the rational homotopy of $G$ the structure of a graded Lie algebra, that is, $\pi_*(G) \otimes \Q$ is a graded vector space together with a linear map of degree  zero (the Samelson product)  satisfying two properties: antisymmetry  and the Jacobi identity.

\begin{figure}[htbp]
       \begin{center}
          \psfrag{G1}{\small $G_1${ \ \ {\fbox{$x_0,y_0$}}}}
          \psfrag{G2}{\small $G_2${ \ \ {\fbox{$y_0,x_1$}}}}
          \psfrag{G3}{\small{{\fbox{$x_1,y_1$}}} \ $G_3$}
          \psfrag{G4}{\small $G_4${ \ \ {\fbox{$y_1,z$}}}}
          \psfrag{G5}{\small $G_5${ \ \ {\fbox{$z,x_0$}}}}
          \psfrag{A1}{\small $A_1${  {\fbox{$y_0$}}}}
          \psfrag{A2}{\small $A_2${  {\fbox{$x_1$}}}}
          \psfrag{A3}{\small $A_3${  {\fbox{$y_1$}}}}
          \psfrag{A4}{\small $A_4${  {\fbox{$z$}}}}
          \psfrag{A5}{\small{{\fbox{$x_0$}}} \ $A_5$}
          \psfrag{i1}{\small $i_1$}
          \psfrag{i2}{\small $i_2$}
          \psfrag{i3}{\small $i_3$}
          \psfrag{i4}{\small $i_4$}
          \psfrag{i5}{\small $i_5$}
          \psfrag{j1}{\small $j_1$}
          \psfrag{j2}{\small $j_2$}
          \psfrag{j3}{\small $j_3$}
          \psfrag{j4}{\small $j_4$}
          \psfrag{j5}{\small $j_5$}
           \resizebox{!}{6.5cm}{\includegraphics{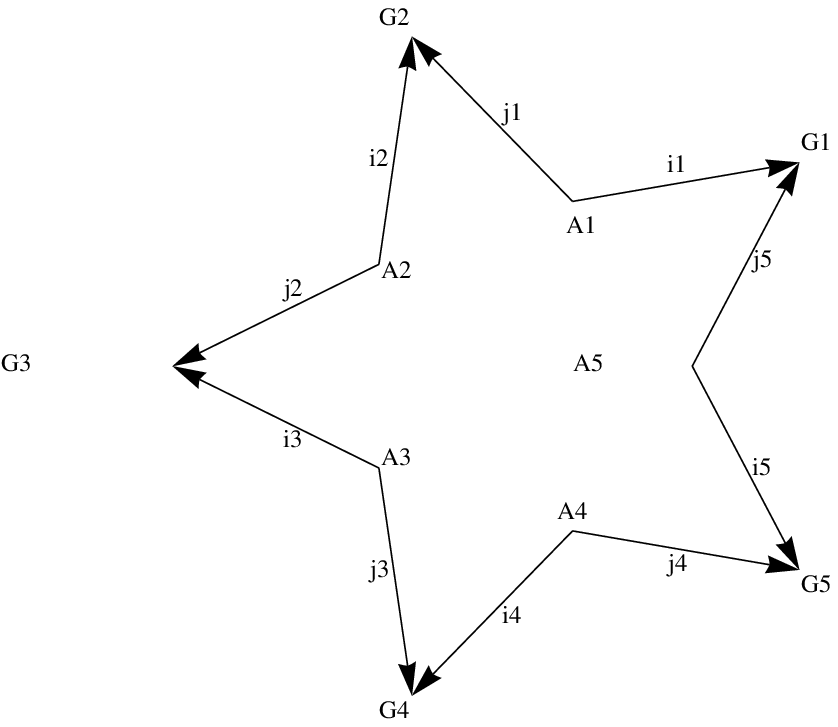}}
           \caption{Homotopy type of $\Symp(\tMucc)$ if $\mu=1$.}
           \label{HomType}
        \end{center}
    \end{figure}

Let $\Gucc$ be the group of symplectomorphisms of $\tMucc$ acting trivially in homology. One of the main results of this paper is the computation of the rational homotopy Lie algebra of the group $\Gucc$,  $\pi_*(\Gucc)\otimes \Q$, that we denote simply by $\pi_*$. Let us write $\mu=\ell + \lambda$ with $\ell \in \N$ and $0< \lambda \leq 1$. The first part of the main result says that if $\mu=1$ then $\pi_*$ is generated by five generators of degree 1,  $x_0,y_0,x_1,y_1,z$, such that all Samelson products between them vanish except  for $[x_0,y_1], [x_0,x_1], [y_0,y_1], [z,y_0]$ and $ [z,x_1]$. In Section \ref{conventions} we will see that these generators are represented by Hamiltonian $S^1$--actions contained in the  toric actions one can put in the symplectic manifold. It follows from the computation of the Pontryagin ring (Corollary  \ref{PontryaginRing}) that in the case $\mu=1$ the symplectomorphism group, as a topological group, has the homotopy type of a homotopy direct limit of Lie groups.  More precisely, it has the homotopy type of the group
$$G=\lim_{\longrightarrow}(G_k,A_k)$$ where $G_k=T^2$ and $A_k=S^1$ for all $1 \leq k \leq 5$ and  the homomorphisms $$i_k: A_k \to G_k, \quad 1\leq k \leq 5$$
$$j_k: A_k \to G_{k+1}, \quad 1\leq k \leq 4 \quad \mbox{and} \quad j_5:A_5 \to G_1$$
 are inclusions in the first and second factors respectively. This pushout is represented in Figure \ref{HomType} where the elements inside the boxes correspond to the generators of the toric actions. 
 
 When $\mu$ gets bigger than 1 the number of generators of $\pi_*$ increases.  More precisely, there is a new degree one generator and, assuming that $ \ell < \mu \leq \ell +1$ and  $\lambda \leq c_2 $ there are also two new generators of degree $4\ell-2$. Then, when $\lambda$ passes $c_2$ there is only one generator of the same degree, while when $\lambda$ passes $c_1$ this generator vanishes and two new generators of degree $4\ell$ appear. Finally, when $\lambda$ passes $c_1+c_2$ there is only one generator of degree $4\ell$ that vanishes when $\mu$ passes the next integer. 
The main theorem contains the complete statement.

\begin{thm}\label{maintheorem}
If   $\mu =1 > c_1 + c_2 > c_1 > c_2 $ then  $\pi_*$ is generated by  $x_0,y_0,x_1,y_1,z$, where all generators have degree 1 and all  Samelson products between them vanish except the ones in the following set 
$$ S= \{ \, [x_0,y_1], [x_0,x_1], [y_0,y_1], [z,y_0], [z,x_1] \, \}.$$

\noindent If $\mu= \ell+\lambda $, where $\ell \in \N$, $\lambda \in (0,1]$ and   $1> c_1 + c_2 > c_1 > c_2 $ then  $\pi_*$ is given as follows.

\begin{enumerate}
\item[(S1)]  If $1 < \mu \leq 2 $ and $\lambda \leq c_2$  then $\pi_*$ is generated by  $x_0,y_0,x_1,y_1,z,t$, where $\deg t=1$  and all Samelson products between them vanish except the ones in the set
$$ S_1= S \cup \{ \,   [x_1,t], [y_0,t]=[x_0,t],  [z,t]=[x_1,z]+[x_1,t]\, \}.$$

\item[(S2)]   if $1 < \mu \leq 2 $ and $c_2< \lambda \leq c_1$  then $\pi_*$ is generated by  $x_0,y_0,x_1,y_1,z,t$, where all Samelson products vanish except the ones in the set 
 $$ S_2= S \cup \{ \, [x_1,t], [y_0,t]=[x_0,t]=[x_1,t]+[x_1,x_0], [z,t]=[x_1,z]+[x_1,t]\,\}.$$
 
 \item[(S3)]  If 
 $(a) \ \ell > 1 \ {\mbox and} \ \lambda \leq c_2, \ {\mbox or} \quad   (b) \ \ell \geq 1   \ {\mbox and} \  c_1 < \lambda \leq c_1+c_2, $
 then $\pi_*$ is generated by  $x_0,y_0,x_1,y_1,z,t,w_{\ell_1},w_{\ell_2}$,  where 
 $$\deg w_{\ell_1}= \deg w_{\ell_2} = \left\{ \begin{array}{ll}
 4 \ell-2 & \mbox{in case }  (a)\\
 4\ell & \mbox{in case } (b),
\end{array} \right. 
 $$
 and all  Samelson products vanish except  the ones in the set 
 $$ S_3= S \cup \{ \,  [y_0,t]=[x_0,t]=[x_1,x_0], [z,t]=[x_1,z]\,\}.$$

\item[(S4)]  If 
$  (a) \ \ell > 1 \ {\mbox and} \ c_2<\lambda \leq c_1,  \ {\mbox or} \quad  (b) \  \ell \geq 1  \  {\mbox and}  \ c_1+c_2 < \lambda,$
 then $\pi_*$ is generated by  $x_0,y_0,x_1,y_1,z,t,w_{\ell}$,  where 
 $$\deg w_{\ell} = \left\{ \begin{array}{ll}
 4 \ell-2 & \mbox{in case } (a)\\
 4\ell & \mbox{in case }  (b),
\end{array} \right. 
 $$
 and all  Samelson products vanish except  the ones in the set $S_3$.
\end{enumerate}
\end{thm}

\noindent  Other important results regarding the homotopy type of the symplectomorphism group of $\tMucc$ are the following. 

\begin{prop}\label{stabilizer}
Consider $c_1,c_2 \in (0,1)$ such that either 
\begin{enumerate}
\item $c_2 < \lambda \leq  c_1 < c_1+c_2 < 1$; 
\item $c_2 < c_1 < c_1+c_2 < \lambda \leq 1$.
\end{enumerate}
Then  the group of symplectomorphisms $\Gucc$ is homotopy equivalent to the group $\Symp_p(\tMuc) \subset \Symp(\tMuc)$ of symplectomorphisms that fix a point $p$ in the manifold $\tMuc$, where  $\tMuc$  is obtained from $M_\mu$ by performing a symplectic blow--up of capacity $c_1$. 
\end{prop}

\begin{prop}\label{Inflation}
Given $c_2 < c_1 < c_2+c_1<1$, there is a homotopy coherent system of maps 
\[
\xymatrix{
\Gucc \ar[r] \ar[dr] & G_{\mu+\epsilon,c_1,c_2} \ar[d] \\
 & G_{\mu+\epsilon+\delta,c_1,c_2}}
\]
defined for all $\mu\geq 1$ and all $\epsilon, \delta >0$. Moreover, 
\begin{enumerate}
\item the homotopy type of $\Gucc$ remains constant as $\mu$ varies in either  $(\ell,\ell+c_2]$, $(\ell+c_2,\ell+c_1]$, $(\ell+c_1,\ell+c_1+c_2]$ or $(\ell+c_2+c_1,\ell+1]$.
\item The map $\Gucc\to G_{\mu+\epsilon,c_1,c_2}$ is $(4\ell-3)$-connected when $\lambda$ passes $c_1$ or $c_2$, and $(4\ell-1)$-connected when $\lambda$ passes $c_1+c_2$ or $\mu$ passes the integer $\ell+1$. In particular,  when $\mu>1$, it induces an isomorphism of fundamental groups.
\item These maps induce surjections $H^{*}(BG_{\mu+\epsilon,c_1,c_2})\to H^{*}(BG_{\mu,c_1,c_2})$ for all coefficients. Consequently, the map $B\Gucc\to BG_\infty$ induces a surjection in cohomology.
\end{enumerate}
\end{prop}

\noindent The analysis of the natural action of the symplectomorphism group $\Gucc$ on the (contractible!) space of compatible integrable complex structures $\jj_{\mu,c_{1},c_{2}}^{\text{int}}$ allows us to obtain an homotopy decomposition of the classifying space  $B\Gucc$ as an iterated homotopy pushout of classifying spaces of isometry subgroups. Given a $G$-space $X$, let write $X_{hG}$ for the corresponding homotopy orbit, that is,
\[X_{hG}:=EG\times_{G} X\]
Then, the homotopy type of $B\Gucc$ is given by

\begin{thm}\label{thm:HomotopyPushout}
If $\ell < \mu \leq \ell+1 $ where $\ell \in \N$ and  $c_2 < c_1 < c_2+c_1<1$  there is a homotopy pushout diagram
$$
\xymatrix{
S^{d}_{hIso} \ar[r] \ar[d]^{j_{\ell}} & B\, Iso \ar[d]^{i_{\ell}} \\
B{G}_{\mu',c_1,c_2} \ar[r] & B\Gucc}
$$
 where $\mu'= \ell + \lambda'$ lies in the interval of Proposition \ref{Inflation} (1) immediately before the interval where $\mu$ lies. The group $Iso$ is the group of isometries of certain complex structures on the manifold $\tMucc$.  The upper horizontal map  is the universal bundle map associated to the representations of the isometries $Iso$  on  the moduli space of infinitesimal deformations of the complex structure $J_\ell$, $H^{0,1}_{J_{\ell}}(T\tMucc)$, $i_{\ell}$ is induced by the inclusion $Iso \into \Gucc$, and the map $BG_{\mu',c_1,c_2}\to B\Gucc $ coincides, up to homotopy, with the one  described in Proposition ~\ref{Inflation}.
\end{thm}
\noindent We give a more precise statement of this theorem in Section \ref{HomotopyDecomposition}, where we explain exactly how the complex structures are given as well as their isometries.

 In order to prove Theorem \ref{maintheorem} we need a fundamental calculation which is  the computation of the rational cohomology of the classifying space of the group $\Gucc$.

\begin{thm} \label{CohomologyRingClassSpace}
If $\mu =1 > c_1 + c_2 > c_1 > c_2 $ then the rational cohomology ring of $B\Gucc$ is given by
$$ H^*(B\Gucc; \Q)= \Q[X_0,Y_0,X_1,Y_1,Z] / I, $$ where
all generators have degree 2 and $I=\langle \, X_0X_1,\,Y_0Y_1,\,X_0Y_1,\,Y_0Z,\,X_1Z \, \rangle$. 

If $\mu= \ell+\lambda $, where $\ell \in \N$, $\lambda \in (0,1]$ and   $1> c_1 + c_2 > c_1 > c_2 $ then 
$$ H^*(B\Gucc; \Q)= \Q[X_0,Y_0,X_1,Y_1,Z,T] / I \cup I_{\lambda,\ell}  $$
where $\deg T=2$ and $I_{\lambda,\ell}$ depends on $\ell \in \N$ and $\lambda \in (0,1]$ as follows. Let 
\begin{eqnarray*}
 A_k & = & k(X_0+X_1+k(Z+Y_0))+(k-1)(T+kY_1) \quad {\mbox  and } \\
 B_k & =  & k(X_0+X_1-k(Z+Y_0))+(k+1)(T-kY_1).
 \end{eqnarray*}
\begin{enumerate}[(a)]
\item If $\lambda \leq c_2$ then 
$$ I_{\lambda,\ell} = \langle \, T \prod_{k=1}^{\ell-1}(A_k\, B_k )\, (X_0+ \ell Y_0), \  T \prod_{k=1}^{\ell-1}(A_k\, B_k) \, [(\ell -1)(T+\ell Y_1)+\ell(X_1+\ell Z)] \, \rangle.$$
\item If $c_2 < \lambda \leq c_1$ then 
$$ I_{\lambda,\ell}= \langle \, T \prod_{k=1}^{\ell-1}(A_k\, B_k) \, A_\ell  \, \rangle.$$
\item If $c_1 < \lambda \leq c_1+c_2$ then 
$$ I_{\lambda,\ell}= \langle \, T \prod_{k=1}^{\ell-1}(A_k\, B_k) \, A_\ell (X_1+T-\ell Y_0), \ T \prod_{k=1}^{\ell-1}(A_k\, B_k) \, A_\ell [\ell X_0- \ell (\ell+1) Y_1 - \ell^2 Z +T] \,  \rangle  .$$
\item If $c_1+c_2 < \lambda$ then 
$$ I_{\lambda,\ell}= \langle \, T \prod_{k=1}^{\ell}(A_k\, B_k) \, \rangle.$$
\end{enumerate}
\end{thm}

\begin{remark}
One of the major differences from the previous cases that were studied in ~\cite{AGK} and ~\cite{ALP} comes from  the stratification of the space of almost complex structures compatible with the symplectic form. In our case we encountered a finer structure on this space since it is no longer sufficient to consider the stratification on the $B$-- ($J$--holomorphic) curves, where $B \in H_2(M)$ represents the class $[S^2 \times \{pt\}]$. One has to consider the stratification on $F$--curves as well, where $F$ represents the class $[\{pt\} \times S^2]$. This fact gives rise to new phenomena such as the existence of several strata of the same codimension and the addition of several strata at once to the space of almost complex structures each time $\mu$ crosses a critical value. Moreover, some of these strata have their K\"ahler isometry groups isomorphic to $S^{1}$ instead of a full $T^2$ as before. Another new feature is that for some values of $\mu$ the group of symplectomorphisms $\Gucc$ is no longer homotopy equivalent to the stabilizer of a point in $\X_{2}$. 
\end{remark}

{\bf Organization of the paper:} The body of the paper is divided into five sections and two Appendices.  In Section \ref{Jholomorphic} we first explain the general setting and then we study the structure of $J$--holomorphic curves in $\tMucc$. The third section is devoted to the study of the homotopy type of the symplectomorphism group $\Gucc$. More precisely, we state the main result about the stability of the symplectomorphism groups, whose proof is delayed to Appendix \ref{se:appendixA}, and we present the main ideas of the proofs of Theorems \ref{Inflation} and \ref{thm:HomotopyPushout}, that follow closely the proofs of Theorems 2.2 and 2.5 in ~\cite{ALP}. In Section \ref{se:LieAlgebra}  we prove Theorem \ref{maintheorem}. We begin this section with the computation of the homotopy groups of $\Gucc$ and recall some facts about rational homotopy theory. At the end of this section we prove Theorem \ref{CohomologyRingClassSpace}. Note that up to this point  we just study the topology of symplectomorphism groups $\Gucc$ in the generic case, that is, when the sizes of the blow-ups satisfy $0<c_2 < c_1 < c_2+c_1 <1$. In Section 5 we give some remarks about the remaining cases.  Appendix B is devoted to a brief explanation of the differential and topological results relative to the stratification of the space $\jj_{\om}:=\tjjucc$ of compatible almost complex structures on $\tMucc$, and to the stratification it induces on the subspace $\jj_{\om}^{int}\subset \jj_{\om}$ of compatible, integrable, complex structures. 

{\bf Acknowledgements.} The first author is grateful to Jarek Kedra for helpful comments and discussions. Both authors would like to thank Gustavo Granja for explaining some facts about the algebraic computations and for clarifying the construction of $A_{\infty}$ replacements for smooth actions. Part of this paper was completed while the second author stayed at MSRI in Spring 2010, and he is very grateful for its hospitality and support. We also would like to thank Olguta Buse for showing and explaining to us her results in a preliminary version of \cite{Bu} about negative inflation. Finally, the authors warmly thank the referee for reading the paper carefully and making several important suggestions that have helped to clarify various arguments.

\section{Symplectic forms on $\X_{3}$ and structure of $J$--holomorphic curves}
\label{Jholomorphic}


\subsection{Symplectic cones of rational surfaces and reduced classes}\label{SymplecticCone}

Since diffeomorphic symplectic forms define homeomorphic symplectomorphism groups, and since symplectomorphisms are invariant under rescalings of symplectic forms, we first describe a fundamental domain for the action of $\Diff\times \R_{*}$ on the space $\Omega_{+}$ of orientation-compatible symplectic forms defined on the $3$-fold blowup $\X_{3}=\CTC$.

Let $\om_{\nu}$ be the standard Fubini-Study form on $\CP^{2}$ rescaled so that $\om_{\nu}(\CP^{1})=\nu$. Given $1\leq n\leq 3$, we write $\om_{\nu;\delta_{1},\ldots,\delta_{n}}$ for the symplectic form on $\X_{n}$ obtained from the $n$-fold symplectic blow-up of $(\CP^{2},\om_{\nu})$ at $n$ disjoint balls of capacities $\delta_{1},\dots,\delta_{n}$. If $\{L,V_{1},\cdots,V_{n}\}$ is the standard basis of $H_{2}(\X_{n};\Z)$ consisting of the class $L$ of a line together with  the classes $V_{i}$ of the exceptional divisors, the Poincaré dual of $[\om_{\nu;\delta_{1},\ldots,\delta_{n}}]$ is then $\nu L-\sum_{i}\delta_{i}V_{i}$. Similarly, given any compatible almost-complex structure $J$ on $\X_{n}$, the first Chern class $c_{1}:=c_{1}(J)\in H^{2}(\X_{n};\Z)$ is Poincaré dual to $K:=3L-\sum_{i} V_{i}$. Let $\mathcal{C}$ be the Poincaré dual of the symplectic cone of $\X_{n}$, that is,
\[
  \mathcal{C} =
  \{A\in H_{2}(\X_{n};\R)~|~A=\PD[\om]\text{~for some~}\om\in\Omega_{+}\},
\]
Recall that the ``uniqueness of symplectic blow-ups'' Theorem proved by McDuff~\cite{MD:Isotopie} implies that the diffeomorphism class of the form $\om_{\nu;\delta_{1},\ldots,\delta_{n}}$ only depends on its cohomology class. Hence, it is sufficient for us to describe a fundamental domain for the action of $\Diff\times \R_{*}$ on $\mathcal{C}$. In fact, since the canonical class $K$ is unique up to 
orientation preserving diffeomorphisms ~\cite{LL:Uniqueness}, we only need to describe the action of the subgroup $\Diff_{K}\subset\Diff$ of diffeomorphisms fixing $K$ on the $K$-symplectic cone
\[
  \mathcal{C}_{K} =
  \{A\in H_{2}(\X_{n};\R)~|~A=\PD[\om]\text{~for some~}\om\in\Omega_{K}\},
\]
where $\Omega_{K}$ is the set of orientation-compatible symplectic forms with $K$ as the symplectic canonical class. 

The key ingredient to understand that action is the notion of a \emph{reduced class}:

\begin{defn}\label{def:ReducedClasses}
Let $n\geq 3$. A class $A=a_{0}L-\sum_{i} a_{i}V_{i}$ is said to be reduced with respect to the basis $\{L, V_{1}, \ldots, V_{n}\}$ if $a_{1}\geq a_{2}\geq\cdots\geq a_{n}\geq 0$ and $a_{0} \geq a_{1}+a_{2}+a_{3}$.
\end{defn}

\begin{thm}[Li--Liu~\cite{LL2}]\label{thm:FundamentalDomain} 
The set $\mathcal{C}_{0}$ of reduced classes $a_{0}L-\sum_{i}a_{i}V_{i}$ with $a_{n}>0$ is a fundamental domain of $\mathcal{C}_{K}(\X_{n})$ under the action of $\Diff_{K}$.
\end{thm}

A more convenient description of the fundamental domain $\mathcal{C}_{0}$ can be given by viewing $\X_{3}$ as the $2$-fold blowup of $S^{2}\times S^{2}$. Recall that by the Classification Theorem of Lalonde and McDuff \cite{LM}, any  rational ruled  4--manifold  is symplectomorphic, after  rescaling, to either

\begin{itemize}

\item the topologically trivial $S^2$--bundle over $S^2$, $M_\mu^0= (S^2 \times S^2, \omega^0_\mu)$, where $\omega^0_\mu$ is the split symplectic form with area $\mu \geq 1 $ for  the first $S^2$--factor, and with  area 1 for the second factor; or

\item the topologically non--trivial $S^2$--bundle over $S^2$, $M_\mu^1=(\PbP, \omega_\mu^1)$, where the symplectic area of the exceptional divisor is $\mu>0$ and the area of the projective line is $\mu+1$ (this implies that the area of the fiber is 1). 

\end{itemize}

It is well known that blowing up $M^0=S^2 \times S^2$ or $M^1=\PbP$ leads to diffeomorphic smooth manifolds ${\widetilde M}^0 \simeq{\widetilde M}^1=\X_{2}$. Hence, blowing up once more yields $\X_{3}$. Let $\{B, F, E_{1}, E_{2}\}$ be the natural basis of $H_2(M^{0}\#\,2\overline{\CP}\,\!^2)$. The diffeomorphism $\X_{3} \to M^{0}\#\,2\overline{\CP}\,\!^2$ can be chosen so that the ordered basis $\{L,~ V_{1},~ V_{2},~ V_{3}\}$ is identified with $\{B+F-E_{1},~ B-E_{1},~ F-E_{1},~E_{2}\}$. 
When one considers this birational equivalence in the symplectic category, the uniqueness of symplectic blow-ups implies that $(\X_{3}, \om_{\nu; \delta_{1}, \delta_{2}, \delta_{3}})$ is symplectomorphic, after rescaling, to $M_{\mu}^{0}$ blown--up with capacities $c_{1}$ and $c_{2}$, where $\mu=(\nu-\delta_{2})/(\nu-\delta_{1})$, $c_{1}=(\nu-\delta_{1}-\delta_{2})/(\nu-\delta_{1})$, and $c_{2}=\delta_{3}/(\nu-\delta_{1})$. It is then easy to see that the normalized homology class $\PD[\om_{\mu,c_{1},c_{2}}]=\mu B + F - c_{1}E_{1}-c_{2}E_{2}$ is reduced if, and only if, $0<c_{2}\leq c_{1}<c_{1}+c_{2}\leq 1$. This proves the following lemma.

\begin{lemma}\label{le:standard}
Every symplectic form on $\X_3$ is, after rescaling, diffeomorphic to a form Poincaré dual to $\mu B+F-c_{1}E_{1}-c_{2}E_{2}$ with $c_{2}\leq c_{1}<c_{1}+c_{2}\leq 1$.
\end{lemma}


\subsubsection{Gromov Invariants}
Let $(X, \omega)$ be a symplectic 4--manifold and $A \in H_2(X;\Z)$ be a homology class such that $k(A)\equiv \frac12(A\cdot A +c_1(A)) \geq 0$. Then the Gromov invariant of $A$ defined by Taubes in \cite{Ta} counts, for a generic complex structure $J$ tamed by $\omega$, the algebraic number of embedded $J$--holomorphic curves in class $A$ passing through $k(A)$ generic points.  Taubes's curves arise as zero sets of sections and so need not be connected. Therefore components of non-negative self-intersection are embedded holomorphic curves of some genus $g \geq 0$, while all other components are exceptional spheres. It follows from Gromov's Compactness Theorem that, given a tamed $J$, any class $A$ with non-zero Gromov invariant $\rm{Gr}(A)$, is represented by a collection of $J$--holomorphic curves or cusp-curves. If $(X, \omega)$ is some blow-up of a rational ruled surface $M_\mu^i$, with canonical homology class $K=-{\rm PD}(c_1)$, then, since $b_2^+=1$, it follows from the wall-crossing formula of Li and Liu for Seiberg--Witten invariants (\cite[Corollary 1.4]{LL}) that for every class $A \in H_2(X;\Z)$  verifying the condition $k(A) \geq 0$ we have 
$$ \rm{Gr}(A)= \pm  \rm{Gr}(K-A)=\pm 1.$$
In particular, if $\rm{Gr}(K-A)=0$ (for example, when $\omega(K-A) \leq 0$), then $\rm{Gr}(A)\neq 0.$  

\subsection{Structure of $J$--holomorphic curves in $\tMucc^i$ }
\label{StructureJholomorphic}
Since we can always obtain $\tMucc^i$ by blowing up a ball in the product $\STS$, we will work only with $\Muo$ and, to simplify the notation, we will omit the superscript ``$^0$'' on related objects. 
Let $\tjjucc$ be the set of almost complex structures in $\tMucc$ compatible with $\omega_\mu$. 

\begin{lemma}
\label{lm:classes}
Given $\mu \geq 1 > c_1+c_2 > c_1 >c_2$, if a class $A=pB+qF-r_1E_1-r_2E_2 \in H_2(\tMucc;\Z)$ has a simple $J$--holomorphic representative, then $p \geq 0$. Moreover, if $p=0$, then $A$ is either $F,F-E_1, F-E_2, F-E_1-E_2$ or $E_1,E_2, E_1-E_2$. Furthermore, if $p=1$ then $r_1,r_2 \in \{ 0,1 \}$. 
\end{lemma}
\begin{proof}
The fact that $p\geq 0$ follows from the adjunction inequality 
$$2g_v(A)=2(p-1)(q-1)-r_1(r_1-1)-r_2(r_2-1) \geq 0,$$ 
the positivity of the symplectic area $\omega(A)=\mu p+ q -c_1r_1-c_2r_2$, and our choice of normalization: $\mu=\omega(B)\geq \omega (F)=1 > \omega(E_1) =c_1 >\omega(E_2) =c_2$. Indeed, if we assume that $p<0$ then $p < \frac12$ and 
\begin{eqnarray*}
 -2g_v(A) & > & (q-1) + r_1(r_1-1)+ r_2(r_2-1)  \\
         & \geq & (c_1r_1+c_2r_2 -\mu p -1)+r_1(r_1-1)+ r_2(r_2-1) \\
         & \geq & (c_1r_1+c_2r_2) +r_1(r_1-1)+ r_2(r_2-1) \\
         & = & r_1(r_1-1+ c_1) + r_2(r_2-1+ c_2)
\end{eqnarray*}
Since $r_i$ is an integer we have $ r_i(r_i-1+ c_i) \geq 0$ and therefore we would have $g_v(A) <0$, which is a contradiction. 

If $p=0$ then the adjunction inequality and the positivity of intersections give
$2(q-1)+r_1(r_1-1)+r_2(r_2-1) \leq 0$ and $q -c_1r_1-c_2r_2 \geq0$, respectively. The values of $q,r_1,r_2$ that satisfy these conditions are the following:
\begin{itemize}
\item $q=1 \ \wedge \ r_1=0 \ \Rightarrow \ r_2=0  \ \vee \ r_2=1 \ \Rightarrow \ A=F \ \vee \ A=F-E_2$;
\item $q=1 \ \wedge \ r_1=1 \ \Rightarrow \ r_2=0  \ \vee \ r_2=1 \ \Rightarrow \ A=F-E_1 \ \vee \ A=F-E_1-E_2$;
\item $q=0 \ \wedge \ r_1=0 \  \Rightarrow \ r_2=-1  \ \Rightarrow \ A=E_2;$
\item $q=0 \ \wedge \ r_1=-1 \ \Rightarrow \ r_2=0  \ \vee \ r_2=1 \ \Rightarrow \ A=E_1 \ \vee \ A=E_1-E_2$.
\end{itemize}
The last statement follows easily from the adjunction inequality.
\end{proof}

\begin{lemma}\label{lm:E2}
Given $\mu \geq 1 > c_1+c_2 > c_1 >c_2$ and any tamed almost complex structure $J \in \tjjucc$, the exceptional class $E_2$ is represented by a unique embedded $J$--holomorphic curve. Consequently, if a class $A=pB+qF-r_1E_1-r_2E_2$ has a simple $J$--holomorphic representative, then $r_2 \geq 0$.
\end{lemma}
\begin{proof}
We know that the set of almost complex structures $J$ for which the exceptional classes have embedded $J$--holomorphic representatives is open and dense. Moreover, by positivity of intersections, any such $J$--holomorphic representative is necessarily unique, and it is regular by the Hofer--Lizan--Sikorav regularity criterion. We must show it cannot degenerate to a cusp--curve. Suppose it is represented by a cusp curve $C= \bigcup_1^N m_iC_i$, where $N \geq 2$ and each component $C_i$ is simple. By Lemma \ref{lm:classes} the homology class of each component $C_i$ is either $F$, $F-E_1$, $F-E_2$, $F-E_1-E_2$ or  $E_1,E_2, E_1-E_2$. The cases $[C_i]=F,[C_i]=E_1,[C_i]=F-E_1$ and $[C_i]=E_2$ are impossible since the symplectic area of each component must be less than $\omega(E_2)=c_2$. But then each component $C_i$  can only represent one of the following classes: $E_1-E_2, F-E_2,F-E_1-E_2 $, which is also impossible, because $m_i>0$.  
\end{proof}

\begin{remark}\label{exceptionalclasses}
The exceptional classes are $E_2, E_1, F-E_1, F-E_2,B-E_1$ and $B-E_2$. Of all these, only $E_2$ cannot be represented by a cusp--curve. Note that $(E_1-E_2)\cdot (E_1-E_2)=(F-E_1-E_2 )\cdot(F-E_1-E_2 )=-2 $  so we cannot guarantee that these two classes always have a $J$--holomorphic representative. 
\end{remark}

\begin{remark}
If the class $E_1-E_2$ is represented by a embedded $J$--holomorphic curve (which is necessarily unique by positivity of intersections) then if a class $A=pB+qF-r_1E_1-r_2E_2$ has also a simple $J$--holomorphic representative, then $r_1 \geq r_2 \geq 0$. If the class $F-E_1-E_2$ is represented then $p \geq r_1+r_2  \geq 0$.
\end{remark}

\begin{lemma}\label{lm:F}
Given $\mu \geq 1 > c_1+c_2 > c_1 >c_2$ and a tamed almost complex structure $ J \in \tjjucc$, the moduli space ${\mathcal{M}}(F, \jj)$ of embedded $J$--holomorphic spheres representing the fiber class $F$ is a non-empty, open, 2--dimensional manifold.  
\end{lemma}
\begin{proof}
This was proved in  ~\cite[Lemma 2.3]{Pi}. 
\end{proof}
Let us define the classes $D_i \in H_2(\tMucc; \Z)$, $i \in \Z$, by setting $D_{4k+1}=B+kF$, $D_{4k}=B+kF-E_2$, $D_{4k-1}=B+kF-E_1$ and $D_{4k-2}=B+kF-E_1-E_2$. Note  that we have $$D_i \cdot D_i= \left\{ \begin{array}{cl}
                         \frac{i}2-1 & \mbox{if $i$ is even} \\
                         \frac{i-1}2 & \mbox{if $i$ is odd} 
                         \end{array} \right.  \ \ \mbox{ and } \ \
c_1(D_i)= \left\{ \begin{array}{cl}
                         \frac{i}2+1 & \mbox{if $i$ is even} \\
                         \frac{i+1}2+1 & \mbox{if $i$ is odd.} 
                         \end{array} \right. 
$$
Hence the adjunction formula $g_v(D_i)=1+\frac12(D_i \cdot D_i-c_1(D_i))=0$ implies that any $J$-holomorphic sphere in class $D_i$ must be embedded. Note also that the virtual dimension of the spaces of unparameterized $J$--curves representing $D_i$ is 
$$\dim {\mathcal{M}}(D_i,J)= \left\{ \begin{array}{cl}
                         i & \mbox{if $i$ is even} \\
                         i+1 & \mbox{if $i$ is odd} 
                         \end{array} \right.  \ \ \mbox{ and } \ \
k(D_i)=\frac12(D_i \cdot D_i+c_1(D_i))= \left\{ \begin{array}{cl}
                         \frac{i}2 & \mbox{if $i$ is even} \\
                         \frac{i+1}2 & \mbox{if $i$ is odd.} 
                         \end{array} \right. 
$$

\begin{lemma}\label{classDi}
Given an integer $i \geq 1$, the Gromov invariant of the class $D_i$ is nonzero. Consequently, for a generic tamed almost complex structure $J$, there exist an embedded $J$--holomorphic sphere representing $D_i$ and passing through a generic set of $k(D_i)$ points in $\tMucc$.
\end{lemma}
\begin{proof}
We mimic the proof of Lemma 2.4 in ~\cite{Pi}. Now the Poincaré dual of the canonical class of $\tMucc$ is $K=-2(B+F)-E_1-E_2$. As in Pinsonnault's paper we can conclude that $\rm{Gr}(D_i)=\pm 1$ and thus, given a generic $J$, $D_i$ is represented by the disjoint union of some embedded  $J$--holomorphic curve of genus $g \geq 0$ with some exceptional $J$--spheres. Since $D_i$ has nonnegative  intersection with the six exceptional classes $E_2, E_1, F-E_1, F-E_2,B-E_1$ and $B-E_2$, this representative is in fact connected and the adjunction formula gives its genus $g=0$.  
\end{proof}

\begin{lemma}\label{lm:decomposition}
The set of tamed almost complex structures on $\tMucc$ for which the classes $D_0=B-E_2$ and $D_{-1}=B-E_1$ are represented by an embedded J-holomorphic sphere is open in dense in $\tjjucc$. If for a given $J$ there are no such spheres, then either the class $D_{-2}=B-E_1-E_2$ is represented by a unique embedded J-holomorphic sphere or there is a unique integer $1 \leq m \leq \ell$ such that one of the classes $D_{-4m+1}=B-mF$, $D_{-4m}=B-mF-E_2$, $D_{-4m-1}=B-mF-E_1$ or $D_{-4m-2}=B-mF-E_1-E_2$ is represented by a unique embedded $J$--holomorphic sphere. 
\end{lemma}
\begin{proof}
Again we adapt the proof of the corresponding lemma in ~\cite{Pi}(Lemma 2.5) to our case. Since the classes $D_0=B-E_2$ and $D_{-1}=B-E_1$ are exceptional, the set of almost complex structures $J$ for which they have an embedded $J$--holomorphic representative is open and dense. Since they are primitive, it follows they have no multiply-covered representatives. Suppose they are represented by a cusp curve $C=\bigcup_{i=1}^Nm_iC_i$, where $N\geq 2.$ By Lemma \ref{lm:classes} and Lemma \ref{lm:E2}, any cusp curve gives a decomposition of the class $B-E_i$ of the form 

\begin{equation}\label{decomposition}
\begin{split} 
B-E_i  \ =  \ (B+mF-r_1E_1- & r_2E_2)+ \sum_{a_i}s_{a_i}E_i+\sum_{b_j}t_{b_j}(F-E_j)+\sum_cu_cF+ \\
& + \sum_dw_d(E_1-E_2)+\sum_ez_e(F-E_1-E_2),
\end{split}
\end{equation}
where the multiplicities $s_{a_i},t_{b_j},u_c,w_d,z_e$ are nonnegative and where $r_1,r_2 \in \{0,1\}$. This implies that $m \leq 0$. If $m=0$ then the class $B-E_1-E_2$ must be  represented by a unique embedded J-holomorphic sphere. If $|m|<0$, since the symplectic area of the component in class $B+mF-r_1E_1-r_2E_2$ must be strictly positive, $|m|$ must be less or equal than $\ell$. Moreover, if $|m|< \ell$, the classes $B-|m|F-E_2$, $B-|m|F-E_1$ and $B-|m|F-E_1-E_2$ have positive area since $c_1+c_2<1$, while if $|m|=\ell$, they have positive area only if $c_2 < \lambda$, $c_1 < \lambda$ and  $c_1+c_2 < \lambda$ respectively. Finally, note that the intersection of  $B-|m|F-r_1E_1-r_2E_2$ with any other class $D_i$, such that $i\leq 0$, is negative, therefore at most one such class can be represented by a $J$--holomorphic curve for a given tamed almost complex structure $J$. 
\end{proof}

\begin{figure}[htbp]
       \begin{center}
        \psfrag{A}{$B-E_2$}
          \psfrag{B}{$E_2$}
          \psfrag{C}{$F-E_2$}
          \psfrag{D}{$B-E_1$}
          \psfrag{E}{$E_1$}
          \psfrag{F}{$F-E_1$}
          \psfrag{G}{$B-E_2$}
          \psfrag{H}{$E_2$}
          \psfrag{I}{$F-E_2-E_1$}
          \psfrag{J}{$E_1$}
          \psfrag{K}{$B-E_1$}
          \psfrag{L}{$F-E_1$}
          \psfrag{M}{$E_1$}
          \psfrag{N}{$B-E_2-E_1$}
          \psfrag{O}{$E_2$}
          \psfrag{P}{$F-E_2$}
          \psfrag{Q}{$B-E_1$}
          \psfrag{T}{$F-E_1-E_2$}
          \psfrag{S}{$E_2$}
          \psfrag{R}{$E_1-E_2$}
          \psfrag{U}{$F-E_1$}
          \psfrag{V}{$E_2$}
          \psfrag{X}{$E_1-E_2$}
          \psfrag{Z}{$B-E_1$}
          \psfrag{AA}{$B-E_1-E_2$}
          \psfrag{BB}{$E_2$}
          \psfrag{CC}{$E_1-E_2$}
          \psfrag{DD}{$F-E_1$}
           \resizebox{!}{7.5cm}{\includegraphics{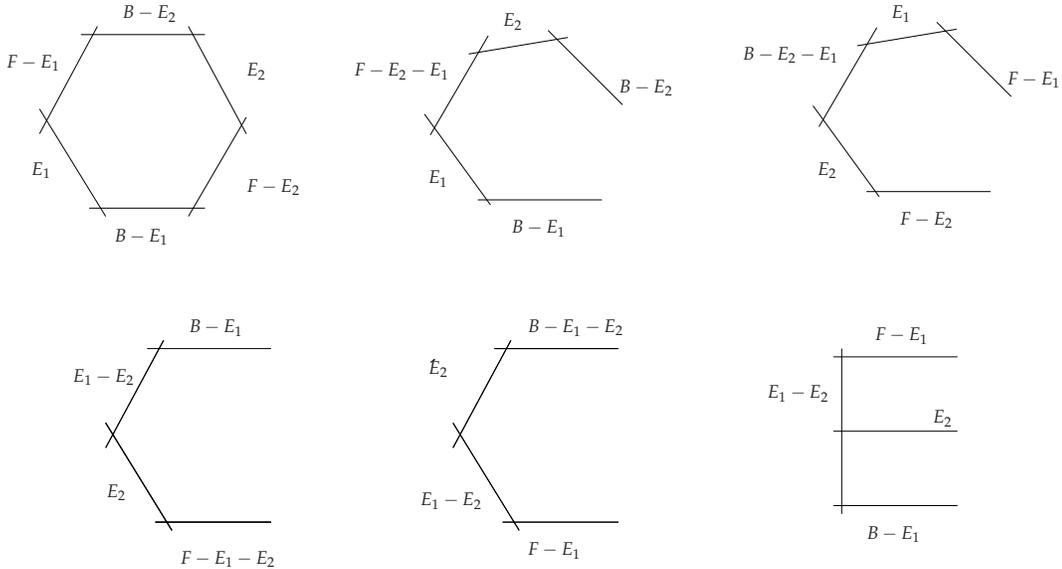}}
           \caption{Configurations of $J$--holomorphic spheres for $\mu=1$: (1) to (6).}
           \label{mu=1}
        \end{center}
      
    \end{figure}

Now note that in the decomposition \eqref{decomposition} of $B-E_i$ some of the multiplicities  $s_{a_i},t_{b_j},u_c,w_d,z_e$ can be zero. Indeed, for instance, if $t_{b_j}>0$ then necessarily we will have $z_e=0$ by positivity of intersections, that is, the classes $F-E_j$ and $F-E_1-E_2$ cannot be represented by $J$-holomorphic spheres with the same $J$. Therefore, if $\mu=1$, using the Lemmas in this section and positivity of intersections, we conclude that the configurations  we can obtain are the ones in Figure \ref{mu=1}.

\medskip

For values of $\mu >1$ more strata appear in the space of almost complex structures $\tjjucc$. The possible configurations 
in this case are the following:

\smallskip

\begin{itemize}
\item (7) to (10) in Figure \ref{configB-MF}, if the class $B-mF$ is represented by a $J$--holomorphic sphere;

 \begin{figure}[htbp]
        \begin{center}
          \psfrag{B}{$B-mF$}
          \psfrag{F1}{$F-E_1$}
          \psfrag{E12}{$E_1-E_2$}
          \psfrag{E2}{$E_2$}
          \psfrag{E1}{$E_1$}
          \psfrag{F12}{$F-E_1-E_2$}
          \psfrag{F2}{$F-E_2$}
           \resizebox{!}{3cm}{\includegraphics{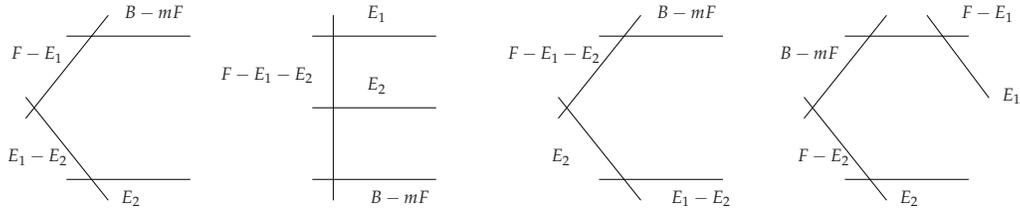}}
           \caption{Configurations (7) to (10)}
           \label{configB-MF}
        \end{center}
 
    \end{figure}

\item (11) and (12) in Figure \ref{configB-MF-E2}, if the class $B-mF-E_2$ is represented by a $J$--holomorphic sphere;

 \begin{figure}[htbp]
        \begin{center}
          \psfrag{A}{$E_1$}
          \psfrag{D}{$B-mF-E_2$}
          \psfrag{E}{$E_2$}
          \psfrag{F}{$F-E_1-E_2$}
          \psfrag{G}{$E_1$}
          \psfrag{H}{$F-E_1$}
          \psfrag{I}{$B-mF-E_2$}
          \psfrag{J}{$E_2$}
          \psfrag{K}{$F-E_2$}
           \resizebox{!}{3.5cm}{\includegraphics{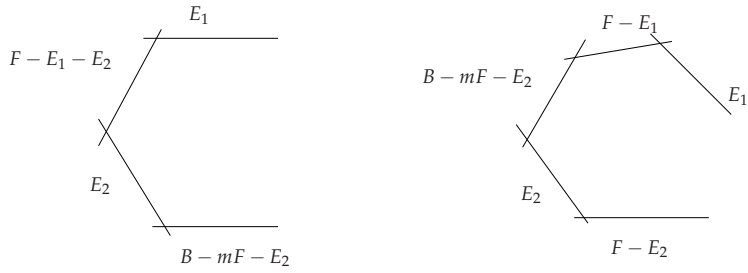}}
         \caption{Configurations (11) and (12)}
          \label{configB-MF-E2}
        \end{center}
     
    \end{figure}

\item (13) to (16) in Figure \ref{configB-MF-E1}, if the class $B-mF-E_1$ is represented by a $J$--holomorphic sphere;

\begin{figure}[htbp]
        \begin{center}
          \psfrag{B}{$B-mF-E_1$}
          \psfrag{F1}{$F-E_1$}
          \psfrag{E12}{$E_1-E_2$}
          \psfrag{E2}{$E_2$}
          \psfrag{E1}{$E_1$}
          \psfrag{F12}{$F-E_1-E_2$}
          \psfrag{F2}{$F-E_2$}
           \resizebox{!}{3cm}{\includegraphics{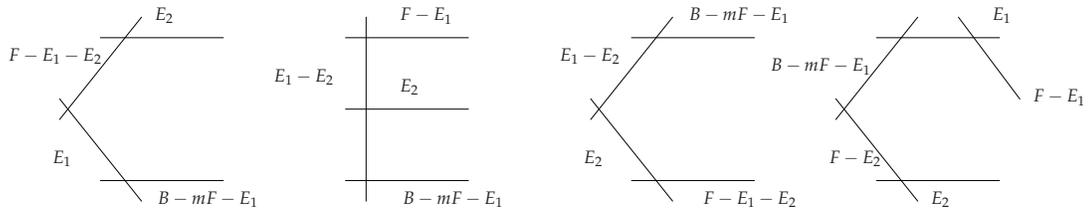}}
           \caption{Configurations (13) to (16)}
           \label{configB-MF-E1}
        \end{center}
    \end{figure}
    
\item and finally, configurations (17) and (18) in Figure \ref{configB-MF-E12}, if the class $B-mF-E_1-E_2$ is represented by a $J$--holomorphic sphere.

 \begin{figure}[htbp]
        \begin{center}
          \psfrag{A}{$F-E_1$}
          \psfrag{D}{$B-mF-E_1-E_2$}
          \psfrag{E}{$E_2$}
          \psfrag{F}{$E_1-E_2$}
          \psfrag{G}{$F-E_1$}
          \psfrag{H}{$E_1$}
          \psfrag{I}{$B-mF-E_1-E_2$}
          \psfrag{J}{$E_2$}
          \psfrag{K}{$F-E_2$}
           \resizebox{!}{3.5cm}{\includegraphics{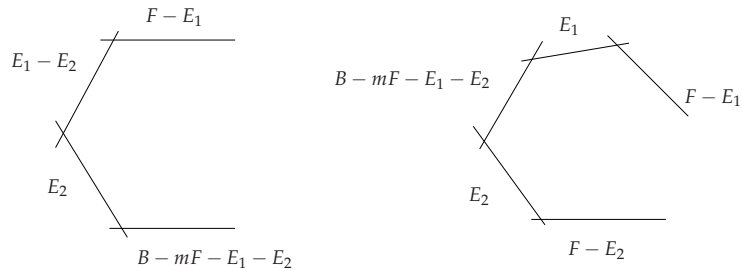}}
         \caption{Configurations (17) and (18)}
          \label{configB-MF-E12}
        \end{center}
    \end{figure}

\end{itemize}

\bigskip

Let $\tjj_{\ucc,m}$ be the space of almost complex structures $J \in \tjjucc$ such that the class $D_{-m}$, $m \geq 0$, is represented by an embedded $J$--holomorphic sphere. Lemma \ref{lm:decomposition} and the above considerations show that the space 
$\tjjucc$ is the disjoint union of its subspaces $\tjj_{\ucc,m}$, that is

\begin{equation}\label{strata}
\tjj_{\ucc}=\tjj_{\ucc,1}\sqcup \hdots \sqcup \tjj_{\ucc,N} \mbox{ where } N= 
\left\{ \begin{array}{cl}
4\ell -1 & \mbox{if}\  \ \lambda \leq c_2 < c_1 < c_1+c_2, \\ 
4\ell & \mbox{if} \ \ c_2 <\lambda \leq c_1 < c_1+c_2, \\
4\ell+1  & \mbox{if} \ \ c_2 < c_1 < \lambda \leq  c_1+c_2, \\
4\ell+2  & \mbox{if}  \ \ c_2 < c_1 < c_1+c_2< \lambda. 
\end{array} \right.
\end{equation}
Note that $\tjj_{\ucc,0} \subset \tjj_{\ucc,1}$.

\begin{lemma}\label{representedclasses}
Given $\mu= \ell+\lambda $, where $\ell \in \N$ and $\lambda \in (0,1]$,   $1> c_1 + c_2 > c_1 > c_2 $ and a tamed almost complex structure $J$ belonging to the stratum $\tjj_{\ucc,m}$, a class $D_i$, $i\neq -m$, is represented by some embedded $J$--holomorphic sphere if and only if $D_i \cdot D_{-m}\geq 0$. In particular, the class $D_i$ always has embedded $J$--holomorphic representatives for all $\tjjucc$ if either
\begin{itemize}
\item $\lambda \leq c_2 $ and $i \geq 4\ell-2$, or
\item $c_2 <\lambda \leq c_1$ and $i \geq 4\ell+1$ or $i=4\ell-1$, or
\item $c_1 < \lambda \leq  c_1+c_2$ and $i \geq 4\ell$, or
\item $c_1+c_2< \lambda$ and $i \geq 4\ell +3$ or $i=4\ell+1$,
\end{itemize}
except if $i=4k-2$ or $i=4k$, $k \in \Z$, and the classes $F-E_1-E_2$ or $E_1-E_2$, respectively, are represented by embedded $J$--holomorphic spheres. 
\end{lemma}
\begin{proof}
This proof mimics the proof of Lemma 2.6 in ~\cite{Pi} with some minor changes. As in that paper, it is obvious, by positivity of intersections, that if $J \in \tjj_{\ucc,m}$ and $D_i$ is represented by an embedded $J$--holomorphic curve then we have $D_i \cdot D_{-m}\geq 0$. Now, suppose we have $D_i$ satisfying the latter inequality and $D_i$ does not have any embedded $J$--holomorphic  representative. This immediately implies that $i \geq 0$. By Lemma ~\ref{classDi} the Gromov invariant of the class $D_i$ is nonzero. Therefore there must be a $J$--holomorphic cusp--curve representing $D_i$  and passing through $k=k(D_i)$ generic points $\{p_1, \hdots, p_k \}$. Such $J$--holomorphic cusp--curve $C=\bigcup_{i=1}^Nm_iC_i$, where $N\geq 2$, defines a decomposition of $D_i$ of the form 
\begin{equation}
\begin{split} 
D_i  \ =  \  D_j + \sum_{a_i}s_{a_i}E_i & +\sum_{b_j}t_{b_j}(F-E_j)+\sum_cu_cF+ \\
& + \sum_dw_d(E_1-E_2)+\sum_ez_e(F-E_1-E_2),
\end{split}
\end{equation}
where the multiplicities $s_{a_i},t_{b_j},u_c,w_d,z_e$ are nonnegative and where $D_j \cdot D_{-m}\geq 0$, unless $j=-m$.
Consider the reduced cusp--curve $\bar C$ obtained from $C$ by removing all but one copy of its repeated components. Then the idea is to show that the dimension of the space of unparameterized and reduced cusp--curves of type $\bar C$  is less than $2k(D_i)$, and therefore there is no such cusp--curve passing through $k(D_i)$ generic points.
The last statement in the lemma follows from checking, for  values of the capacities $c_1,c_2$ varying in the four possible ranges, to which class $D_k$ corresponds  the highest codimensional stratum in $\tjj_{\ucc,m}$ and imposing that $D_k \cdot D_i \geq 0$. This gives the values of $i$ for each case. The exception is due to the fact that $D_{4k-2}\cdot (F-E_1-E_2)=D_{4k}\cdot (E_1-E_2)=-1$.
\end{proof}

\section{Homotopy type of symplectomorphisms groups}

\subsection{Connectedness of $\Gucc$}\label{Connexity}

This subsection is devoted to the proof of the connectedness of $\Gucc$ for all values of $\mu, c_{1}, c_{2}$. Since the arguments are similar to the case of the one point blow-up of $S^{2}\times S^{2}$, we only give a sketch of the proof, referring to the literature where needed.

\begin{figure}[htbp]
       \begin{center}
          \psfrag{A}{$c_1$}
          \psfrag{B}{$c_1$}
          \psfrag{C}{$1-c_2$}
          \psfrag{D}{$1$}
          \psfrag{E}{$\mu-c_2$}
          \psfrag{F}{$\mu$}
           \resizebox{!}{4cm}{\includegraphics{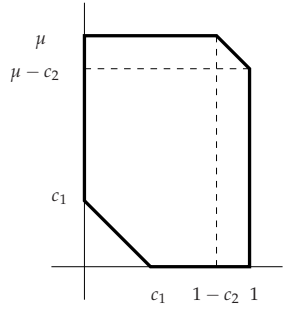}}
           \caption{Moment polytope $\Delta_0$}
           \label{fig:ProductAction}
        \end{center}
    \end{figure}

\begin{lemma}\label{lemma:ExistenceToricAction}
Any symplectic form on $\X_{3}$ admits a Hamiltonian $T^{2}$-action.
\end{lemma}
\begin{proof}
 From Lemma \ref{le:standard} one knows that any symplectic form on $\X_{3}$ is diffeomorphic to a form Poincaré dual to $\mu B+F-c_{1}E_{1}-c_{2}E_{2}$ with $c_{2}\leq c_{1}<c_{1}+c_{2}\leq 1$. Moreover, we can perform the simultaneous blow-up of sizes $c_{1}, c_{2}$ equivariantly with respect to the  product action of $T^{2}=S^{1}\times S^{1}$ on $M_{\mu}^{0}$ whose moment polytope $\Delta_{0}$ is given in Figure~\ref{fig:ProductAction}.
\end{proof}

Let $T^{2}_{0}$ be the toric action on $\tMucc$ described in Lemma~\ref{lemma:ExistenceToricAction}. Let write $\tjj_{\text{good}}\subset\tjj_{\ucc}$ for the subset of almost complex structures for which the exceptional classes $\{E_{1}, E_{2}, B-E_{1}, B-E_{2}, F-E_{1}, F-E_{2}\}$ are simultaneously represented by $J$-holomorphic embedded spheres. By Lemma~A.12 in~\cite{KKP}, $\tjj_{\text{good}}$ is open, dense, and path-connected in $\tjj_{\ucc}$. Let $\mathnormal{Conf}$ be the space of configurations of embedded symplectic spheres $S_{1}\cup S_{2}\cup S_{3}\cup S_{4}$ intersecting transversely and positively, representing the classes $[S_{1}]=E_{1}$, $[S_{2}]=F-E_{1}$, $[S_{3}]=B-E_{2}$, $[S_{4}]=E_{2}$, and for which there exists a $J\in\tjj_{\ucc}$ making them $J$-holomorphic. Let $\Gamma_{0}\in Conf$ be the configuration left invariant by the action $T^{2}_{0}$. There is a fibration
\[\tjj(\Gamma_{0})\to\tjj_{\text{good}}\to Conf\]
where $\tjj(\Gamma_{0})$ is the contractible space of almost complex structures for which $\Gamma_{0}$ is represented by $J$-holomorphic embedded spheres. Hence, the spaces $\tjj_{\text{good}}$ and $Conf$ are weakly homotopy equivalent and, in particular, $Conf$ is connected. Now, standard arguments show that the configuration space $Conf$ is weakly homotopy equivalent to its subspace $Conf^{0}$ of configurations whose spheres intersect orthogonally (see~\cite{LP} Proposition 4.9. A similar point occurs in the proof of ~\cite{MS:J-HolomorphicCurves2} Theorem 9.4.7). The Symplectic Neighborhood Theorem implies that the symplectomorphism group $\Gucc$ acts transitively on $Conf^{0}$ yielding a fibration
\[\Gucc(\Gamma_{0})\to \Gucc \to Conf^{0}\]
The connectedness of $\Gucc$ will follow from the following lemma.
\begin{lemma}
The stabilizer $\Gucc(\Gamma_{0})$ is weakly homotopy equivalent to its subgroup $T^{2}_{0}$.
\end{lemma}
\begin{proof}
Let write $p_{i}=S_{i}\cap S_{i+1}$, $i=1,2,3$. The restriction of a symplectomorphism $\phi\in\Gucc(\Gamma_{0})$ to the sphere $S_{2}$ yields a fibration
\[H_{1}\to\Gucc(\Gamma_{0})\to\Symp(S_{2}, p_{1},p_{2})\simeq S^{1}\]
where $H_{1}$ is the subgroup of symplectomorphisms which restrict to the identity on $S_{2}$, and where $\Symp(S_{2}, p_{1},p_{2})$ is the group of symplectomorphisms of the sphere $S_{2}$ fixing two points.
Similarly, the restriction of the differential $d\phi$ to the normal bundle $\nu(S_{2})$ gives a fibration
\[ H_{2}\to H_{1}\to \mathcal{G}(\nu(S_{2}))\simeq \Map(S_{2}; S^{1})\simeq S^{1}\]
where $\mathcal{G}(\nu(S_{2}))$ is the gauge group of the symplectic bundle $\nu(S_{2})$ and where $H_{2}\subset H_{1}$ is the subgroup of symplectomorphisms whose differentials restrict to the identity on $T_{S_{2}}\tM$. By a Moser-type argument, $H_{2}$ is weakly homotopy equivalent to its subgroup 
\[H_{3}=\{\phi\in H_{2}~:~\phi \text{~is the identity near~} S_{2}\}\]
Restricting elements of $H_{3}$ to $S_{1}$ defines a fibration over a contractible space
\[H_{4}\to H_{3}\to \{\psi\in\Symp(S_{1})~:~\psi=\id\text{~near~}p_{1}\}\simeq \{\id\}\]
and looking at the action of elements of the fiber $H_{4}$ on the normal bundle $\nu(S_{1})$ yields
\[
H_{5}\to H_{4}\to\Map((S_{1},p_{1}); (S_{1},\id))\simeq \{\id\}
\]
with fiber $H_{5}\simeq H_{6}:=\{\phi\in H_{4}~:~\phi=\id\text{~near~} S_{1}\}$. This shows that $H_{2}\simeq H_{6}$. Then restricting in a similar way the elements of $H_{6}$ to $S_{3}$ and then to the normal bundle $\nu(S_{3})$, and finally restricting to $S_{4}$ and $\nu(S_{4})$, we see that the subgroup $H_{2}$ is weakly homotopy equivalent to the subgroup $H_{c}\subset H_{2}$ of symplectomorphisms with compact supports in the complement of $\Gamma_{0}$.

\begin{claim*}
$H_{c}$ is contractible. 
\end{claim*} 
Indeed, we can identify $\tMucc\setminus \Gamma_{0}$ with the preimage of the region ${\Delta}_{0}^{0}$ under the moment map of the $T^{2}_{0}$ action. Therefore, the open set $\tMucc\setminus \Gamma_{0}$ admits a Liouville vector field $L$ lifting the Euler vector field $\xi(p)=-p$ of $\R^{2}$ (and horizontal with respect to the affine structure induced by the toric structure). The vector field $L$ is transverse to arbitrarily small tubular neighborhoods of $\Gamma_{0}$, and its flow defines a retraction of $\tMucc\setminus \Gamma_{0}$ onto a standard $T^{2}_{0}$-invariant ball of capacity $0<\epsilon<1$ whose moment map image is the lower left corner of $\Delta_{0}$. This shows that $\tMucc\setminus \Gamma_{0}$ is symplectically star-shaped. Then, by Gromov's Theorem on the contractibility of compactly supported symplectomorphisms of the standard $4$-ball (see~\cite{MS:SymplecticTopology} Corollary~13.16), we conclude that $H_{c}$ is contractible.
\end{proof}

\subsection{Stability of symplectomorphism groups}\label{sc:stability}

Consider the natural action of the identity component of the diffeomorphism group of $\tMucc$ on the space $\Ssucc$  of all symplectic forms on $\tMucc$ in the cohomology class $[\omega_{\mu,c_1,c_2}]$ that are isotopic to $\omega_{\mu,c_1,c_2}$. By Moser's Theorem the group $\Diff_0(\tMucc)$ acts transitively on $\Ssucc$ via an action that we will write as 
$$\phi \cdot \omega = \phi_*(\omega) = (\phi^{-1})^* \omega,$$
so that $\Ssucc$ is simply the homogeneous space $\Diff_0(\tMucc)/\Gucc $. So this defines a fibration
$$ \Gucc = \Symp(\tMucc) \cap \Diff_0(\tMucc) \to \Diff_0(\tMucc) \to \Ssucc,$$
whose fiber is the symplectomorphism group $\Gucc$.
Therefore, to prove that the groups $\Gucc$ and ${ G}_{\uccp}$ are homotopy equivalent, it is enough to construct a homotopy equivalence between the spaces $\Ssucc$ and $\Ss_{\uccp}$ such that the following diagram 
$$ \xymatrix{\Gucc \ar[r] & \Diff_0 \ar[r] \ar@{=}[d] & \Ssucc \ar@{<->}[d] \\
             {\widetilde G}_{\uccp} \ar[r] & \Diff_0   \ar[r] & \Ss_{\uccp}}$$
commutes up to homotopy. Since we do not know any way of constructing such a map directly, we will use a clever idea introduced by McDuff in ~\cite{MD-Haefliger}. It consists in considering instead a larger space $ \Xxucc$ of pairs
$$ \Xxucc = \{ (\omega, J) \in \Ssucc \times \Aaucc \, : \, \omega \mbox{ tames } J \},$$ 
where $\Aaucc$ denotes the space of almost complex structures that are tamed by some form in $\Ssucc$.\footnote{Recall the $\omega$ is said to tame $J$ if $\omega(v,Jv) >0$ for all $v \neq 0$. }
Then both projection maps $\Xxucc \to \Aaucc$, $\Xxucc \to \Ssucc$ are fibrations with contractible fibers, and so are homotopy equivalences.  
Therefore, to show that the groups $\Gucc$ and ${ G}_{\uccp}$ are homotopy equivalent, it is sufficient to find a homotopy equivalence $\Aaucc \to \Aa_{\uccp}$, that commutes, up to homotopy, with the action of $\Diff_0$. Indeed, we prove a stronger statement, namely:

\begin{prop}\label{prop:stability}
Let $\mu=\ell +\lambda \geq 1$, where $\ell \in \N$ and $\lambda \in (0,1]$. Given $\mu' \in (\ell,\ell+1]$, write $\mu'=\ell+\lambda'$, $\lambda' \in (0,1]$. Consider $c_1,c_2,c_1',c_2'\in (0,1)$ with $c_2<c_1, \ c_1+c_2<1$ and $c_2'<c_1', \ c_1'+c_2'<1$ such that either
\begin{itemize}
\item $\lambda \leq c_2 <  c_1 < c_1+c_2 $ and $\lambda' \leq c_2' <  c_1' < c_1'+c_2' $; or 
\item $c_2 < \lambda \leq c_1 < c_1+c_2 $ and $c_2' < \lambda' \leq c_1' < c_1'+c_2' $; or
\item $c_2 <c_1 < \lambda \leq c_1+c_2$ and $c_2' <c_1' < \lambda' \leq c_1'+c_2'$; or 
\item $c_2 <c_1 < c_1+c_2 < \lambda$ and $c_2' <c_1' < c_1'+c_2' < \lambda'$. 
\end{itemize}
Then the spaces $\Aaucc$ and $\Aa_{\uccp}$ are equal.
\end{prop}
The proof of this proposition is postponed to Appendix ~\ref{se:appendixA}.

\subsection{The stratification of $\Aaucc$}\label{sc:stratification}

\begin{defn}\label{def:strata}
Given $\ell >0$, let $\Aa_{\mu,\ell, r_1,r_2}$ be the subset of $J \in \Aaucc$ consisting of elements that admit a $J$--holomorphic curve in class $A=B-\ell F-r_1E_1-r_2E_2$ where $r_1,r_2 \in \{0,1 \}$. Further define
$$ \Aa_{\mu,0,r_1,r_2}=\Aa_{\ucc}- \bigcup_{\ell >0}\Aa_{\mu,\ell, r_1,r_2},  \quad r_1,r_2 \in \{0,1 \}.$$
\end{defn}   
\noindent From Lemmas ~\ref{lm:classes} and ~\ref{lm:F} it easily follows that  $\Aa_{\mu,0,r_1,r_2}$ is the set of $J$ for which the only classes with $J$--holomorphic representatives are $pB+qF-r_1E_1-r_2E_2$ with $p,q \geq 0$. It is also clear  that the sets $\Aa_{\mu,\ell, r_1,r_2}$ are disjoint for $\ell \geq 0$. Moreover, $\Aa_{\mu,\ell, 0,0}$ is nonempty only if $\mu > \ell$, and  the sets $\Aa_{\mu,\ell, 1,0}$, $\Aa_{\mu,\ell, 0,1}$ and 
$\Aa_{\mu,\ell, 1,1}$ are nonempty only if $\lambda >c_1$, $\lambda >c_2$ and $\lambda >c_1+c_2$, respectively. 

The proof of the following lemma follows exactly the same ideas of Abreu  in ~\cite{Ab} and McDuff in ~\cite{MD-Haefliger} and ~\cite{MD-Stratification}.

\begin{lemma}\label{codimensionAa}
For all $\ell \in \N$ and $\mu \geq 1 > c_1+c_2 > c_1 >c_2$, $\Aa_{\mu,\ell, r_1,r_2}$ is a disjoint union of Fréchét suborbifolds of  $\Aaucc$, all  of codimension $n(\ell,r_1,r_2)=2-2c_1(A)=4\ell -2+2r_1+2r_2$, where $r_1,r_2 \in \{0,1\}$.
\end{lemma}

\begin{thm}\label{thm:stability}
Let $\mu=\ell +\lambda \geq 1$, where $\ell \in \N$ and $\lambda \in (0,1]$. Given $\mu' \in (\ell,\ell+1]$, write $\mu'=\ell+\lambda'$, $\lambda' \in (0,1]$. Consider $c_1,c_2,c_1',c_2'\in (0,1)$ with $c_2<c_1, \ c_1+c_2<1$ and $c_2'<c_1', \ c_1'+c_2'<1$ such that either
\begin{itemize}
\item $\lambda \leq c_2 <  c_1 < c_1+c_2 $ and $\lambda' \leq c_2' <  c_1' < c_1'+c_2' $; or 
\item $c_2 < \lambda \leq c_1 < c_1+c_2 $ and $c_2' < \lambda' \leq c_1' < c_1'+c_2' $; or
\item $c_2 <c_1 < \lambda \leq c_1+c_2$ and $c_2' <c_1' < \lambda' \leq c_1'+c_2'$; or 
\item $c_2 <c_1 < c_1+c_2 < \lambda$ and $c_2' <c_1' < c_1'+c_2' < \lambda'$. 
\end{itemize}
Then the symplectomorphism groups $\Gucc$ and $G_{\uccp}$ are homotopy equivalent.
Moreover, for fixed values of $c_1$ and $c_2$, as $\lambda$ passes $c_2$, $c_1$ or $c_1+c_2$, or $\mu$ passes the integer $\ell+1$, $\ell \geq 0$, the groups $\pi_i(\Gucc)$, with $i\leq 4\ell -3$ in the first two cases and $i\leq 4\ell -1$ in the last two, do not change. 
\end{thm}
\begin{proof}
The first statement is an immediate consequence of Proposition \ref{prop:stability}.  From Lemma \ref{codimensionAa} one knows that when $\mu$ passes $\ell+c_2$ the topology of $\Aaucc$ changes by the addition of strata of codimension $n(\ell)=4\ell$. Hence its homotopy groups $\pi_i$  for $i\leq 4\ell-2$ do not change. It follows, using the long exact homotopy sequence of the fibration 
$$ \Gucc \longrightarrow \Diff_0(\X_3) \longrightarrow \Ssucc \simeq \Aaucc,$$
that the groups $\pi_i(\Gucc)$ for $i \leq 4\ell-3$ do not change. The argument to prove the other cases is similar and this concludes the proof of the second statement. 
\end{proof}

\subsection{Homotopy type}
\label{HomotopyType}

\noindent{\it Proof of Proposition \ref{stabilizer}.}
This is Proposition 4.21 of Lalonde and Pinsonnault in \cite{LP}. Let us recall their proof. The main idea is to prove that the stabilizer  $\Symp_p(\tMuc)$ is homotopy equivalent to the group $\Symp^{U(2)}(\tMucc,E_2)$ of the simplectomorphisms in $\Symp(\tMuc)$ that act linearly in a small neighborhood of the exceptional fiber whose homology class is $E_2$. Every symplectomorphism in $\Symp^{U(2)}(\tMucc,E_2)$ gives rise to a symplectomorphism of $\tMuc$ acting linearly near the embedded ball $B_{c_2}$ and fixing its center $p$. Conversely, every homotopy class of the stabilizer $\Symp_p(\tMuc)$ can be realized by a family of simplectomorphisms that act linearly on a ball $B_{c_2'}$ of sufficiently small capacity $c_2'$ centered at $p$. Hence we can lift such a representative to the group $\Symp^{U(2)}(\tM_{\mu, c_1,c_2'},E_2)$.  If $r$ is the restriction map, then the directed system of homotopy maps
$$ r_{c_2', c_2}: \Symp^{U(2)}(\tM_{\mu, c_1,c_2'},E_2) \to \Symp^{U(2)}(\tMucc,E_2)$$
stabilizes if $c_2 \leq c_2'$ are sufficiently small, that is, they are homotopy equivalences, by the Stability Theorem \ref{thm:stability}.  This system together with the maps 
$$ g_{c_2}: \Symp^{U(2)}(\tMucc,E_2) \to  \Symp_p(\tMuc)$$
yields a  commutative triangle for each pair $c_2 \leq c_2'$. It is then obvious that each map $g_{c_2}$ is a weak homotopy equivalence. On the other hand one knows, by  Lemma 2.3 in \cite{LP}, that  $\Symp^{U(2)}(\tMucc,E_2)$ is homotopy equivalent to $\Symp(\tMucc,E_2)$. By Lemma  \ref{lm:E2} we know that the exceptional curve $E_2$ cannot degenerate, that is, for every $J$ in the space $\tjjucc $ of $\omega_\mu$--tamed almost complex structures, there is a unique  embedded $J$--holomorphic representative of $E_2$ and no cusp--curve in class $E_2$. Therefore, by Lemma 2.4 in  ~\cite{LP} we can conclude that $\Gucc$ retracts onto its subgroup $\Symp(\tMucc,E_2)$ and this concludes the proof. \qed

Note that Proposition \ref{Inflation} follows from the results in Sections \ref{sc:stability} and \ref{sc:stratification}.  The proof is an analog of Theorem 2.2 in ~\cite{ALP}.

\section{The Lie algebra  $\pi_*(\Gucc)\otimes \Q$}\label{se:LieAlgebra}

In this section we show that the rational homotopy Lie algebra $\pi_*(\Gucc)\otimes \Q$  is finitely generated  and  we describe geometrically its generators. Moreover,  we show that they  are represented by  Hamiltonian $S^1$--actions on $\tMucc^i$. 
\subsection{Homotopy groups of $\Gucc$} 

We begin by reviewing some facts about rational homotopy theory. First, the rational dichotomy discovered by Félix  \cite{FHT} states that if $X$ is an $n$-dimensional simply connected space with the rational cohomology of finite type and finite category, then either its rational homotopy is finite dimensional or the dimensions of the rational homotopy groups grow exponentially. In the first case, the space is called {\it rationally elliptic} and in the second {\it rationally hyperbolic}. 

\begin{prop}[see Part VI in ~\cite{FHT}]
If $X$ is rationally elliptic, then
\begin{enumerate}
\item $\dim \pi_{even}(X) \otimes \Q \leq \dim \pi_{odd}(X) \otimes \Q \leq \, cat(X),$
\item $\chi(X) \geq 0$, where $\chi$ denotes the Euler characteristic.  
\end{enumerate}
\end{prop}

If $X$ is a 4--dimensional simply connected finite CW--complex, then it follows from the basic properties of the Lusternik--Schinrelmann category that $cat(X) \leq 2$.   The above proposition implies that if $\pi_2(X) \otimes \Q \geq 3$ then $X$ is rationally hyperbolic. Thus every simply connected 4--manifold $X$ whose  $b^2(X):= \dim H^2(X; \Q) \geq 3 $ is rationally hyperbolic. We conclude that if $n \geq 2$ then ${\CP^2\#\,n\,\overline{\CP}\,\!^2}$, and in particular $\tMuc^i$ (since it is diffeomorphic to ${\CP^2\#\,2\,\overline{\CP}\,\!^2}$ as a smooth manifold), is rationally hyperbolic.

Recall also that for any topological space $X$ with rational homology of finite type, the {\it Poincaré series} for $X$ is the formal power series \begin{displaymath}P_X= \sum_0^\infty \dim H_n(X; \Q)z^n.\end{displaymath}

Denote by $r_n$ the integers $r_n= \dim \pi_n(\Omega X)\otimes \Q = \dim \pi_{n+1}(X) \otimes \Q$. Then, as explained in \cite{FHT} \S 33, the Poincaré--Birkoff--Witt Theorem identifies $P_{\Omega X}$ as the formal power series

\begin{equation}\label{powerseries}
P_{\Omega X} = \frac{\prod_{n=0}^\infty (1+z^{2n+1})^{r_{2n+1}}}{\prod_{n=1}^\infty (1-z^{2n})^{r_{2n}}}.
\end{equation}

Note that \eqref{powerseries} provides algorithms for computing the integers $\dim \pi_{n+1}(X) \otimes \Q$, $1 \leq n \leq N $, from the integers $\dim H_n(\Omega X;\Q )$, $1 \leq n \leq N $. The latter can be computed using the Serre spectral sequence of the path fibration $\Omega X \to PX \to X$, where $PX$ is the space of paths in $X$ starting in $x_0$, and $p: PX \to X$ sends each path to its endpoint. For $X=\X_{2}:={\CP^2\#\,2\,\overline{\CP}\,\!^2}$ the spectral sequence gives $H_0(\Omega \X_{2};\Q)= \Q$, $H_1(\Omega \X_{2};\Q)= \Q^3$ and if we write $h_n= \dim H_n(\Omega \X_{2};\Q)$ then it is easy to check that $h_n=3h_{n-1}-h_{n-2}$ when $n \geq 2$.  Therefore, using \eqref{powerseries}, we can compute the dimension of 
the homotopy groups of $\X_{2}$ and we obtain, for example, $r_1=3$, $r_2=r_3=5$, $r_4=10$, $r_5=24$ and $r_6=352$.

\begin{prop}\label{homotopygroups}
Let $\ell < \mu \leq \ell+1$ with $\ell \geq 1$ or $\mu=1$. Recall that $\mu= \ell + \lambda$ and $0< c_2 < c_1 < c_2+c_1< 1$. Then the rational homotopy groups of $\Gucc$ are given by 
\begin{align*}
(1) \quad \pi_n & =  \left\{ \begin{array}{ll}
                   \Q^5 & \mbox{if $n=1$} \\
                   \Q^{r_n} & \mbox{if $n \geq 2$},
                   \end{array} \right.  \quad   \mbox{when $ c_1+c_2 < \mu=1$};\\
(2) \quad \pi_n & = \left\{ \begin{array}{ll}
                   \Q^6 & \mbox{if $n=1$} \\
                   \Q^{r_{4\ell-2}+1} & \mbox{if $n=4\ell-2$} \\
                   \Q^{r_n} & \mbox{if $n \neq 1,4\ell-2$},
                   \end{array} \right. \quad \mbox{ when $\mu > 1$ and  $ c_2 \leq \lambda \leq c_1 $;}\\
 \end{align*}
 \begin{align*}                                    
(3) \quad \pi_n & = \left\{ \begin{array}{ll}
                   \Q^6 & \mbox{if $n=1$} \\
                   \Q^{r_{4\ell}+1} & \mbox{if $n=4\ell$} \\
                   \Q^{r_n} & \mbox{if $n \neq 1,4\ell$},
                   \end{array} \right.\quad \mbox{ when $\mu >1$ and  $c_1+c_2 < \lambda $.} 
\end{align*}
\end{prop}
\begin{proof}
We consider the {\it evaluation fibration} 
$$\Symp_p(\tMuc) \to \Symp(\tMuc) \to \tMuc $$
and the Stability Theorem \ref{thm:stability}. Then the rational homotopy groups of $\Gucc$ are easily computed by combining Proposition \ref{stabilizer} and the computation of the homotopy groups of the symplectomorphism group $ \Symp(\tMuc)$:
\begin{prop}[Pinsonnault \cite{Pi}, Corollary 3.10] The following holds. 
\begin{enumerate}
\item when $\mu=1$, the only non--trivial rational homotopy group of the topological group $\Symp(\tMuc)$ is $\pi_1\simeq \Q^2$.  
\item When $ \mu>1$, the non--trivial rational homotopy groups of the topological group $\Symp(\tMuc)$ are $\pi_1\simeq \Q^3$ and $\pi_k \simeq \Q$ where
$$ k= \left\{ \begin{array}{ll}
4\ell -2  &  \mbox{ if } \lambda \leq c_1; \\
4\ell  &  \mbox{ if } \lambda > c_1.
\end{array} \right. $$
 \end{enumerate}
\end{prop}
If $\mu =1$ the above homotopy fibration gives an exact sequence which immediately yields $\pi_n=\Q^{r_n}$ if $n >1$, and 
the nontrivial term 
$$0 \to \Q^3 \to \pi_1 \to \Q^2 \to 0$$ gives $\pi_1 = \Q^5$. 
In cases (2) and (3) we use  again the homotopy fibration together with the Stability Theorem \ref{thm:stability}, Proposition \ref{stabilizer} and Corollary 3.10 in ~\cite{Pi}.
\end{proof}

\subsection{Conventions}\label{conventions} \hfill

We now  want to compute the graded Lie algebra $ \pi_*(\Gucc)\otimes \Q$. In order to carry this  computation we need to describe the generators of this algebra in more detail and for that we follow the conventions  of Anjos, Lalonde and Pinsonnault in ~\cite{ALP}:
\begin{enumerate}
\item Let $T^{4}\subset\U(4)$ act in the standard way on $\C^{4}$. Given an integer $n\geq 0$, the action of the subtorus $T^{2}_{n}:=(ns+t,t,s,s)$ is Hamiltonian with moment map
\[(z_{1},\ldots,z_{4})\mapsto(n|z_{1}|^{2}+|z_{3}|^{2}+|z_{4}|^{2},|z_{1}|^{2}+|z_{2}|^{2}).\]
We identify $M_{\mu}^{0}=(S^{2}\times S^{2},\mu\sigma\oplus\sigma)$ with each of the toric Hirzebruch surface $\F_{2k}^{\mu}$, $0\leq k\leq \ell$ (where as usual $\ell < \mu \le \ell+1$), defined as the symplectic quotient $\C^{4}/\!/T^{2}_{2k}$ at the regular value $(\mu+k,1)$ endowed with the residual action of the torus $T(2k):=(0,u,v,0)\subset T^{4}$. The image $\Delta(2k)$ of the moment map $\phi_{2k}$ is the convex hull of
\[
\{(0,0),(1,0),(1,\mu+k),(0,\mu-k)\}.
\]
Similarly, we identify $\Mul=(\NTB,\om_{\mu})$ with the toric Hirzebruch surface $\F_{2k-1}^{\mu}$, $1\leq k\leq \ell$, defined as the symplectic quotient $\C^{4}/\!/T^{2}_{2k-1}$ at the value $(\mu+k,1)$. The moment polygon $\Delta(2k-1):=\phi_{2k-1}(\F_{2k-1})$ of the residual action of the $2$-torus $(0,u,v,0)$ is the convex hull of 
\[
\{(0,0),(1,0),(1,\mu+k),(0,\mu-k+1)\}.
\]
Note that the group $\Symp_{h}(M_{\mu})$ of symplectomorphisms acting trivially on homology being connected, any two identifications of $\F_{n}^{\mu}$ with $M_{\mu}^{i}$ are isotopic and lead to isotopic identifications of $\Symp_{h}(\F_{n}^{\mu})$ with $\Symp(M_{\mu}^{i})$. 

\item We identify the symplectic blow-up $\tMuco$ with the equivariant blow-up of the Hirzebruch surfaces $\F_{n}^{\mu}$ for appropriate parameter $\mu$ and capacity $c_1\in(0,1)$.
\begin{enumerate}
\item We define the even torus action $\tT(2k)$ as the equivariant blow-up of the toric action of $T(2k)$ on $\F_{2k}^{\mu}$ at the fixed point $(0,0)$ with capacity $c_1$. The image $\tilde{\Delta}(2k)$ of the moment map $\tilde{\phi}_{2k}$ is the convex hull of
\[
\{(1,\mu+k),(0,\mu-k), (0,c_1), (c_1,0), (1,0)\}.
\]

\item The odd torus action $\tT(2k-1)$, $k\geq 1$, is obtained by blowing up the toric action of $T(2k-1)$ on $\F_{2k-1}^{\mu-c_1}$ at the fixed point $(0,0)$ with capacity $1-c_1$.  Its moment polygon $\tDelta(2k-1)$ is the convex hull of 
\[
\{(1,\mu-c_1+k),(0,\mu-c_1-k+1), (0,1-c_1),(1-c_1,0),(1,0)\}.
\]
\end{enumerate}
 Note that when $c_1<\ccrit := \mu - \ell$, $\tMuco$ admits exactly $2\ell+1$ inequivalent toric structures $\tT(0),\ldots, \tT(2\ell)$, while when $c_1 \ge \ccrit$, it admits only $2\ell$ of those, namely $\tT(0),\ldots, \tT(2\ell-1)$. 

\item The K\"ahler isometry group of $\F_{n}^{\mu}$ is $N(T^{2}_{n})/T^{2}_{n}$ where $N(T^{2}_{n})$ is the normalizer of $T^{2}_{n}$ in $\U(4)$. There is a natural isomorphism\footnote{In the untwisted case, we assume $\mu>1$ so that the permutation of the two $S^{2}$ factors is not an isometry.} $N(T_{0})/T^{2}_{0}\simeq \SO(3)\times \SO(3):=K(0)$, while for $k\geq 1$, we have $N(T^{2}_{2k})/T_{2k}\simeq S^{1}\times \SO(3):=K(2k)$ and $N(T^{2}_{2k-1})/T^{2}_{2k-1}\simeq \U(2):=K(2k-1)$. The restrictions of these isomorphisms to the maximal tori are given in coordinates by
\begin{align*}
(u,v)&\mapsto (-u,v)\in T(0):=S^{1}\times S^{1}\subset K(0)\\
(u,v)&\mapsto (u,ku+v)\in T(2k):=S^{1}\times S^{1}\subset K(2k)\\
(u,v)&\mapsto (u+v,ku+(k-1)v)\in T(2k-1):=S^{1}\times S^{1}\subset K(2k-1).
\end{align*}
These identifications imply that the moment polygon associated to the maximal tori $T(n)=S^{1}\times S^{1}\subset K(n)$  and $\tT(n)$ are  the images of $\Delta(n)$  and  $\tilde{\Delta}(n)$, respectively, under the transformation $C_{n}\in\GL(2,\Z)$ given by
\[
\begin{matrix}
 C_{0}=\begin{pmatrix}-1 & 0 \\ 0 & 1   \end{pmatrix}
&~C_{2k}=\begin{pmatrix} 1 & 0 \\ -k & 1   \end{pmatrix}
&~C_{2k-1}=\begin{pmatrix} 1-k & 1 \\ k & -1 \end{pmatrix}
\end{matrix}
\]
Our choices imply that under the blow-down map, $\tT(n)$ is sent to the maximal torus of $K(n)$, for all $n\geq 0$. Again, because $\Symp(\tMuco)$ is connected (see~\cite{LP,Pi}), all choices involved in these identifications give the same maps up to homotopy.

\item We identify the symplectic blow--up $\tMucc$ with the equivariant two blow--up of the Hirzebruch surfaces $\F_n^\mu$ for appropriate parameter $\mu$ and capacities $c_2<c_1< c_1+c_2 <1$.
In this case, depending on which point we blow--up, we obtain inequivalent toric structures. 
We define the torus actions $T_i(2k),T_i(2k-1) $, $i=1,\hdots,5$ as the equivariant blow--ups of the toric action of $\tT(n)$ on ${\widetilde{\F}}_{2k}^\mu$ and ${\widetilde{\F}}_{2k-1}^{\mu-c_1}$ respectively, with capacity $c_2$, at each one of the five fixed points, which correspond to the vertices of the moment polygon $\tDelta(n)$,  described in item 2(a) and (b).

\item The cohomology ring of $BT_i(n)$ is isomorphic to $\Q[x_{n,i},y_{n,i}]$ where $|x_{n,i}|=|y_{n,i}|=2$. We identify the generators $x_{n,i}$, $y_{n,i}$ with the cohomology classes induced by the circle actions whose moment maps are, respectively, the first and the second component of the moment map associated to $T_i(n)$.  Note that since we work only with topological groups up to rational equivalences, we will also denote by $\{x_{n,i}, y_{n,i}\}$ the generators in $\pi_{1}(T_i(n))$ and in $\pi_{2}(BT_i(n))$.
\end{enumerate}

\begin{remark}
Each stratum $\tjjucc$, whose configuration is given in Figures 2 to 6 corresponds to a toric structure on $\tMucc$, unique up to equivariant symplectomorphisms (see Lemma \ref{le:isometries} in Appendix B). These toric structures are given by the tori $T_i(n)$, $n \in \N$. More precisely, the relationship between the tori $T_i(n)$  and various configurations is the following: 
\begin{itemize}
\item $T_i(0), \ i=1, \hdots, 5$ correspond to configurations (1), (3), (5), (4), (2), in Figure 2, respectively; 
\item $T_i(2k-1),\  i=1, \hdots, 5$ and $ k \geq 1$  correspond to configurations (10), (12), (11), (9), (7), in Figures 3 and 4, respectively; 
\item $T_i(2k),\  i=1, \hdots, 5$ and $ k \geq 1$  correspond to configurations (16), (18), (17), (15), (13), in Figures 5 and 6, respectively.
\end{itemize}
\end{remark}

McDuff proved in \cite{MT}  that the maps $T_i(n)\to \Gucc$ induce injective maps of fundamental groups. More precisely, this follows from  Theorems 1.3 and 1.24 in that paper, as explained in Remark 1.25, since those two theorems imply that if  $(M, \omega, T, \Phi)$ is a symplectic toric 4--manifold  with moment polytope $\Delta$ which has five or more facets then the map $\pi_1(T) \to \pi_1(\Symp_0(M, \omega))$ is an injection. 

Let us define 
$$x_0:=x_{0,1}\quad  y_0:=y_{0,1}\quad  x_1:=x_{0,2}, \quad y_1:=y_{0,4}-x_{0,4} \quad \mbox{and} \quad z:=y_{0,4}.$$

Recall that, in  \cite{Ka}, Karshon associates a labeled graph (cf. Sections 2.1 and 2.2) to each compact four dimensional Hamiltonian $S^1$--space and then using these graphs she gives a classification (cf. Theorems 5.1 and 4.1) of such spaces with isolated fixed points. Moreover, she shows that a complete list of such spaces  (up to isomorphism) is provided by the list of Delzant polytopes, and two Delzant polytopes give isomorphic spaces if and only if their graphs, as constructed in Section 2.2, coincide. 
Therefore, constructing these graphs for the $S^1$--actions contained in the tori $T_i(n)$ one can obtain the following relations between them. 

\begin{lemma}\label{classification}
Let  $\mu \geq 1 > c_1+c_2 > c_1 >c_2$.   If  $k=0$ then the following identifications hold in $\pi_1(\Gucc)\otimes \Q$
$$x_{0,5}=x_0, \quad y_{0,5}=z,  \quad x_{0,3}=x_1, \quad y_{0,3}=x_1+y_1 \quad  \mbox{and} \quad y_{0,2}=y_0.$$
Moreover, for  all admissible values $k,j \geq 1$, we have the following identifications
\begin{gather*}\label{identifications}
x_{1,4}= x_{1,5}=y_1, \quad  \quad  \quad   x_{1,1}=y_0-x_0,\quad  \quad  \quad x_{1,2}=y_0-x_1,\\
(k-1)x_{2j-1,i}+ky_{2j-1,i}=(j-1)x_{2k-1,i}+jy_{2k-1,i}  \quad \mbox{where} \quad i=3,4,5 \\
(k-1)x_{2j-1,1}+ky_{2j-1,1}=(j-1)x_{2k-1,2}+jy_{2k-1,2},\\
kx_{2k,i}+y_{2k,i}=(k+1)x_{2k-1,i}+ky_{2k-1,i}  \quad \mbox{where} \quad i=1,2,4,\\
kx_{2k,3}+y_{2k,3}=(k+1)x_{2k-1,5}+ky_{2k-1,5}, \quad 
kx_{2k,5}+y_{2k,5}=(k+1)x_{2k-1,3}+ky_{2k-1,3}, \\
x_{2k,1}=x_{2k,5}, \quad  \quad  \quad x_{2k,2}=x_{2k,3},\\
y_{2k,1}=kx_1+y_0, \quad\quad y_{2k}^2=kx_0+y_0, \quad\quad y_{2k,3}=(k+1)x_1+y_1, \quad\quad y_{2k,5}=kx_0+z, \\
jx_{2k,i}-y_{2k,i}=kx_{2j,i}-y_{2j,i} \quad \mbox{where} \quad i=3,4,5\\
jx_{2k,1}-y_{2k,1}=kx_{2j,2}-y_{2j,2},\\
x_{2k-1,1}+y_{2k-1,1}=x_{2k-1,5}+y_{2k-1,5},  \quad
x_{2k-1,2}+y_{2k-1,2}=x_{2k-1,3}+y_{2k-1,3},\\
kx_{2k,4}+y_{2k,4}=kx_{2k,5}+y_{2k,5}.
\end{gather*}
\end{lemma}

\begin{remark}\label{rmk:identifications}
Let us define $t:=x_{1,1}+y_{1,1}$. Then, from the previous lemma, it is not hard to check that the images of the generators of $\pi_{1}(T_i(n))$ in $\pi_1(\Gucc)\otimes \Q$ are determined by the image of the generators $x_0,y_0,x_1,y_1,z \in \pi_{1}(T_i(0))$ where $i=1,2,4$  and $t \in \pi_1(T_1(1))$. Indeed, for $k \geq 1$, we have 
\begin{align*}
& x_{2k-1,1}=y_0-kx_0,  & & y_{2k-1,1} =k(t-x_1)+x_1+kx_0-y_0,\\
& x_{2k-1,2}=y_0-kx_1,  & & y_{2k-1,2} =kt+x_0-y_0,\\
& x_{2k-1,3}=z-kx_0,  & & y_{2k-1,3}  =k(t-x_1)+(k+1)x_0-z,\\
& x_{2k-1,4}=ky_1+(1-k)z, & &y_{2k-1,4} =k(t-x_1+z)-(k+1)y_1,\\
& x_{2k-1,5}=y_1+(1-k)x_1,  & & y_{2k-1,5} =kt-y_1,\\
& & & \\
& x_{2k,1} = k(t-x_1)-x_0, &  & y_{2k,1} = kx_1+y_0,\\
& x_{2k,2} =k(t-x_1)-x_1 , & & y_{2k,2}  = kx_0+y_0,\\
& x_{2k,3} =k(t-x_1)-x_1 , & & y_{2k,3}  = (k+1)x_1+y_1, \\
& x_{2k,4} = k(t-x_1)+y_1-z, & &y_{2k,4}  = (k+1)z-ky_1,\\
& x_{2k,5} = k(t-x_1)-x_0, & &y_{2k,5}  = kx_0+z,\\
\end{align*}

\end{remark}

Besides the circle actions coming from these tori $T_i(n), i=1, \hdots,5$, there are circle actions obtained by blowing up the toric actions $\tT(n) $ at an interior point of a curve representing the exceptional class $E_1$, where $c_1 > c_2$, such as in configurations  (6), (8) and (14) in Figures \ref{mu=1}, \ref{configB-MF} and \ref{configB-MF-E1}. We will denote the corresponding generators in homotopy by $a_0, a_{2k-1}$ and $a_{2k}$, respectively. Using again the classification of Hamiltonian $S^1$--actions, obtained by Karshon in ~\cite{Ka}, it is easy to compute the relations satisfied by these generators. 

\begin{lemma}\label{s1action}
If  $\mu \geq 1 > c_1+c_2 > c_1 >c_2$ then the following identifications hold in $\pi_1(\Gucc) \otimes \Q$. 
$$a_0=y_1, \ \ a_{2k}=(k+1)x_{2k,3}+y_{2k,3}=(k+1)x_{2k,4}+y_{2k,4}\mbox{ \ \ and \ } $$
$$ a_{2k-1}=(k+1)x_{2k-1,3}+ky_{2k-1,3}=(k+1)x_{2k-1,4}+ky_{2k-1,4},$$
where $k \geq 1$. Remark \ref{rmk:identifications} then implies that
$ a_{2k-1}=k^2(t-x_1) +z \mbox{ \ and \ } a_{2k}=k(k+1)(t-x_1)+y_1.$
\end{lemma}

\noindent We will see next, in the proof of Theorem \ref{maintheorem}, that the generators  $x_0,y_0,x_1,y_1,z,t$ have a fundamental role in the description of the homotopy Lie algebra $\pi_*(\Gucc)\otimes \Q$.
\medskip

\subsection{Homotopy Decomposition}
\label{HomotopyDecomposition}
All of arguments of \cite{AGK} on $M_\mu$  and of \cite{ALP} on $\tMuc$, regarding the homotopy decomposition of the classifying space of the group $\Gucc$ apply as well for the group of symplectomorphisms of the blow--up $\tMucc$ although there are some differences in this case. Depending on the size of  $\lambda$, as we saw in Section  \ref{StructureJholomorphic}, now one may add two or four distinct strata to the space of almost complex structures at once when $\mu$  passes a critical value instead of just one stratum each time as in the previous cases. Moreover, not all strata correspond to a toric structure on the manifold as before. Now there also are some strata which correspond just to a circle action on $\tMucc$, unique up to equivariant symplectomorphisms. More precisely, the strata represented in Section \ref{StructureJholomorphic} by configurations (6), (8) and (14) are the ones that correspond to $S^1$--actions on $\tMucc$ and which were classified in the previous section.  

The next proposition is a particular case of Theorem \ref{thm:HomotopyPushout} since it gives the homotopy decomposition in the case $c_1 < \lambda \leq c_1+c_2$. The homotopy pushout diagram should be obvious for the other cases. Recall that given a $G$-space $X$, we write $X_{hG}$ for the homotopy orbit $X_{hG}:=EG\times_{G} X$.

\begin{prop}\label{SpecialCaseHomotopyPushout}
If $\ell < \mu \leq \ell+1 $ where $\ell \in \N$, $c_2 < c_1 < c_2+c_1<1$ and $c_1 < \lambda \leq c_1+c_2$, there is a homotopy pushout diagram
$$
\xymatrix{
S^{4\ell-1}_{h{T_1}(2\ell)} \sqcup S^{4\ell+3}_{h{T_4}(2\ell)} \sqcup S^{4\ell+1}_{h{T_5}(2\ell)}\sqcup S^{4\ell+1}_{h{S^1}(2\ell)}\ar[r] \ar[d]^{j_{2\ell}} & B{T_1}(2\ell) \sqcup B{T_4}(2\ell)  \sqcup B{T_5}(2\ell)\sqcup BS^1(2\ell)\ar[d]^{i_{2\ell}} \\
B{G}_{\mu',c_1,c_2} \ar[r] & B\Gucc}
$$
where $\mu'= \ell + \lambda'$ satisfies $ c_2 <  \lambda'  \leq c_1$,
so that $G_{\mu',c_1,c_2}$ is the group associated with a stratification having four strata less than the stratification associated with $G_{\mu,c_1,c_2}$. The upper horizontal map  is the universal bundle map associated to the representations of ${T_i}(2\ell)$ and ${S^1}(2\ell)$  on $H^{0,1}_{J_{i,2\ell}}(T\tMucc)$\footnote{$H^{0,1}_{J_{i,2\ell}}(T\tMucc)$ denotes the moduli space of infinitesimal deformations of the complex structure $J_{i,2\ell}$.} and on $H^{0,1}_{J_{2\ell}}(T\tMucc)$ respectively, $i_{2\ell}$ is induced by the inclusions ${T_i}(2\ell)\into \Gucc$ and ${S^1}(2\ell)\into \Gucc$, and the map $BG_{\mu',c_1,c_2}\to B\Gucc $ coincides, up to homotopy, with the one given by inflation described in Proposition ~\ref{Inflation}.
\end{prop}

The proof of this proposition follows the same steps of the proof of Theorem 2.5 in ~\cite{ALP} . However, since in our case the structure of the space of almost complex structures is more complicated, the analogue of Appendix A in  ~\cite{ALP} requires a more careful explanation which is given in Appendix B. Nevertheless a similar argument applies to prove that the action of the symplectomorphism group on the space of compatible almost complex structures is homotopically equivalent to its restriction to the subset of compatible integrable complex structures, which is one of the key steps in the proof. Another key step is the replacement of the action of $\Gucc$ on complex structures by an equivalent $A_{\infty}$ homotopy action, as explained in Appendix C of~\cite{AGK}. This second step is necessary in order to deal with two different issues: the fact that tubular neighborhoods of strata of equal codimensions may intersect nontrivially, and the nonexistence of $\Gucc$-invariant tubular neighborhoods of the strata. In particular, this justifies our use of disjoint unions in the upper-left corner of the pushout diagram of Proposition~\ref{SpecialCaseHomotopyPushout}, in spite of the fact the normal neighborhoods may intersect. Again, this clever argument was first introduced by Abreu--Granja--Kitchloo in~\cite{AGK}.

\subsection{Proof of Theorem \ref{maintheorem}}\label{ProofMainTheorem}

We divide the proof of the theorem into three parts: in the first part we consider the case $\mu=1$, in the second part we give the proof of statements (S2) and (S4)  and the third part is devoted to the proof of statements (S1) and (S3).  In the case $\mu=1$ and statements (S2) and (S4) the symplectomorphism group has the homotopy type of a stabilizer so the argument is similar in all three cases, while in the first and third cases we will need a different technique that involves the computation of the rational cohomology ring of the classifying space of the group $\Gucc$.

\subsubsection{Proof of the case $\mu=1$}
Note that in this case there are six distinct configurations of symplectic embedded surfaces in $\tMucc$ which intersect transversally and positively, namely the ones described in Figure  \ref{mu=1}.  Therefore there are five torus actions, namely $T_i(0)$ with $i=1,\ldots,5$ and one $S^1$--action that corresponds to configuration (6). 
Denote by $\Lambda$ the algebra described in the first statement of the theorem, that is, the Lie graded algebra generated by $x_0,y_0,x_1,y_1,z$ such that all Samelson products between the generators vanish except the five ones indicated there. Let $\lambda_n$ be the rank of $\Lambda$ in dimension $n$. It is clear that the Samelson products $[x_0,y_0], \ [x_1,y_0], \ [x_0,z]$ vanish, since these classes have commuting representatives in $T_i(0)$ when $i=1,2,5$, respectively. The relation $[x_1,y_1]=0$ follows from the third and fourth  identifications in Lemma  \ref{identifications} together with the fact that the classes $x_{0,3},y_{0,3}$ obviously have commuting representatives in $T_3(0)$. Finally, the relation $[z,y_1]=0$ follows from the definition of the classes $z$ and $y_1$ together with the fact that the classes $x_{0,4},y_{0,4}$  have commuting representatives in $T_4(0)$. 
Recall, from Proposition \ref{homotopygroups}, that $\pi_n(\Gucc)\otimes \Q = \Q^{r_n}$ if $n \geq 2$ where $r_n=\dim \pi_{n}(\Omega(\PbP))$.  We need to show that the graded Lie algebras $\Lambda$ and  $ \pi_*(\Gucc)\otimes \Q$ are isomorphic. First we will show that they have the same rank in each dimension, i.e.,  $\lambda_n=r_n$ for all $n \geq 2$ and $\lambda_1=5$ (cf. Proposition \ref{homotopygroups}). It is clear that the  five generators of $\Lambda$ give $\lambda_1=5$. 

Now we need to recall some facts about Lie graded algebras.
Let $L$ be such an algebra and $TL$ the tensor algebra on the graded vector space L. Let $I$ be the ideal in the (associative) graded algebra $TL$ generated by the elements of the form $x \otimes y - (-1)^{\deg x \deg y} y \otimes x - [x,y]$, $x, y \in L$, where $[\cdot,\cdot]$ denotes the Lie bracket in $L$. The graded algebra $TL/I$ is called the {\it universal enveloping algebra} of $L$ and is denoted by $\Uu L$.  We need to use the following important result.

\begin{thm}[Milnor--Moore]
If $X$ is a simply connected topological space then 
\begin{enumerate}
\item $\pi_*(\Omega X)\otimes \Q$ is a graded Lie algebra, $L_X$;
\item The Hurewicz homomorphism extends to an isomorphism of graded algebras 
$$ \Uu L_X \stackrel{\cong}{\longrightarrow}  H_*(\Omega X; \Q).$$
\end{enumerate}
\end{thm}

\medskip

Going back to our proof,  the  Poincaré--Birkoff--Witt Theorem (see \eqref{powerseries}) implies  that the numbers  $\lambda_n$ defined above satisfy the equation

\begin{equation}\label{PBWequation}
\sum_{n=0}^{+\infty}h'_nz^n= \frac{\prod_{n=0}^\infty (1+z^{2n+1})^{\lambda_{2n+1}}}{\prod_{n=1}^\infty (1-z^{2n})^{\lambda_{2n}}},
\end{equation}
where the left hand side of the equation is the Poincaré series of the enveloping algebra $\mathcal U \Lambda$.

\begin{claim}
$h_n'=5h_{n-1}$ if $n \geq 1$ and $h_0'=1$ where $h_n=\dim H_n(\Omega({\CP^2\#\,2\,\overline{\CP}\,\!^2}); \Q)$.
\end{claim}
\begin{proof} The proof is done  by induction on $n$. Note that indeed we have $h'_0=1$, $h'_1=5h_0=5$, which corresponds to the 5 generators $x_0,y_0,x_1,y_1,z$, and $h'_2=5h_1=15$, which corresponds to the classes $x_0y_0, x_0x_1,x_0y_1,x_0z,y_0x_1,y_0y_1,y_0z,x_1y_1,x_1z,y_1z,x_1x_0,y_1x_0,y_1y_0,zy_0$ and $zx_1$.  Moreover to compute $h'_{n+1}$ with $n \geq 2$ we begin with $5h'_n$ elements (the number of classes in dimension $n$ times the 5 generators, where we can assume that we always multiply on the left), then subtract $h'_n$ to that number because the square of each generator vanishes, and finally subtract $h'_n+h'_{n-1}$ due to the relations. In particular, we can show that $h'_{n-1}$ products vanish using the fact that the square of each generator vanishes together with the other relations. Finally we obtain $h'_{n+1}=3h'_n -h'_{n-1}$ which gives $h'_{n+1}=5(3h_{n-1}- h_{n-2})=5h_n$. 
\end{proof}
Using this claim equation \eqref{PBWequation} becomes 
 \begin{equation*}
1+\sum_{n=1}^{+\infty}5h_{n-1}z^n =  \frac{\prod_{n=0}^\infty (1+z^{2n+1})^{\lambda_{2n+1}}}{\prod_{n=1}^\infty (1-z^{2n})^{\lambda_{2n}}}.
\end{equation*}
Since the Betti numbers $h_n$ satisfy $h_n=3h_{n-1}-h_{n-2}$ for $n \geq 2$ with $h_0=1$ and $h_1=3$, it is easy to verify that the left hand side of the last formula is equal to $$ (1+z)^2\sum_{n=0}^{+\infty}h_nz^n.$$ Hence, again by the  Poincaré--Birkoff--Witt Theorem applied to $X={\CP^2\#\,2\,\overline{\CP}\,\!^2}$, we conclude that $\lambda_n=r_n$ for $n \geq 2$.

Now consider the evaluation fibration
$$\xymatrix@1{ \Omega X \ar[r]^-{f} & \Symp_p(\tMuc) \ar[r]^-{g} & \Symp(\tMuc)},$$
where the base is equivalent to $T^2$ by Pinsonnault \cite{Pi}. 
Passing to rational homotopy we obtain 
$$\xymatrix@1{ \pi_*(\Omega X) \otimes \Q \ar[r] & \pi_*(\Symp_p(\tMuc))\otimes \Q  \ar[r] & \Q[x_0,y_0]}$$
where $ \Q[x_0,y_0]$ is the Lie graded algebra  with generators $x_0,y_0$ of degree 1 such that all Samelson products vanish. Note that the maps $\pi_*(f) \otimes \Q$ and $\pi_*(g) \otimes \Q$ are morphisms of graded Lie algebras and the second one is clearly surjective. Recall, from Proposition \ref{stabilizer}, that $\Gucc$ is homotopy equivalent to the stabilizer $\Symp_p(\tMuc)$. Therefore we can consider the following diagram 
\begin{equation}\label{algebra_diagram}
\xymatrix{ \pi_*(\Omega X) \otimes \Q \ar[r] \ar@{=}[d] & \pi_*(\Gucc)\otimes \Q  \ar[r] & \Q[ x_0,y_0]   \ar@{=}[d] \\
\pi_*(\Omega X) \otimes \Q \ar@{^{(}->}[r]^-{j} & \Lambda \ar[u]_i \ar[r] & \Q[x_0,y_0] }
\end{equation}
where the map $i$ is a morphism of algebras that simply maps the generators of $\Lambda$ injectively into $\pi_*(\Gucc)\otimes \Q$. Clearly there is a map  $j$ that makes the diagram commutative, because the map $\pi_n(f)\otimes \Q$ is injective  if $n=1$, it is an isomorphism if $n \geq 2$,  and $\lambda_n =r_n$ if $n\geq 2$ while $r_1=3$. So, for example, the three elements of degree one in $\pi_*(\Omega X) \otimes \Q$ are mapped to $x_1,y_1,z \in \Lambda$.  Hence, we can conclude that the map $i$ is an isomorphism of algebras. This concludes the proof for the case $\mu=1$.

\subsubsection{Proof of statements (S2) and (S4)}
If $1 < \mu \leq 2 $ and $ \lambda \leq  c_2$  then there are ten distinct configurations of symplectic embedded surfaces in $\tMucc$, namely the six that existed already when $\mu=1$  plus the four configurations described in Figure \ref{configB-MF}. Therefore we will have the five tori $T_i(0)$ we had before plus three new torus actions represented by the tori $T_i(1)$ where $i=1,4,5$  and one $S^1$--action that corresponds to  configuration (8) in Figure \ref{configB-MF}. Recall, from Proposition \ref{homotopygroups}, that when $\mu >1$ there is one more generator of the fundamental group. We may identify it with the element $t=x_{1,1}+y_{1,1} \in \pi_1(T_1(1))$ of Remark \ref{rmk:identifications}.  Since the classes $x_{1,i},y_{1,i}$ now have commuting representatives in $T_i(1)$  with $i=1,4,5$, it follows from Remark \ref{rmk:identifications} that 
\begin{gather*}
0=[x_{1,1},y_{1,1}]=[y_0-x_0,t+x_0-y_0]=[y_0,t]-[x_0,t],\\
0=[x_{1,4},y_{1,4}]=[y_1,t-x_1-2y_1+z]=[y_1,t], \\
0=[x_{1,5},y_{1,5}]=[y_1,t-y_1]=[y_1,t].
\end{gather*}
Moreover, from Lemma \ref{s1action}, one knows that the generator of the circle action corresponding to configuration (8) is $a_1$ and it commutes with $x_{1,4}$ and $y_{1,4}$. This gives the following relations
\begin{gather*}
0=[a_1,x_{1,4}]=[t-x_1+z,y_1]=[t,y_1],\\
0=[a_1,y_{1,4}]=[t-x_1+z,t-x_1+z-2y_1]=-2[t,x_1]+2[t,z]-2[z,x_1].
\end{gather*}
Note that this yields the relations of statement (S1) in the theorem. 

 Once  $\lambda$ passes $c_2$ with $\mu \in (1,2]$ then we need to consider two more configurations: (11) and (12) of Figure \ref{configB-MF-E2}. Hence there are two more torus actions, $T_2(1)$ and $T_3(1)$. The classes $x_{1,i},y_{1,i}$, $i=2,3$, now have commuting representatives in $T_i(1)$, so that gives the following relations
\begin{gather*}
0=[x_{1,2},y_{1,2}]=[y_0-x_1,t+x_0-y_0]=[y_0,t]-[x_1,t]-[x_1,x_0]\quad \mbox{and}\\
0=[x_{1,3},y_{1,3}]=[z-x_0,t-x_1+2x_0-z]=[z,t]-[z,x_1]-[x_0,t]+[x_0,x_1].
\end{gather*}
Since $a_1$  commutes with $x_{1,3}$ and $y_{1,3}$ it follows that 
\begin{gather*}
0=[a_1,x_{1,3}]=[t-x_1+z,z-x_0]=[t,z]-[t,x_0]-[x_1,z]+[x_1,x_0],\\
0=[a_1,y_{1,3}]=[t-x_1+z,t-x_1+2x_0-z]=-2[t,x_1]+2[t,x_0]-2[x_0,x_1],
\end{gather*}
which gives the relations in statement (S2).  Therefore for this algebra we have $\pi_2=\Q^6$. Note that for the case $\mu=1$ the homotopy groups satisfy $\pi_1=\pi_2=\Q^5$. Here we assume that the extra classes in $\pi_1$ and $\pi_2$   are represented  by $t$ and $[x_1,t]$, respectively. Moreover, it is not hard  to verify that all Samelson products in dimension three involving one of these classes do not add more generators to $\pi_3$, since they either vanish or are the linear combination of some other generators. Therefore we may conclude that the rank in each dimension of this algebra agrees with the rank of the homotopy groups described in Proposition \ref{homotopygroups}, when $1<\mu\leq 2$ and $c_2 <\lambda \leq c_1$. A similar argument to the one used in the case $\mu=1$ shows that we have indeed an isomorphism of algebras. More precisely, we just need to substitute the algebra 
$ \Q[x_0,y_0] $ of diagram \eqref{algebra_diagram} by the homotopy Lie algebra $\pi_*(\Symp(\tMuc))\otimes \Q$, that now includes the generator $t$. This concludes the proof of statement (S2).

Once $\lambda$ passes $c_1$ with $\mu \in (1,2]$ there are four configurations we should consider: (13) to (16) in Figure \ref{configB-MF-E1}. Hence the classes  $x_{2,i},y_{2,i}$, $i=1,4,5$, now have commuting representatives in $T_i(2)$ which gives the following relations
\begin{gather*}
0=[x_{2,1},y_{2,1}]=[t-x_1-x_0,x_1+y_0]=[t,x_1]+[t,y_0]-[x_0,x_1],\\
0=[x_{2,4},y_{2,4}]=[t-x_1+y_1-z,2z-y_1]=2[t,z]-2[x_1,z], \\
0=[x_{2,5},y_{2,5}]=[t-x_1-x_0,x_0+z]=[t,x_0]+[t,z]-[x_1,x_0]-[x_1,z].
\end{gather*}
The generator of the circle action corresponding to configuration (14) is $a_2$ and it commutes with $x_{2,4}$ and $y_{2,4}$.  That is, 
\begin{gather*}
0=[a_2,x_{2,4}]=[2t-2x_1+y_1,t-x_1+y_1-z]=-4[t,x_1]-2[t,z]+2[x_1,z],\\
0=[a_2,y_{2,4}]=[2t-2x_1+y_1,2z-y_1]=4[t,z]-4[x_1,z].
\end{gather*}
These yields  $[t,x_1]=0$, and hence the relations in statements (S3) and (S4). Note that 
for  $\lambda >c_1+c_2$ and  $\mu \in (1,2]$ or for all $\mu >2$, the Samelson products $[x_{2k-1,i},y_{2k-1,i}]$ and $[x_{2k,i},y_{2k,i}]$  vanish if we assume the above relations. So we do not get more relations, coming from these Samelson products, between the generators $x_0,y_0,x_1,y_1,z,t \in \pi_1(\Gucc)\otimes \Q$. Therefore the rank on each dimension of the algebra generated by these six elements, satisfying the relations in (S3) and (S4), agrees with the rank on each dimension of $\Lambda$, except in dimension 1. By comparison with Proposition \ref{homotopygroups}, it follows that in case (S4) we just need to add one more generator $w_\ell$ of the right degree, either $4\ell-2$ or $4\ell$, so the rank on each dimension of these algebras is the same to  the rank of the homotopy groups described in the proposition. Moreover we can use again diagram \eqref{algebra_diagram}, where now the Lie graded homotopy algebra of $\Symp(\tMuc)$ is given by $\Q[x_0,y_0,t,w_\ell]$,  to show that there is indeed an isomorphism of algebras and this concludes the proof of statement (S4). 

\subsubsection{Proof of statements (S1) and (S3)}
In order to conclude the proof of statements (S1) and (S3) we need to use the computation of  the rational cohomology ring of the classifying space of $\Gucc$, given in Theorem \ref{CohomologyRingClassSpace}, whose proof is postponed  to Section \ref{algebraiccomputations}
Note that these statements correspond to statements (a) and (c) of this Theorem (statement (S1) is just a particular case, since the two extra generators are given by the Samelson products $[x_1,t]$ and $[x_0,t]=[y_0,t]$).  In each one of these cases there are two relations in the cohomology ring of $B\Gucc$:  in  degree $4\ell$  in case (a); and $4\ell+2$ in case (c). Therefore the minimal model of the classifying space will have two new generators: in degree $4\ell-1$ and $4\ell+1$, respectively.  So the homotopy groups of $\Gucc$ have two new elements in degree $4\ell-2$ in case (a), and $4\ell$ in case (c), as we claim. Adapting diagram \eqref{algebra_diagram} again, one concludes the proof. 

Note that cases (b) and (d) of Theorem \ref{CohomologyRingClassSpace} agree with the result we obtained in statement (S4) of the main theorem, that is, this one relation in the cohomology ring of the classifying space corresponds precisely, to the { \it jumping} generator $w_\ell$ in the homotopy Lie algebra. 

\subsection{The limit at infinity}\label{limitinfinity}
As D. McDuff showed in ~\cite{MD-Haefliger} we can understand the limit $G_{\infty}= \lim_\mu \Gucc$ by studying the space $\cup_\mu \Aaucc$ since the  homotopy limit  $\lim_\mu \Aaucc$ is homotopy equivalent to the union $\Aa_\infty=\cup_\mu \Aaucc$ and $G_\infty$ is homotopy equivalent to the homotopy fiber of the evaluation map  $\Diff_0 \to \Aa_\infty$. More precisely, the map  $\Diff_0 \to \Ssucc$ lifts to a map 
$$\Diff_0 \to \Xxucc: \quad \phi \mapsto (\phi_*(\omega_\mu),\phi_*(J_0))$$ where $J_0$ is the standard complex structure contained in the open stratum. Composing with the projection to $\Aaucc$ we get a map $\Diff_0 \to \Aaucc$. Taking the homotopy limit we obtain the desired homotopy fibration.  Each $J \in \Aa_\infty$ gives a foliation of $\tMucc$ by $J$--spheres in class $F$, so the limit $\Aa_\infty$ is given by the stratification on $F$--curves. This implies that $\Aa_\infty$ is a disjoint union of sets of almost complex structures  for which the classes $E_1-E_2$ and $F-E_1-E_2$ are both represented, or one  of them is represented or none  of them is represented by $J$--holomorphic spheres, that is, 
$$ \Aa_\infty = \Aa_{none} \sqcup \Aa_{E_1-E_2} \sqcup \Aa_{F- E_1-E_2} \sqcup \Aa_{both}.$$
So each one of these four sets corresponds to one foliation by $F$--curves
and the limit at infinity is given by $$ G_\infty=\{ \phi \in \Diff_0: \phi \mbox{ leaves the $F$--foliation invariant} \}.$$

Note that, by Proposition \ref{Inflation}, we know on one hand that there are maps, defined up to homotopy,  $\Gucc \to G_{\mu + \epsilon, c_1, c_2}$ and on  the other hand the homotopy type of $\Gucc$ remains constant as $\mu$ varies in certain intervals, so we can pick one $\mu$ in each one of these intervals and consider the  sequence of the corresponding topological groups. The mapping telescope of this sequence is precisely the limit $G_\infty$. Moreover, it follows easily from the definition of the mapping telescope that we can consider only the groups that are stabilizers, as given in Proposition \ref{stabilizer} and drop the remaining ones from the sequence and the resulting mapping telescope is still homotopy equivalent to the limit $G_\infty$.  This observation together with  statement (S4) of Theorem \ref{maintheorem} (the one for which the proof is complete) implies that $\pi_*(G_\infty) \otimes \Q$ is generated, as a Lie graded algebra, by $x_0,y_0,x_1,y_1,z,t$ where all Samelson products vanish except the ones of the set $S_3$.

\begin{prop}\label{RingInfinity}
The rational cohomology ring of the classifying space $BG_\infty$ is given by 
$$ H^*(BG_\infty;\Q)= \Q[X_0,Y_0,X_1,Y_1,Z,T] /I$$
where all generators have degree 2 and $I$ is the ideal $ \langle \, X_0X_1, Y_0Y_1, X_0Y_1, Y_0Z, X_1Z \, \rangle$.
\end{prop}
\begin{proof}
It follows from the previous section that $\pi_{n+1}(BG_\infty)\otimes \Q = \pi_n(G_\infty) \otimes \Q = \Q^{r_n}$ if $n \geq 2$, where $r_n = \dim \pi_n (\Omega (\PbP))$ was defined in Section  \ref{stabilizer}, and  $H^2 (BG_\infty; \Q)=\pi_2(BG_\infty) \otimes \Q= \pi_1(G_\infty) \otimes \Q=\Q^6$. Let denote these six generators by $X_0,X_1, Y_0,Y_1,Z,T$. Since $r_2=5$ there are five generators of degree 3, namely $W_1,W_2,W_3,W_4$ and $W_5$, in the minimal model of $BG_\infty$. The differential of this generators is in degree 4 which implies that the differential is quadratic and therefore it is dual to the Samelson product in $\pi_*(G_\infty) \otimes \Q $ (see Section 13.(e) in \cite{FHT}). Since the non-vanishing Samelson products in  $\pi_2(G) \otimes \Q$ are given by 
$$ [x_0,x_1], \quad [y_0,y_1], \quad [x_0,y_1], \quad [y_0,z], \quad[x_1,z],$$ we can assume that the differential on the minimal model satisfies 

\begin{equation}\label{relations}
d (W_1)= X_0X_1, \quad d(W_2)=Y_0Y_1,  \quad d(W_3)=X_0Y_1,  \quad d(W_4) =Y_0Z,  \quad d(W_5)= X_1Z.
\end{equation}
These are the relations that generate the ideal $I$. In order to show that there are no more generators of $I$  consider the Serre spectral sequence of the fibration 
$$ G_\infty  \longrightarrow  EG_\infty \longrightarrow   BG_\infty.$$
The restriction to the primitive part of  $H^q(G_\infty)$ of the differential 
$$d_2^{0,q}:H^0(BG_\infty) \otimes H^q(G_\infty) \longrightarrow H^2(BG_\infty) \otimes H^{q-1}(G_\infty)$$
is dual to the Samelson product 
$$ H_2(BG_\infty) \otimes H_{q-1}(G_\infty) \simeq \pi_1(G_\infty) \otimes \pi_{q-1}(G_\infty) \otimes \Q \longrightarrow \pi_q(G_\infty) \otimes \Q.$$
So not only we can confirm that \eqref{relations} gives the generators of $I$, but also that we can check that, in higher dimensions, the differential of these primitive generators of the cohomology ring of $G_\infty$ correspond to the iterative Samelson products in higher degrees. These give relations that are multiples of  the ones in degree 4. For instance, the degree 3 Samelson product  $[z,[x_0,x_1]]$ gives rise to the relation $ZX_0X_1=0$ in degree 6, which belongs to the ideal $I$. Moreover, the  spectral sequence collapses at the second term and this concludes the proof. 
\end{proof}

\subsection{Algebraic computations}\label{algebraiccomputations}

This section is devoted to the computation of the rational cohomology ring of the classifying space of $\Gucc$ when $\mu \geq 1 > c_1+c_2 > c_1 >c_2$.
The idea to compute the cohomology ring of the classifying space $B\Gucc$ is the same as used in \cite[Section 5]{AGK} or \cite[Appendix B]{ALP} to compute the cohomology of $B\Symp(M_{\mu,c_1})$. 
Since we know that the map $\psi_{\lambda,\ell}: B\Gucc  \to BG_\infty$ induces a surjection in rational cohomology, it follows from Proposition \ref{RingInfinity} that $H^*(B\Gucc; \Q)\simeq \Q[X_0,Y_0,X_1,Y_1,Z,T]/(I \cup I_{\lambda,\ell})$ where $I_{\lambda,\ell}$ is the kernel of $\psi_{\lambda,\ell}^*$. The homotopy decomposition of $B\Gucc$ yields, at the rational cohomology level, an extended pullback diagram, that, for instance, in the case $c_1< \lambda \leq c_1+c_2$ is given by

\begin{equation}\label{ExtendedDiagramCohomology}
\xymatrix{
\Q[a_{2k}]/{\langle {e}_{2k}\rangle } \bigoplus_{i=1,4,5}\Q[x_{2k,i},y_{2k,i}]/{\langle {e}_{2k,i}\rangle }  & \ar[l]^-{\pi^*}  \Q[a_{2k}] \bigoplus_{i=1,4,5}\Q[x_{2k,i},y_{2k,i}]&\\
 H^*(BG_{\mu',c_1,c_2}); \Q)\ar[u]^{j_{2k}^{*}} & \ar[l]^{} \ar[u]^{i_{2k}^{*}} H^*(B\Gucc; \Q)  & \\
  &  & \ar[ul]_-{\psi_{\lambda,\ell}^*} \ar[ull]^{\psi_{\lambda',\ell'}^*}  \ar[uul]_{\psi_{2k,1}^*} H^*(BG_\infty;\Q)}
\end{equation}
where $e_{2k}$ is the Euler class of the representation of $S^1(2\ell)$ and $e_{2k,i}$ are the Euler classes of the representations of $T_i(2\ell)$.
So, in order to compute the ideal $I_{\lambda,\ell}$, one has to understand the maps 
\begin{align*}
& \psi^*_0: H^*(BG_\infty; \Q) \to \oplus_{i=1}^5 \Q[x_{0,i},y_{0,i}] \oplus \Q[a_0] \\
& \psi^*_{n,1}: H^*(BG_\infty; \Q) \to \oplus_{i=1,4,5} \Q[x_{n,i},y_{n,i}] \oplus \Q[a_n] \\
& \psi^*_{n,2}: H^*(BG_\infty; \Q) \to \oplus_{i=2,3} \Q[x_{n,i},y_{n,i}] 
\end{align*}
 for all $ n> 0$. For that it is enough to consider the relations in $\pi_1(G_\infty)$ between the generators of $\pi_1(T_i(n))$, $i=1, \dots, 5$.
So, from Lemma \ref{classification} and Remark \ref{rmk:identifications}, we get the definition of these maps. Recall that we use the same  notation for the generators of $\pi_1(T_i(n))$ and of $H^*(T_i(n);\Q)$.

\newcolumntype{F}{%
>{\color{blue}
\columncolor{white}[.6\tabcolsep]}l|}
\newcolumntype{A}{%
>{\color{black}
\columncolor[rgb]{0.96,0.96,0.9}[0.9\tabcolsep]}c}

\newcolumntype{P}{%
>{\color{black}
\columncolor[rgb]{0.9,0.9,0.96}[0.9\tabcolsep]}c}

\newcolumntype{G}{%
>{\color{blue}
\columncolor{white}[.6\tabcolsep]}c}

\begin{prop}\label{kernels}
The map $\psi^*_0$ is given by the following expressions
$$  \psi_0^*(X_0)=x_0 \quad \psi_0^*(Y_0)=y_0 \quad \psi_0^*(X_1)=x_1 \quad \psi_0^*(Y_1)=y_1\quad \psi_0^*(Z)=z \quad \psi_0^*(T)=0 $$ and the maps $\psi^*_{n,i}$, where $i=1,2$  and $n >0$, are given by Tables \ref{Table1} and \ref{Table2}.
\begin{table}[!ht]
\setlength{\extrarowheight}{.2cm}
\begin{tabular}{l!{\vrule}ccc}
%
\multicolumn{1}{F}{} &
%
\multicolumn{1}{G}{$ \psi^*_{2k-1,1}(...)=$} &
\multicolumn{1}{G}{$\psi^*_{2k-1,2}  (...)=$} &
\\ \hline

\multicolumn{1}{F}{$X_0$} &
\multicolumn{1}{A}{$(-kx_{2k-1,1}+ky_{2k-1,1},0,0,0)$}   &
\multicolumn{1}{P}{$(y_{2k-1,2}, -kx_{2k-1,3}+(k+1)y_{2k-1,3}) $} &
\\ \hline

\multicolumn{1}{F}{$Y_0$} &
\multicolumn{1}{A}{$(x_{2k-1,1}-y_{2k-1,1},0,0,0)$}   &
\multicolumn{1}{P}{$(x_{2k-1,2}-y_{2k-1,2},0)$} &
\\ \hline

\multicolumn{1}{F}{$X_1$} &
\multicolumn{1}{A}{$((1-k)y_{2k-1,1},-ky_{2k-1,4},(1-k)x_{2k-1,5},-k^2a_{2k-1}) $}   &
\multicolumn{1}{P}{$(-kx_{2k-1,2},-ky_{2k-1,3})$} &
\\ \hline

\multicolumn{1}{F}{$Y_1$} &
\multicolumn{1}{A}{$(0,kx_{2k-1,4}-(k+1)y_{2k-1,4},x_{2k-1,5}-y_{2k-1,5},0)$}   &
\multicolumn{1}{P}{$(0,0) $} &
\\ \hline

\multicolumn{1}{F}{$Z$} &
\multicolumn{1}{A}{$(0,(1-k)x_{2k-1,4}+ky_{2k-1,4},0,a_{2k-1}) $}   &
\multicolumn{1}{P}{$(0,x_{2k-1,3}-y_{2k-1,3}) $} &
\\ \hline

\multicolumn{1}{F}{$T$} &
\multicolumn{1}{A}{$(ky_{2k-1,1},ky_{2k-1,4},ky_{2k-1,5},k^2a_{2k-1}) $}   &
\multicolumn{1}{P}{$(ky_{2k-1,2},ky_{2k-1,3}) $} &

\end{tabular}
\bigskip

\caption{Definition of the maps $\psi^*_{2k-1,1}$ and $\psi^*_{2k-1,2}$.}
\label{Table1}
\end{table}

\begin{table}[!ht]
\setlength{\extrarowheight}{.2cm}
\begin{tabular}{l!{\vrule}ccc}
%
\multicolumn{1}{F}{} &
%
\multicolumn{1}{G}{$ \psi^*_{2k,1}(... )=$} &
\multicolumn{1}{G}{$\psi^*_{2k,2}  (...)=$} &
\\ \hline

\multicolumn{1}{F}{$X_0$} &
\multicolumn{1}{A}{$(-x_{2k,1},0,-kx_{2k,5}+ky_{2k,5},0)$}   &
\multicolumn{1}{P}{$(ky_{2k,2}, 0) $} &
\\ \hline

\multicolumn{1}{F}{$Y_0$} &
\multicolumn{1}{A}{$(y_{2k,1},0,0,0)$}   &
\multicolumn{1}{P}{$(y_{2k,2},0)$} &
\\ \hline

\multicolumn{1}{F}{$X_1$} &
\multicolumn{1}{A}{$ -(kx_{2k,1}- ky_{2k,1},kx_{2k,4},kx_{2k,5},k(k+1)a_{2k}) $}   &
\multicolumn{1}{P}{$-((k+1)x_{2k,2},(k+1)x_{2k,3}-(k+1)y_{2k,3}) $} &
\\ \hline

\multicolumn{1}{F}{$Y_1$} &
\multicolumn{1}{A}{$(0,x_{2k,4}-ky_{2k,4},0,a_{2k})$}   &
\multicolumn{1}{P}{$(0,y_{2k,3})$} &
\\ \hline

\multicolumn{1}{F}{$Z$} &
\multicolumn{1}{A}{$(0,-x_{2k,4}+(k+1)y_{2k,4},y_{2k,5},0) $}   &
\multicolumn{1}{P}{$(0,0) $} &
\\ \hline

\multicolumn{1}{F}{$T$} &
\multicolumn{1}{A}{$(kx_{2k,1},kx_{2k,4},kx_{2k,5},k(k+1)a_{2k}) $}   &
\multicolumn{1}{P}{$(kx_{2k,2},kx_{2k,3}) $} &

\end{tabular}
\bigskip

\caption{Definition of the maps $\psi^*_{2k,1}$ and $\psi^*_{2k,2}$.}
\label{Table2}
\end{table}

Their kernels are the ideals
\begin{align*}
& \ker \psi^*_0  =\langle \,T\, \rangle, \\
& \ker \psi^*_{2k-1,1}  = \langle \, X_0+kY_0, \, k(X_1+kZ)+(k-1)(T+kY_1) \, \rangle, \\
 & \ker \psi^*_{2k-1,2}  = \langle \, Y_1, kY_0+X_1+T,\, kX_0+k^2Z-T \, \rangle \\
 & \ker \psi^*_{2k,1}   = \langle \, -kY_0+X_1+T, \, kX_0-k(k+1)Y_1-k^2Z+T \,  \rangle, \\
&  \ker \psi^*_{2k,2}  = \langle \, Z,X_0-kY_0, \, kX_1+(k+1)T-k(k+1)Y_1\, \rangle.
\end{align*}
\end{prop}
In order to finish the proof we now argue by induction on $m$ where $m=N-1$ if $N \geq 3$ and  $N$ is defined in \eqref{strata}, while $m=1$  if $\mu=1$. Fix some $c_1,c_2$. The proof for the case $\mu=1$ is similar to the proof of Proposition \ref{RingInfinity}, it is only necessary to remove the generator $T$. Now assume the result holds for some $\mu' >1$ for which the induction step corresponds to $m-1$ and consider a $\mu > \mu'$ corresponding to $m$.  If $c_1+c_2 < \lambda$ then $c_1 < \lambda' \leq c_1+c_2$ . Then diagram \eqref{ExtendedDiagramCohomology} implies that $I_{\lambda,\ell} \subset I_{\lambda',\ell} \cap \ker \psi^*_{2l,2} $. Note that in this case the homotopy Lie algebra of $\Gucc$ is known. In particular we know there is a new generator in degree $4\ell$ which should correspond to a relation in degree $4\ell+2$ in the cohomology ring of the classifying space. Since the relations for the case $m-1$ are already in degree $4\ell+2$ and the generators of the ideals $ \ker \psi^*_{2l,1} = \langle b_1,b_2 \rangle$  and $ \ker \psi^*_{2l,2}=\langle c_1,c_2, c_3 \rangle$ of Proposition \ref{kernels} satisfy the relation $\ell b_1+b_2=-\ell^2c_1+\ell c_2+c_3$, we conclude that  $I_{\lambda,\ell} = I_{\lambda',\ell} \cap \ker \psi^*_{2l,2} $ and the statement follows. The argument is the same to prove the case when $c_2 < \lambda \leq c_1$. In the cases $\lambda \leq c_2$ or $c_1 < \lambda \leq c_1+c_2$  we know that  $I_{\lambda',\ell}$ and $ \psi^*_{2l-1,2} $ or  $\psi^*_{2l,1}$ , respectively, are coprime. Therefore  $I_{\lambda',\ell} \cap \ker \psi^*_{2l-1,1}= I_{\lambda',\ell} . \ker \psi^*_{2l-1,1} $ and $I_{\lambda',\ell} \cap \ker \psi^*_{2l,1}= I_{\lambda',\ell} . \ker \psi^*_{2l,1} $.  Since in the next step of induction the relations are already  in degrees  $4\ell$ and $4\ell+2$, respectively, it follows that the same holds for these cases. So $I_{\lambda,\ell}=I_{\lambda',\ell} . \ker \psi^*_{2l-1,1} $ or $I_{\lambda,\ell}=I_{\lambda',\ell} . \ker \psi^*_{2l,1} $, depending on the size of $\lambda$. This concludes the proof.

\subsection{Consequences of Theorem \ref{maintheorem}}
As an immediate corollary of Theorem \ref{maintheorem} we  obtain the Pontryagin ring with rational coefficients of the group of symplectomorphisms of $\tMucc$. Recall that if $G$ is a topological group then the Samelson product in $\pi_*(G)$ is related to the Pontryagin product in $H_*(G;\Z)$, which is induced by the product in $G$, by the formula
$$[x,y]=xy-(-1)^{\deg x \deg y}yx, \quad \quad x,y \in \pi_*(G),$$
where we suppressed the Hurewicz homomorphism $\rho:\pi_*(G) \to H_*(G;\Z)$ to simplify the expression. Therefore we will use $[x,y]$ to denote both the Samelson product and the commutator in homology. 
Denote by $\Q \langle x_1,\hdots,x_n \rangle $ the free non-commutative algebra over $\Q$ with generators $x_j$.

\begin{cor}\label{PontryaginRing}
If $\mu = 1 > c_1+c_2 > c_1 >c_2$ then the Pontryagin ring of $\Gucc$ is given by
$$ H_*(\Gucc;\Q)=\Q \langle x_0,y_0,x_1,y_1,z \rangle / R $$ where
where all the generators have degree 1 and  $R$ is the set of relations
$$R=\{ [x_0,y_0]=[x_0,z]=[x_1,y_1]=[x_1,y_0]=[z,y_1]=0, \, x_i^2=y_i^2=z^2=0, \ i=0,1 \}.$$

If $\mu= \ell+\lambda $, where $\ell \in \N$,  $\lambda \in (0,1]$ and   $1> c_1 + c_2 > c_1 > c_2 $ then the Pontryagin ring is given as follows.
\begin{enumerate}
\item If $1 < \mu \leq 2$ and $\lambda \leq c_2$ then 
$$  H_*(\Gucc;\Q)=\Q \langle x_0,y_0,x_1,y_1,z,t \rangle / R_1$$
where $\deg t=1$ and 
$$R_1=  R \cup  \{ \, [y_1,t]=0, \, [x_0,t]=[y_0,t], \, [z,t]=[x_1,t]+[x_1,z] \, \}.$$

\item If $1 < \mu \leq 2$ and $c_2  < \lambda \leq c_1$ then 
$$  H_*(\Gucc;\Q)=\Q \langle x_0,y_0,x_1,y_1,z,t \rangle / R_2$$
where 
$$R_2= R_1 \cup   \{ \, [x_0,t]=[y_0,t]=[x_1,t]+[x_1,x_0] \, \}.$$

\item If 
 $(a) \ \ell > 1 \ {\mbox and} \ \lambda \leq c_2, \ {\mbox or}\quad   (b) \ \ell \geq 1   \ {\mbox and} \  c_1 < \lambda \leq c_1+c_2, $
 then 
$$  H_*(\Gucc;\Q)=\Q \langle x_0,y_0,x_1,y_1,z,t,w_{\ell_1},w_{\ell_2} \rangle / R_3$$
where  
$$\deg w_{\ell_1}= \deg w_{\ell_2} = \left\{ \begin{array}{ll}
 4 \ell-2 & \mbox{in case }  (a)\\
 4\ell & \mbox{in case } (b),
\end{array} \right. 
 $$
and $R_3$ is the set of relations
$$R_3=R_2 \cup \,  \{\, w_{\ell_1}^2=w_{\ell_2} ^2=0,\ [x_1,t]=0,\,\mbox{and $w_{\ell_1},w_{\ell_2} $ commute with all  generators} \, \}.$$
 
\item If 
$  (a) \ \ell > 1 \ {\mbox and} \ c_2<\lambda \leq c_1,  \ {\mbox or} \quad   (b) \  \ell \geq 1  \  {\mbox and}  \ c_1+c_2 < \lambda,$

 then
 $$  H_*(\Gucc;\Q)=\Q \langle x_0,y_0,x_1,y_1,z,t,w_\ell \rangle / R_3$$
 
 where
 $$\deg w_{\ell} = \left\{ \begin{array}{ll}
 4 \ell-2 & \mbox{in case } (a)\\
 4\ell & \mbox{in case }  (b).
\end{array} \right. $$

\end{enumerate}
\end{cor}

\begin{proof}
The result  follows immediately from   Theorem \ref{maintheorem} together with the Milnor--Moore Theorem applied to $X=B\Gucc$.
\end{proof}

\begin{remark} By the Cartan--Serre Theorem, if the rational homology of $\Gucc$ is finitely generated in each dimension, the rational cohomology is a Hopf algebra for the cup product and coproduct induced by the product in $\Gucc$, and it is generated as an algebra by elements that are dual to the spherical classes in homology. In particular, the number of generators of odd dimension $d$ appearing in the anti-symmetric part of the rational cohomology algebra  is equal to the dimension of $\pi_d(G)\otimes \Q$, and the number of generators of even dimension $d$ appearing in the symmetric part $S$ is equal to the dimension of $\pi_d(G)\otimes \Q$. We can then conclude from  Proposition \ref{homotopygroups} and Theorem \ref{maintheorem} that the rational cohomology algebra of $\Gucc$ is infinitely generated. 
\end{remark}
\medskip

\section{Remaining cases}\label{RemainingCases}
Up to this point we just studied the generic case, that is, we studied the topology of the symplectomorphism group $\Gucc$ when the sizes of the blow-ups satisfy $c_2 < c_1 < c_1+c_2 < 1$. In this section we study the remaining cases, which are the following: 
\begin{enumerate}
\item[(R1)] $\mu \geq 1, c_1=c_2$ and $c_1+c_2 <1$;
\item[(R2)] $\mu \geq 1, c_2 < c_1$ and $c_1+c_2=1$;
\item[(R3)] $\mu \geq 1$ and $c_2=c_1=\frac12$.
\end{enumerate}

\subsection{The cases (R1) and (R2)}
The argument to study these cases is similar to the one used to study the generic case. In the case (R1) the major difference is that since we have $c_1=c_2$ the class $E_1-E_2$ cannot be represented by a $J$--holomorphic sphere for all $J$ in the space of almost complex structures, that is, not all configurations from (1) to (18) in Section \ref{Jholomorphic} can be realized by $J$--holomorphic spheres. Therefore when the size of $\mu$ increases 
there are still new strata being added to the space of almost complex structures, but less than in the generic case. When $\mu=1$ the isometries are given by the tori $T_1(0), T_2(0)$ and $T_5(0)$. Furthemore, it is easy to describe the corresponding isometries of the complex structures in each stratum as the size of $\mu$ increases. Assuming that $\ell < \mu \leq \ell+1$,
\begin{itemize}
\item if $\lambda \leq c_2$ we add the isometries $T_1(2\ell-1)$ and $S^1(2\ell-1)$;
\item if  $c_2 <\lambda \leq c_1$ we add two more tori $T_2(2\ell-1)$ and $T_3(2\ell-1)$; 
\item if $c_1<\lambda \leq c_2+c_1$ there are again two tori $T_1(2\ell)$ and $T_5(2\ell)$; and
\item if  $c_1 +c_2<\lambda $ there is only one more torus $T_2(2\ell)$.
\end{itemize}
Repeating all the steps of the argument in the generic case, assuming the changes on the stratification of the space of almost complex structures we just discussed, it is easy to obtain the rational homotopy algebra in this case. A key point in the argument for the generic case is that the group $\Gucc$  has the homotopy type of a stabilizer for some values of $\lambda$ (cf. Proposition \ref{stabilizer}). This still holds if $c_1=c_2$ and $c_1+c_2 < 1$ when $ \lambda$ lies in the interval $(c_1+c_2, 1]$.  The idea of the proof is the same as in Proposition \ref{stabilizer}. The only difference is that to prove stability we cannot perform inflation along the curves $E_1-E_2$ because they are not represented. Finally we conclude that the homotopy algebra is just a subalgebra of the algebra of Theorem \ref{maintheorem} which is isomorphic to the subalgebra generated by all the generators except $y_1$. Note that this generator of the homotopy algebra in the generic case corresponds to a circle action that only appears as an isometry of the complex structures $J$  for which the class $E_1-E_2$ is represented by a $J$-holomorphic sphere. 

There is an analog argument for the case (R2). Here, since $c_1+c_2=1$ the class $F-E_1-E_2$ cannot be represented by a $J$--holomorphic sphere for all $J$. It follows that, in this case, the homotopy algebra of the symplectomorphism group $\Gucc$ is isomorphic to the subalgebra generated by all generators of the algebra of the generic case (cf. Theorem~ \ref{maintheorem} ) except $z$. 

\subsection{The case (R3)}
If $c_1=c_2= \frac12$ then both the classes $E_1-E_2$ and $F-E_1-E_2$ cannot be represented by $J$-- holomorphic spheres. It follows that  the stratification of the space of almost complex structures changes  only when $\mu $ passes an integer or $\lambda$ passes $\frac12$. We call these the {\it critical} values. If $\mu=1$ there is only one stratum which implies that the symplectomorphism group $\Gucc$ is homotopy equivalent to the group of isometries  of the complex structure in that stratum, which is a torus. Then, when we increase the size of $\mu$,  two strata are added once $\lambda$ or $\mu$ pass a critical value. More precisely, when $\mu$ passes an integer the isometries corresponding to the new strata are $T_2(2\ell-2)$ and $T_1(2\ell-1)$ and when $\lambda$ passes $\frac12$ the corresponding isometries are $T_1(2\ell)$ and $T_2(2\ell-1)$.  As before the stratification yields an homotopy decomposition of $B\Gucc$ that we describe for the case $\lambda >\frac12$ (there is a similar pushout diagram for the case $\lambda \leq \frac12$): 
\begin{prop}
If $\mu= \ell+\lambda $, where $\ell \in \N$ and $\lambda \in (0,1]$, $c_2 =c_1= \frac12$ and $ \lambda > \frac12$ there is a homotopy pushout diagram
$$
\xymatrix{
S^{4\ell-1}_{h{T_1}(2\ell)} \sqcup S^{4\ell-3}_{h{T_2}(2\ell-1)} \ar[r] \ar[d]^{j_{2\ell, 2\ell-1}} & B{T_1}(2\ell) \sqcup B{T_2}(2\ell-1)  \ar[d]^{i_{2\ell,2\ell-1}} \\
B{G}_{\mu',c_1,c_2} \ar[r] & B\Gucc}
$$
 where $\mu'= \ell + \lambda'$ satisfies $ \lambda'  \leq \frac12$,
so that $G_{\mu',c_1,c_2}$ is the group associated with a stratification having two strata less than the stratification associated with $G_{\mu,c_1,c_2}$. The upper horizontal map  is the universal bundle map associated to the representations of ${T_1}(2\ell)$ and ${T_2}(2\ell-1)$ on $H^{0,1}_{J_{1,2\ell}}(T\tMucc)$ and $H^{0,1}_{J_{2,2\ell-1}}(T\tMucc)$ respectively, $i_{2\ell, 2\ell-1}$ is induced by the inclusions ${T_1}(2\ell)\into \Gucc$ and ${T_2}(2\ell-1)\into \Gucc$, and the map $BG_{\mu',c_1,c_2}\to B\Gucc $ coincides, up to homotopy, with the one given by inflation described in Proposition ~\ref{Inflation}.
\end{prop}

In this case the group $\Gucc$ does not have the homotopy type of a stabilizer so we cannot compute its homotopy groups as before. In particular we need to find another way to study the limit at $G_\infty$. First note that the space $\Aa_\infty$ in this case is given by all tamed almost complex structures for which the set of classes $\{E_1, E_2, F-E_1, F-E_2 \}$ is always represented, so the $F$-- curves can only break as $E_i \cup (F-E_i)$ which gives us a foliation of $\tMucc$ by $F$--curves only with two singular fibers which are precisely $E_1 \cup (F-E_1)$ and $E_2 \cup (F-E_2)$. Since 
$$G_\infty = \{\phi \in \Diff_0\, : \, \phi \ \mbox {leaves the $F$-foliation invariant} \}$$ we conclude that 
$G_\infty$ is homotopy equivalent to the group of fiber preserving diffeomorphisms of $S^2\times S^2$ that fix two points, $\FDiff_{p,q}$, that we will denote simply by $G_{p,q}$. Next we observe that there is a fibration 
\begin{equation}\label{limit}
\xymatrix{ G_{p,q} \ar[r] & G_p \ar[r]^-{ev_q}& S^2}
\end{equation}
where   $ev_q(\phi)= \phi(q)$ and $G_p$ is the group of fiber preserving diffeomorphisms of $S^2\times S^2$ that fix one point, acting on the fiber containing $q$. This group has the rational homotopy type of $S^1\times S^1\times S^1$ (cf. ~\cite{ALP}). Therefore from the long exact homotopy sequence of fibration \eqref{limit} we obtain the homotopy groups of $G_{p,q}$ which are given by $\pi_1(G_{p,q})= \Q^4$, $\pi_2(G_{p,q})= \Q$ and the homotopy groups vanish if $n \geq 3$. It follows that the minimal model of the classifying space of $BG_{p,q}$ has four generators of degree 2 and one of degree 3, that we denote by $X_0,Y_0,X_1,T$ and $W$ respectively. In order to compute its differential 
we analyse the Serre spectral sequence of the fibration 
$$S^2 \to BG_{p,q} \to BG_p.$$
It is not hard to check that, for dimensional reasons, all differentials in the spectral sequence vanish. Therefore it splits and we can conclude that the differential of the minimal model is non--zero on the degree 3 generator $W$. Hence the rational cohomology ring of the classifying space $BG_{p,q}$ is given by 
$$ H^*(BG_{p,q}; \Q)=\Q[X_0,Y_0,X_1,T] / R$$
where $R$ contains a relation coming from a Samelson product  in $G_{p,q}$. Finally, to understand this relation we study the Serre spectral sequence of the fibration 
$$G_{p,q} \to EG_{p,q} \to BG_{p,q}.$$ 
By analysing the differential in the second page of the spectral sequence we can show that the relation can be written as $X_0X_1=0$. 
The final step now is to compute the rational cohomology ring of the symplectomorphism groups $B\Gucc$. Since we use the same technique as in the generic case  which involves considering diagrams in cohomology, as diagram \eqref{ExtendedDiagramCohomology}, we do not show the algebraic computations.  

\begin{thm} 
If $\mu=1$ and $c_1=c_2=\frac12$ then the rational cohomology ring of $B\Gucc$ is given by
$$ H^*(B\Gucc; \Q)= \Q[X_0,Y_0] $$ where $\deg X_0=\deg Y_0=2$. 

If $\mu= \ell+\lambda $, where $\ell \in \N$ and $\lambda \in (0,1]$ and $c_1=c_2=\frac12$ then 
$$ H^*(B\Gucc; \Q)= \Q[X_0,Y_0,X_1,T] / I \cup I_{\lambda,\ell}  $$
where $\deg X_1 =\deg T=2$, $I=\langle X_0X_1\rangle$ and $I_{\lambda,\ell}$ depends on $\ell$ and $\lambda$ as follows. Let 
$$
 A_k =  k(X_0+X_1+kY_0)+(k-1)T\quad {\mbox  and } \quad 
 B_k =   k(X_0+X_1-kY_0)+(k+1)T.
 $$
\begin{enumerate}[(a)]
\item If $\lambda \leq \frac12$ then 
$$ I_{\lambda,\ell} = \langle \, T \prod_{k=1}^{\ell-1}(A_k\, B_k )\, (X_0+ \ell Y_0), \  T \prod_{k=1}^{\ell-1}(A_k\, B_k) \, [(\ell -1)T+\ell X_1] \, \rangle.$$
\item If $\lambda >\frac12$ then 
$$ I_{\lambda,\ell}= \langle \, T \prod_{k=1}^{\ell-1}(A_k\, B_k) \, A_\ell (X_1+T-\ell Y_0), \ T \prod_{k=1}^{\ell-1}(A_k\, B_k) \, A_\ell [\ell X_0 +T] \,  \rangle  .$$
\end{enumerate}
\end{thm}

 Using the Serre spectral sequence of the fibration 
$$\Gucc \to E\Gucc \to B\Gucc$$ 
it follows from the previous computation that the rational cohomology ring of the group $\Gucc$ is finitely generated. Indeed, looking at the second page of the spectral sequence it is not hard to see that, for dimensional reasons, there must be a generator $w \in H^*(\Gucc ; \Q)$ of degree two. Moreover all elements in the second page of the spectral sequence are killed by the differential except two that correspond to the relations on the rational cohomology ring of the classifying space and which explain why there must be two more generators of even degree. We use the same notation for the generators as in the Pontryagin ring (or homotopy algebra) so it is easier 
to understand the correspondence between dual elements in this ring and the cohomology ring. 
\begin{prop}
If $\mu= \ell+\lambda $, where $\ell \in \N$ and $\lambda \in (0,1]$ and $c_1=c_2=\frac12$  then 
$$H^*(\Gucc ; \Q) =\Lambda (x_0,y_0,x_1,t) \otimes {\mathrm S}(w, w_{\ell_1}, w_{\ell_2})$$
where $\Lambda (x_0,y_0,x_1,t)$ is an exterior algebra on generators of degree one and ${\mathrm S}(w, w_{\ell_1}, w_{\ell_2})$ is a polynomial algebra on the generators $w, w_{\ell_1}, w_{\ell_2}$ such that $\deg w =2$ and 
$$\deg w_{\ell_1}= \deg w_{\ell_2}= \left\{ \begin{array}{cl}
4\ell-2 & \mbox{ if $\lambda \leq \frac12$} \\
4\ell & \mbox{ if $\lambda > \frac12$}.
\end{array} \right.
$$
\end{prop}

Note that the generator $w$ corresponds to the relation $X_0X_1=0$ in the cohomology ring of the classifying space (or to the Samelson product $[x_0,x_1]$ in the rational homotopy algebra) while the generators $w_{\ell_1}, w_{\ell_2}$ correspond to the two other relations in degrees $4\ell$ or $4\ell+2$ depending on $\lambda \leq \frac12$ or $\lambda > \frac12$,  respectively. 

\begin{remark} We note that Kedra's result in ~\cite{Ke} cannot be applied to this case (R3) and this is the only one for which the rational cohomology ring of the group $\Gucc$ is finitely generated.  For all the other cases the computation of the homotopy groups  implies that there are indeed an infinite number of generators of the cohomology ring, since the rank of the homotopy groups never vanishes.  
\end{remark}

\appendix

\section{Stability of symplectomorphism groups} \label{se:appendixA}
This appendix is devoted to prove Proposition \ref{prop:stability}.

The key tool in the proof  is the Inflation Lemma of Lalonde--McDuff and a generalization of it by Olguta Buse.  The Inflation Lemma states that in the presence of a $J$--holomorphic curve with nonnegative self-intersection, it is possible to deform a symplectic form $\tau_0$ through a family of forms $\tau_t$ that tame the almost complex structure $J$. Buse recently showed that, under some restrictions on the parameter $t$, actually the $J$--holomorphic curved used in the inflation procedure is allowed to have negative self-intersection.

\begin{lemma}[Inflation Lemma, see Lalonde ~\cite{La} and ~\cite{MD-Haefliger}]
Let $J$ be an $\tau_0$--tame almost complex structure on a symplectic 4--manifold $(M, \tau_0)$ that admits a $J$--holomorphic curve $Z$ with $Z\cdot Z \geq 0$. Then there is a family $\tau_t$, $t \geq 0$, of symplectic forms that all tame $J$ and have cohomology class $[\tau_t]=[\tau_0]+t {\rm PD}(Z)$, where ${\rm PD}(Z)$ is the Poincaré dual to the  homology class $[Z]$. 
\end{lemma}

\begin{lemma}[Buse ~\cite{Bu}]
Let $J$ be an $\tau_0$--tame almost complex structure on a symplectic 4--manifold $(M, \tau_0)$ that admits a $J$--holomorphic curve $Z$ with $Z\cdot Z =-m, \, m \in \N$. Then for all $\epsilon >0$ there is a family $\tau_t$ of symplectic forms,  all taming $J$, which satisfy $[\tau_t]=[\tau_0]+t {\rm PD}(Z)$ for all $0 \leq t \leq \frac{\tau_0(Z)}{m} - \epsilon$.
\end{lemma}

The proof of the proposition uses the same approach of the proof of Proposition 3.1 in ~\cite{Pi}, i.e., it is done in three steps. First we prove the following claim.
\begin{claim*}
The spaces $\Aaucc$ and $\Aa_{\mu',c_1,c_2}$ are equal if:
\begin{enumerate}
\item $\lambda \leq c_2 < c_1 < c_1+c_2$ and $ \mu' \in (\ell, \mu]$; or
\item $c_ 2 < \lambda \leq c_1 < c_1+c_2$ and $ \mu' \in (\ell + c_2, \mu]$; or
\item $c_ 2 < c_1 <\lambda \leq c_1+c_2$ and $  \mu' \in (\ell + c_1, \mu]$; or 
\item $c_ 2 < c_1 < c_1+c_2 < \lambda$ and $ \mu' \in (\ell + c_1+c_2, \mu]$.
\end{enumerate}
\end{claim*}

In the second step we show that the spaces $\Aaucc$ and  $\Aa_{\mu,c_1',c_2}$ are equal when the two parameters $c_1$ and $ c_1'$ satisfy 
 $c_1 \leq c_1'< \lambda$ or $ \lambda  \leq c_1 \leq c_1'$. 
The third step is similar to the second step and asserts that 
$\Aaucc=\Aa_{\mu,c_1,c_2'}$ whenever the parameters $c_2$  and $c_2'$ verify $c_2 \leq c_2'< \lambda$ or  $ \lambda  \leq c_2 \leq c_2'$.

{\it Step 1}: We begin by showing that for any $\mu \geq 1$ and any $\epsilon >0 $ one has $\Aaucc \subset \Aa_{\mu+\epsilon, c_1,c_2} $. The proof of this fact goes as in the case of 1--blow--up done in \cite{Pi}.  The idea is that it is always possible to inflate the form $\omucc$ along some embedded $J$--holomorphic sphere representing the class $F$, to get a one parameter family of symplectic forms in classes $[\omucc]+\epsilon \, \textrm{PD}(F)=\om_{\mu+\epsilon,c_1,c_2}$. All these forms tame $J$. 

We next show that for $\mu=\ell+\lambda \geq 1$ and $\lambda \leq c_2 < c_1 < c_1+c_2$ we have $\Aaucc \subset \Aa_{\mu',c_1,c_2}$ whenever $ \mu' \in (\ell, \mu]$. Since $\lambda \leq c_2$, there are almost complex structures  giving configurations (7) and (9) in Figure 2. Lemma \ref{representedclasses} implies that  the classes $D=E_1-E_2$, $D_{4\ell+1}=B+\ell F$ and $D_{4\ell-1}=B+\ell F-E_1$ are all represented by embedded $J$--holomorphic spheres. Therefore we can inflate $\omucc$ first along the curves $D_{4\ell+1}$  and $D_{4\ell-1} $ getting a two parameter family of symplectic forms 
$$ \om_{be}=\frac{\omucc+bd_{4\ell+1}+ed_{4\ell-1}}{1+b+e}$$
where $d_i$ is a 2--form representing the Poincaré dual of the homology class $D_i$, PD$(D_i)$ and  $b,e \geq 0$. We then have $\om_{be}(F)=1$,
$$ \om_{be}(B)=\frac{\mu+\ell(b+e)}{1+b+e}, \  \ \om_{be}(E_1)=\frac{c_1+e}{1+b+e} \ \ \mbox{and} \ \ \om_{be}(E_2)=\frac{c_2}{1+b+e}.$$ 
We can then inflate $\om_{be}$ along the curve $E_1-E_2$ getting a family of symplectic forms 
$$\om_a=\om_{be}+a{\rm PD}(E_1-E_2),$$
for all $0\leq a < \frac{\om_{be}(E_1-E_2)}{2}$, where ${\rm PD}(E_1-E_2)$ is a 2--form representing the Poincaré dual of the homology class $E_1-E_2$.
The area of the exceptional classes $E_i$, $i=1,2$ can be made equal to $c_i$ by setting $e=(c_1+c_2)b/(1-c_1-c_2) \geq 0$ and $a=c_2b/(1-c_1-c_2+b)$. Note that  $$2a <  \om_{be}(E_1-E_2) = \frac{(c_1-c_2)(1-c_1-c_2)+ b(c_1+c_2)}{1-c_1-c_2+b},$$ for all $b \geq 0$.  This gives a one parameter family of forms $\omega_b$ verifying $\om_b(F)=1$, $\om_b(E_1)=c_1$, $\om_b(E_2)=c_2$ and $$ \om_b(B)=\frac{\mu(1-c_1-c_2)+b\ell}{1-c_1-c_2 +b}.$$ Letting $ b \to \infty$, this shows that the almost complex structure $J$ is tamed by a symplectic form in class $[\om_{\ell+\epsilon,c_1,c_2}]$ for all $0 <\epsilon \leq \lambda$. Therefore $\Aaucc \subset \Aa_{\mu',c_1,c_2}$, for all $\mu' \in (\ell, \mu]$. 

In order to finish the proof of the claim, that is, to show that $\Aaucc \subset \Aa_{\mu',c_1,c_2}$, we need to consider all other possible cases which includes all configurations up to number (18).  For configurations (13) and (15) we use the same inflation process, while for the remaining cases the process is simpler since we do not need to use negative inflation. Nevertheless one needs  a different choice of curves, since it has to satisfy Lemma \ref{representedclasses}. The Tables \ref{Table3}, \ref{Table4} and \ref{Table5} summarize these choices as well as the values of the parameters. 

\begin{table}[!ht]
\setlength{\extrarowheight}{.25cm}
\begin{tabular}{l!{\vrule}ccccccc}
%
\multicolumn{1}{F}{Configuration} &
%
\multicolumn{1}{P}{$\lambda \leq c_2,\ (7) \mbox{ or } (9)$} &
\multicolumn{1}{P}{$c_1 < \lambda \leq c_1+c_2,\ (13) \mbox{ or } (15)$} &
\\ \hline

\multicolumn{1}{F}{Curves} &
\multicolumn{1}{A}{$D_{4\ell-1},D_{4\ell+1},E_1-E_2$}   &
\multicolumn{1}{A}{$D_{4\ell+3},D_{4\ell+1},F-E_1-E_2$} &
\\ \hline

\multicolumn{1}{F}{$\omega_{be}$} &
\multicolumn{1}{P}{$\frac{\omucc+bd_{4\ell+1}+ed_{4\ell-1}}{1+b+e}$}   &
\multicolumn{1}{P}{$\frac{\omucc+bd_{4\ell+1}+ed_{4\ell+3}}{1+b+e}$} &
\\ \hline

\multicolumn{1}{F}{$\omega_a$} &
\multicolumn{1}{A}{$\om_{be}+a {\rm PD}(E_1-E_2)$}   &
\multicolumn{1}{A}{$\om_{be}+a {\rm PD}(F-E_1-E_2)$} &
\\ \hline

\multicolumn{1}{F}{$a$} &
\multicolumn{1}{P}{$\frac{c_2b}{1-c_1-c_2+b}$}   &
\multicolumn{1}{P}{$\frac{c_2b}{1-c_1+c_2+b}$} &
\\ \hline

\multicolumn{1}{F}{$e$} &
\multicolumn{1}{A}{$\frac{(c_1+c_2)b }{1-c_1-c_2}$}   &
\multicolumn{1}{A}{$\frac{(c_1-c_2)b }{1-c_1+c_2}$} &
\\ \hline

\multicolumn{1}{F}{$\om_b(B)$} &
\multicolumn{1}{P}{$\frac{\mu(1-c_1-c_2)+b\ell}{1-c_1-c_2+ b}$}   &
\multicolumn{1}{P}{$\frac{\mu(1-c_1+c_2)+b(\ell+c_1)}{1-c_1+c_2+ b}$} &

\end{tabular}
\bigskip

\caption{Inflation process to show that $\Aaucc \subset \Aa_{\mu',c_1,c_2}$ when $\mu' \leq \mu$.}
\label{Table3}
\end{table}

\begin{table}[!ht]
\setlength{\extrarowheight}{.25cm}
\begin{tabular}{l!{\vrule}ccccccc}
%
\multicolumn{1}{F}{Configuration} &
%
\multicolumn{1}{P}{$\lambda \leq c_2, \ (8) \mbox{ or } (10)$} &
\multicolumn{1}{P}{$c_2 < \lambda \leq c_1, \ (11) \mbox{ or } (12)$} &
\\ \hline

\multicolumn{1}{F}{Curves} &
\multicolumn{1}{A}{$D_{4\ell-1},D_{4\ell},D_{4\ell+1}$}   &
\multicolumn{1}{A}{$D_{4\ell-1},D_{4\ell+1},D_{4\ell+4}$} &
\\ \hline

\multicolumn{1}{F}{$\omega_{abe}$} &
\multicolumn{1}{P}{$\frac{\omucc+ a d_{4\ell}+bd_{4\ell+1}+ed_{4\ell-1}}{1+a+b+e}$}   &
\multicolumn{1}{P}{$\frac{\omucc+ a d_{4\ell+4}+bd_{4\ell+1}+ed_{4\ell-1}}{1+a+b+e}$} &
\\ \hline

\multicolumn{1}{F}{$a$} &
\multicolumn{1}{A}{$\frac{c_2b}{1-c_1-c_2}$}   &
\multicolumn{1}{A}{$\frac{c_2b}{1-c_1-c_2}$} &
\\ \hline

\multicolumn{1}{F}{$e$} &
\multicolumn{1}{P}{$\frac{c_1b }{1-c_1-c_2}$}   &
\multicolumn{1}{P}{$\frac{c_1b }{1-c_1-c_2}$} &
\\ \hline

\multicolumn{1}{F}{$\om_b(B)$} &
\multicolumn{1}{A}{$\frac{\mu(1-c_1-c_2)+b\ell}{1-c_1-c_2+ b}$}   &
\multicolumn{1}{A}{$\frac{\mu(1-c_1-c_2)+b(\ell+c_2)}{1-c_1+c_2+ b}$} &

\end{tabular}
\bigskip

\caption{Inflation process to show that $\Aaucc \subset \Aa_{\mu',c_1,c_2}$ when $\mu' \leq \mu$.}
\label{Table4}
\end{table}

\begin{table}[!ht]
\setlength{\extrarowheight}{.25cm}
\begin{tabular}{l!{\vrule}ccccccc}
%
\multicolumn{1}{F}{Configuration} &
%
\multicolumn{1}{P}{$c_1 < \lambda \leq c_1+c_2, \ (14) \mbox{ or } (16)$} &
\multicolumn{1}{P}{$c_1+c_2 < \lambda, \ (17) \mbox{ or } (18)$} &
\\ \hline

\multicolumn{1}{F}{Curves} &
\multicolumn{1}{A}{$D_{4\ell+1},D_{4\ell+2},D_{4\ell+3} $}   &
\multicolumn{1}{A}{$D_{4\ell+1},D_{4\ell+3},D_{4\ell+6}$} &
\\ \hline

\multicolumn{1}{F}{$\omega_{abe}$} &
\multicolumn{1}{P}{$\frac{\omucc+ a d_{4\ell+2}+bd_{4\ell+1}+ed_{4\ell+3}}{1+a+b+e}$}   &
\multicolumn{1}{P}{$\frac{\omucc+ a d_{4\ell+6}+bd_{4\ell+1}+ed_{4\ell+3}}{1+a+b+e}$} &
\\ \hline

\multicolumn{1}{F}{$a$} &
\multicolumn{1}{A}{$\frac{c_2b}{1-c_1}$}   &
\multicolumn{1}{A}{$\frac{c_2b}{1-c_1}$} &
\\ \hline

\multicolumn{1}{F}{$e$} &
\multicolumn{1}{P}{$\frac{(c_1-c_2)b }{1-c_1}$}   &
\multicolumn{1}{P}{$\frac{(c_1-c_2)b }{1-c_1}$} &
\\ \hline

\multicolumn{1}{F}{$\om_b(B)$} &
\multicolumn{1}{A}{$\frac{\mu(1-c_1)+b(\ell+c_1)}{1-c_1+ b}$}   &
\multicolumn{1}{A}{$\frac{\mu(1-c_1)+b(\ell+c_1+c_2)}{1-c_1+ b}$} &

\end{tabular}
\bigskip

\caption{Inflation process to show that $\Aaucc \subset \Aa_{\mu',c_1,c_2}$ when $\mu' \leq \mu$.}
\label{Table5}
\end{table}

{\it Step 2}: First we show that $\Aa_{\mu,c_1',c_2} \subset \Aaucc$ whenever $c_1 \leq c_1'$. We consider the case $\lambda \leq c_2$ and almost complex structures $J$ such that the embedded $J$--holomorphic spheres satisfy configurations (8) or (10). Given $J\in \Aa_{\mu,c_1',c_2}$ we need to show that there is a symplectic form that tames $J$ and for which $B,F,E_1$ and $E_2$ have area $\mu,1,c_1$ and $c_2$ respectively. By Lemmas \ref{lm:F} and  \ref{representedclasses} there exist embedded $J$--holomorphic spheres representing the classes $D_{4\ell+1}=B+\ell F$, $D_{4\ell}=B+\ell F-E_2$ and $F$. The Inflation Lemma implies that the Poincaré duals of $D_{4\ell+1}$, $D_{4\ell}$ and $F$ are represented by 2-forms $d_{4\ell+1}$, $d_{4\ell}$ and $f$ such that the form 
\begin{equation}\label{omegac1}
 \om_{abe}=\frac{\om_{\mu,c_1',c_2}+ ad_{4\ell}+bd_{4\ell+1}+ef}{1+a+b}
\end{equation}
is symplectic for all $a,b,e \geq 0$ and all tame $J$.  We then have $ \om_{abe}(F)=1$, 
$$ \om_{abe}(B)=\frac{\mu + (a+b)\ell+e}{1+a+b} \ \ \mbox{and} \ \  \om_{abe}(E_2)=\frac{c_2+a}{1+a+b}.$$
Setting $e=\lambda b/(1-c_2) \geq 0$ and $a=c_2b/(1-c_2) \geq 0$, the resulting family of symplectic forms $\om_b$ satisfies $\om_b(F)=1$, $\om_b(B)=\mu$, $\om_b(E_2)=c_2$ and  $\om_b(E_1)=c_1'(1-c_2)/(1-c_2+b)$. Since $b$ can take any non negative value, this implies that $J \in \Aaucc$ for all $c_1 \leq c_1'$ in this particular case.  For all other cases the choice of curves to use in the inflation process is shown in Tables \ref{Table6}, \ref{Table7} and \ref{Table8}. Note that for almost complex structures $J$ such as in configurations (7), (9), (13) and (15) one again uses negative inflation. 

\begin{table}[!ht]
\setlength{\extrarowheight}{.25cm}
\begin{tabular}{l!{\vrule}ccccccc}
%
\multicolumn{1}{F}{Configuration} &
%
\multicolumn{1}{P}{$\lambda \leq c_2,\ (7) \mbox{ or } (9)$} &
\multicolumn{1}{P}{$c_1 < \lambda \leq c_1+c_2,\ (13) \mbox{ or } (15)$} &
\\ \hline

\multicolumn{1}{F}{Curves} &
\multicolumn{1}{A}{$F, D_{4\ell+1}, E_1-E_2$}   &
\multicolumn{1}{A}{$F, D_{4\ell+1}, F-E_1-E_2$} &
\\ \hline

\multicolumn{1}{F}{$\omega_{be}$} &
\multicolumn{1}{P}{$\frac{\om_{\mu,c_1',c_2}+bd_{4\ell+1}+ef}{1+b}$}   &
\multicolumn{1}{P}{$\frac{\om_{\mu,c_1',c_2}+bd_{4\ell+1}+ef}{1+b}$} &
\\ \hline

\multicolumn{1}{F}{$\omega_a$} &
\multicolumn{1}{A}{$\om_{be}+a {\rm PD}(E_1-E_2)$}   &
\multicolumn{1}{A}{$\om_{be}+a {\rm PD}(F-E_1-E_2)$} &
\\ \hline

\multicolumn{1}{F}{$a$} &
\multicolumn{1}{P}{$\frac{c_2b}{1+b}$}   &
\multicolumn{1}{P}{$\frac{c_2b}{1+b}$} &
\\ \hline

\multicolumn{1}{F}{$e$} &
\multicolumn{1}{A}{$\lambda b $}   &
\multicolumn{1}{A}{$(\lambda-c_2)b $} &
\\ \hline

\multicolumn{1}{F}{$\om_b(E_1)$} &
\multicolumn{1}{P}{$\frac{c_1'-c_2b}{1+ b}$}   &
\multicolumn{1}{P}{$\frac{c_1'+c_2b}{1+ b}$} &

\end{tabular}
\bigskip

\caption{Inflation process to show that $\Aa_{\mu,c_1',c_2} \subset \Aaucc$ when $c_1 \leq c_1'$.}
\label{Table6}
\end{table}

In  the first column of Table \ref{Table4} the parameter $b$ satisfies $0 \leq b \leq \frac{c_1'-c_2}{2c_2}$. This implies that $2a \leq \om_{be}(E_1-E_2)$ as required by negative inflation and $c_2 \leq \om_b(E_1) \leq c_1'$ as desired. 

\begin{table}[!ht]
\setlength{\extrarowheight}{.25cm}
\begin{tabular}{l!{\vrule}ccccccc}
%
\multicolumn{1}{F}{Configuration} &
%
\multicolumn{1}{P}{$\lambda \leq c_2, \ (8) \mbox{ or } (10)$} &
\multicolumn{1}{P}{$c_2 < \lambda \leq c_1, \ (11) \mbox{ or } (12)$} &
\\ \hline

\multicolumn{1}{F}{Curves} &
\multicolumn{1}{A}{$F,D_{4\ell},D_{4\ell+1}$}   &
\multicolumn{1}{A}{$F,D_{4\ell+1},D_{4\ell+4}$} &
\\ \hline

\multicolumn{1}{F}{$\omega_{abe}$} &
\multicolumn{1}{P}{$\frac{\om_{\mu,c_1',c_2}+ a d_{4\ell}+bd_{4\ell+1}+ef}{1+a+b}$}   &
\multicolumn{1}{P}{$\frac{\om_{\mu,c_1',c_2}+ a d_{4\ell+4}+bd_{4\ell+1}+ef}{1+a+b}$} &
\\ \hline

\multicolumn{1}{F}{$a$} &
\multicolumn{1}{A}{$\frac{c_2b}{1-c_2}$}   &
\multicolumn{1}{A}{$\frac{c_2b}{1-c_2}$} &
\\ \hline

\multicolumn{1}{F}{$e$} &
\multicolumn{1}{P}{$\frac{\lambda b }{1-c_2}$}   &
\multicolumn{1}{P}{$\frac{(\lambda-c_2)b }{1-c_2}$} &
\\ \hline

\multicolumn{1}{F}{$\om_b(E_1)$} &
\multicolumn{1}{A}{$\frac{c_1'(1-c_2)}{1-c_2+ b}$}   &
\multicolumn{1}{A}{$\frac{c_1'(1-c_2)}{1-c_2+ (\lambda-c_2)b}$} &

\end{tabular}
\bigskip

\caption{Inflation process to show that $\Aa_{\mu,c_1',c_2} \subset \Aaucc$ when $c_1 \leq c_1'$.}
\label{Table7}
\end{table}

\begin{table}[!ht]
\setlength{\extrarowheight}{.25cm}
\begin{tabular}{l!{\vrule}ccccccc}
%
\multicolumn{1}{F}{Configuration} &
%
\multicolumn{1}{P}{$c_1 \leq c_1' < \lambda \leq c_1+c_2, \ (14) \mbox{ or } (16)$} &
\multicolumn{1}{P}{$c_1+c_2 < \lambda, \ (17) \mbox{ or } (18)$} &
\\ \hline

\multicolumn{1}{F}{Curves} &
\multicolumn{1}{A}{$F, D_{4\ell+1},D_{4\ell+2} $}   &
\multicolumn{1}{A}{$F, D_{4\ell+1},D_{4\ell+6}$} &
\\ \hline

\multicolumn{1}{F}{$\omega_{abe}$} &
\multicolumn{1}{P}{$\frac{\om_{\mu,c_1',c_2}+ a d_{4\ell+2}+bd_{4\ell+1}+ef}{1+a+b}$}   &
\multicolumn{1}{P}{$\frac{\om_{\mu,c_1',c_2}+ a d_{4\ell+6}+bd_{4\ell+1}+ef}{1+a+b}$} &
\\ \hline

\multicolumn{1}{F}{$a$} &
\multicolumn{1}{A}{$\frac{c_2b}{1-c_2}$}   &
\multicolumn{1}{A}{$\frac{c_2b}{1-c_2}$} &
\\ \hline

\multicolumn{1}{F}{$e$} &
\multicolumn{1}{P}{$\frac{(\lambda-c_2)b }{1-c_2}$}   &
\multicolumn{1}{P}{$\frac{(\lambda-2c_2)b }{1-c_2}$} &
\\ \hline

\multicolumn{1}{F}{$\om_b(E_1)$} &
\multicolumn{1}{A}{$\frac{c_1'(1-c_2)+c_2b}{1-c_2+ b}$}   &
\multicolumn{1}{A}{$\frac{c_1'(1-c_2)+c_2b}{1-c_2+ b}$} &

\end{tabular}
\bigskip

\caption{Inflation process to show that $\Aa_{\mu,c_1',c_2} \subset \Aaucc$ when $c_1 \leq c_1'$.}
\label{Table8}
\end{table}

We now suppose that $\lambda \leq c_2$  and prove that  $ \Aaucc \subset \Aa_{\mu,c_1',c_2}$ whenever $c_1 \leq c_1'$. From Lemma \ref{representedclasses} we know that, for the configurations (7) and (9) we  can inflate along the classes $D_{4\ell-1}=B+\ell F-E_1$, $ F$ and $E_1-E_2$. Then, setting first
$$\om_{be}=\frac{\omucc+bd_{4\ell-1}+ef}{1+b}$$ and then 
$$\om_a=\om_{b,e}+a{\rm PD}(E_1-E_2),$$ 
where $a=c_2b$ and $e=\lambda b$, the resulting one parameter family of symplectic forms satisfy $\om_b(F)=1$, $\om_b(B)=\mu$, $\om_b(E_2)=c_2$ and  $\om_b(E_1)=(c_1+(1-c_2)b)/(1+b)$. Letting $b\to \infty$, this yields the desired result. Note that since $2c_2<1$ one has  $2a < \om_{be}(E_1-E_2)$, for all $b \geq 0$, as required by negative inflation. All remaining cases that show that $ \Aaucc \subset \Aa_{\mu,c_1',c_2}$ when $c_1 \leq c_1'$,  are summarized in the Tables \ref{Table9}, \ref{Table10} and \ref{Table11}. 

\begin{table}[!ht]
\setlength{\extrarowheight}{.25cm}
\begin{tabular}{l!{\vrule}ccccccc}
%
\multicolumn{1}{F}{Configuration} &
%
\multicolumn{1}{P}{$\lambda \leq c_2,\ (7) \mbox{ or } (9)$} &
\multicolumn{1}{P}{$c_1 \leq c_1'< \lambda \leq c_1+c_2,\ (13) \mbox{ or } (15)$} &
\\ \hline

\multicolumn{1}{F}{Curves} &
\multicolumn{1}{A}{$F, D_{4\ell-1}, E_1-E_2$}   &
\multicolumn{1}{A}{$ D_{4\ell+1}, D_{4\ell+3},F-E_1-E_2$} &
\\ \hline

\multicolumn{1}{F}{$\omega_{be}$} &
\multicolumn{1}{P}{$\frac{\omucc+bd_{4\ell-1}+ef}{1+b}$}   &
\multicolumn{1}{P}{$\frac{\omucc+bd_{4\ell+3}+ed_{4\ell+1}}{1+b+e}$} &
\\ \hline

\multicolumn{1}{F}{$\omega_a$} &
\multicolumn{1}{A}{$\om_{be}+a {\rm PD}(E_1-E_2)$}   &
\multicolumn{1}{A}{$\om_{be}+a {\rm PD}(F-E_1-E_2)$} &
\\ \hline

\multicolumn{1}{F}{$a$} &
\multicolumn{1}{P}{$\frac{c_2b}{1+b}$}   &
\multicolumn{1}{P}{$\frac{c_2b}{\lambda-c_2+b}$} &
\\ \hline

\multicolumn{1}{F}{$e$} &
\multicolumn{1}{A}{$\lambda b $}   &
\multicolumn{1}{A}{$\frac{(1-\lambda+c_2)b}{\lambda-c_2} $} &
\\ \hline

\multicolumn{1}{F}{$\om_b(E_1)$} &
\multicolumn{1}{P}{$\frac{c_1+(1-c_2)b}{1+ b}$}   &
\multicolumn{1}{P}{$\frac{c_1(\lambda -c_2)+\lambda b}{\lambda-c_2+ b}$} &

\end{tabular}
\bigskip

\caption{Inflation process to show that $ \Aaucc \subset \Aa_{\mu,c_1',c_2}$ when $c_1 \leq c_1'$.}
\label{Table9}
\end{table}

\begin{table}[!ht]
\setlength{\extrarowheight}{.25cm}
\begin{tabular}{l!{\vrule}ccccccc}
%
\multicolumn{1}{F}{Configuration} &
%
\multicolumn{1}{P}{$\lambda \leq c_2, \ (8) \mbox{ or } (10)$} &
\multicolumn{1}{P}{$c_2 < \lambda \leq c_1, \ (11) \mbox{ or } (12)$} &
\\ \hline

\multicolumn{1}{F}{Curves} &
\multicolumn{1}{A}{$F,D_{4\ell},D_{4\ell-1}$}   &
\multicolumn{1}{A}{$F,D_{4\ell-1},D_{4\ell+4}$} &
\\ \hline

\multicolumn{1}{F}{$\omega_{abe}$} &
\multicolumn{1}{P}{$\frac{\omucc+ a d_{4\ell}+bd_{4\ell-1}+ef}{1+a+b}$}   &
\multicolumn{1}{P}{$\frac{\omucc+ a d_{4\ell+4}+bd_{4\ell-1}+ef}{1+a+b}$} &
\\ \hline

\multicolumn{1}{F}{$a$} &
\multicolumn{1}{A}{$\frac{c_2b}{1-c_2}$}   &
\multicolumn{1}{A}{$\frac{c_2b}{1-c_2}$} &
\\ \hline

\multicolumn{1}{F}{$e$} &
\multicolumn{1}{P}{$\frac{\lambda b }{1-c_2}$}   &
\multicolumn{1}{P}{$\frac{(\lambda-c_2)b }{1-c_2}$} &
\\ \hline

\multicolumn{1}{F}{$\om_b(E_1)$} &
\multicolumn{1}{A}{$\frac{c_1(1-c_2)+(1-c_2)b}{1-c_2+ b}$}   &
\multicolumn{1}{A}{$\frac{c_1(1-c_2)+(1-c_2)b}{1-c_2+ b}$} &

\end{tabular}
\bigskip

\caption{Inflation process to show that $ \Aaucc \subset \Aa_{\mu,c_1',c_2}$ when $c_1 \leq c_1'$.}
\label{Table10}
\end{table}

\begin{table}[!ht]
\setlength{\extrarowheight}{.25cm}
\begin{tabular}{l!{\vrule}ccccccc}
%
\multicolumn{1}{F}{Configuration} &
%
\multicolumn{1}{P}{$c_1 \leq c_1' < \lambda \leq c_1+c_2, \ (14) \mbox{ or } (16)$} &
\multicolumn{1}{P}{$c_1 \leq c_1' < c_1+c_2 < \lambda, \ (17) \mbox{ or } (18)$} &
\\ \hline

\multicolumn{1}{F}{Curves} &
\multicolumn{1}{A}{$D_{4\ell+1},D_{4\ell+2},D_{4\ell+3}$}   &
\multicolumn{1}{A}{$D_{4\ell+1},D_{4\ell+3},D_{4\ell+6}$} &
\\ \hline

\multicolumn{1}{F}{$\omega_{abe}$} &
\multicolumn{1}{P}{$\frac{\omucc+ a d_{4\ell+2}+bd_{4\ell+3}+ed_{4\ell+1}}{1+a+b+e}$}   &
\multicolumn{1}{P}{$\frac{\omucc+ a d_{4\ell+6}+bd_{4\ell+3}+ed_{4\ell+1}}{1+a+b+e}$} &
\\ \hline

\multicolumn{1}{F}{$a$} &
\multicolumn{1}{A}{$\frac{c_2b}{\lambda-c_2}$}   &
\multicolumn{1}{A}{$\frac{c_2b}{\lambda-2c_2}$} &
\\ \hline

\multicolumn{1}{F}{$e$} &
\multicolumn{1}{P}{$\frac{(1-\lambda)b }{\lambda-c_2}$}   &
\multicolumn{1}{P}{$\frac{(1-\lambda+c_2)b }{\lambda-2c_2}$} &
\\ \hline

\multicolumn{1}{F}{$\om_b(E_1)$} &
\multicolumn{1}{A}{$\frac{c_1(\lambda-c_2)+\lambda b}{\lambda-c_2+ b}$}   &
\multicolumn{1}{A}{$\frac{c_1(\lambda-2c_2)+(\lambda-c_2) b}{\lambda-2c_2+ b}$} &

\end{tabular}
\bigskip

\caption{Inflation process to show that $ \Aaucc \subset \Aa_{\mu,c_1',c_2}$ when $c_1 \leq c_1'$.}
\label{Table11}
\end{table}

{\it Step 3}: The proof of this step is similar to the last one, although it is simpler when one shows  that $ \Aa_{\mu,c_1,c_2'} \subset \Aaucc$  whenever $c_2 \leq c_2'$, since we can inflate along the curve $E_2$. More precisely,  in this case it is always possible to inflate the form $\om_{\mu,c_1,c_2'}$ along an embedded $J$--holomorphic sphere representing the class $E_2$, to get a one parameter family of symplectic forms in classes $[\om_{\mu,c_1,c_2'}]+\epsilon {\rm PD}(E_2)=\om_{\mu, c_1,c_2'-\epsilon}$, where $0 \leq \epsilon < c_2'$. All these forms tame $J$. 

Next, in Tables \ref{Table12} to \ref{Table13}, we show that $\Aaucc \subset \Aa_{\mu,c_1,c_2'}$  whenever $c_2 \leq c_2'$.

\begin{table}[!ht]
\setlength{\extrarowheight}{.25cm}
\begin{tabular}{l!{\vrule}ccccccc}
%
\multicolumn{1}{F}{Configuration} &
%
\multicolumn{1}{P}{$\lambda \leq c_2 \leq c_2' < c_1,\ (7) \mbox{ or } (9)$} &
\multicolumn{1}{P}{$c_1 < \lambda \leq c_1+c_2,\ (13) \mbox{ or } (15)$} &
\\ \hline

\multicolumn{1}{F}{Curves} &
\multicolumn{1}{A}{$F, D_{4\ell-1}, E_1-E_2$}   &
\multicolumn{1}{A}{$F, D_{4\ell+1}, F-E_1-E_2$} &
\\ \hline

\multicolumn{1}{F}{$\omega_{be}$} &
\multicolumn{1}{P}{$\frac{\omucc+bd_{4\ell-1}+ef}{1+b}$}   &
\multicolumn{1}{P}{$\frac{\omucc+bd_{4\ell+1}+ef}{1+b}$} &
\\ \hline

\multicolumn{1}{F}{$\omega_a$} &
\multicolumn{1}{A}{$\om_{be}+a {\rm PD}(E_1-E_2)$}   &
\multicolumn{1}{A}{$\om_{be}+a {\rm PD}(F-E_1-E_2)$} &
\\ \hline

\multicolumn{1}{F}{$a$} &
\multicolumn{1}{P}{${\frac{(1-c_1)b}{1+b}}$}   &
\multicolumn{1}{P}{${\frac{c_1b}{1+b}}$} &
\\ \hline

\multicolumn{1}{F}{$e$} &
\multicolumn{1}{A}{$\lambda b $}   &
\multicolumn{1}{A}{$(\lambda-c_1)b $} &
\\ \hline

\multicolumn{1}{F}{$\om_b(E_2)$} &
\multicolumn{1}{P}{$\frac{c_2+(1-c_1)b}{1+ b}$}   &
\multicolumn{1}{P}{$\frac{c_2+c_1b}{1+ b}$} &

\end{tabular}
\bigskip

\caption{Inflation process to show that $\Aaucc \subset \Aa_{\mu,c_1,c_2'}$ when $c_2 \leq c_2'$.}
\label{Table12}
\end{table}

Note that in the first column of Table \ref{Table12}  if $c_1 \geq 1/2 \Leftrightarrow 1-c_1 \leq c_1$ then the negative inflation requirement $2a < \om_{be}(E_1-E_2)$ is satisfied for all $b \geq 0$, and this gives $c_2 \leq \om_b(E_2)< 1-c_1$.  On the other hand, if $c_1 \leq 1/2 \Leftrightarrow 1-c_1 \geq c_1 $ then the condition $2a < \om_{be}(E_1-E_2)$ implies that the parameter $b$ satisfies $0 \leq b \leq \frac{c_1-c_2}{1-2c_1}$ and therefore one has $c_2 \leq \om_{be}(E_2) \leq c_1$ as desired. There is a similar remark regarding the parameter $b$ in  the second column.

\begin{table}[!ht]
\setlength{\extrarowheight}{.25cm}
\begin{tabular}{l!{\vrule}ccccccc}
%
\multicolumn{1}{F}{Configuration} &
%
\multicolumn{1}{P}{$ \lambda \leq c_2 \leq c_2', \ (8) \mbox{ or } (10)$} &
\multicolumn{1}{P}{$c_2 \leq c_2' < c_1 < \lambda \leq c_1+c_2, \ (14) \mbox{ or } (16)$} &
\\ \hline

\multicolumn{1}{F}{Curves} &
\multicolumn{1}{A}{$F, D_{4\ell-1},D_{4\ell}$}   &
\multicolumn{1}{A}{$F, D_{4\ell+1},D_{4\ell+2}$} &
\\ \hline

\multicolumn{1}{F}{$\omega_{abe}$} &
\multicolumn{1}{P}{$\frac{\omucc+ a d_{4\ell}+bd_{4\ell-1}+ef}{1+a+b}$}   &
\multicolumn{1}{P}{$\frac{\omucc+ a d_{4\ell+2}+bd_{4\ell+1}+ef}{1+a+b}$} &
\\ \hline

\multicolumn{1}{F}{$a$} &
\multicolumn{1}{A}{$\frac{(1-c_1)b}{c_1}$}   &
\multicolumn{1}{A}{$\frac{c_1b}{\lambda-c_1}$} &
\\ \hline

\multicolumn{1}{F}{$e$} &
\multicolumn{1}{P}{$\frac{\lambda b}{c_1}$}   &
\multicolumn{1}{P}{$\frac{(\lambda-c_1)b}{1-c_1}$} &
\\ \hline

\multicolumn{1}{F}{$\om_b(E_2)$} &
\multicolumn{1}{A}{$\frac{c_2c_1+(1-c_1)b}{c_1+b}$}   &
\multicolumn{1}{A}{$\frac{c_2(1-c_1)+c_1 b}{1-c_1+b}$} &

\end{tabular}
\bigskip

\caption{Inflation process to show that $\Aaucc \subset \Aa_{\mu,c_1,c_2'}$ when $c_2 \leq c_2'$.}
\label{Table13}
\end{table}

For configurations (11) and (12) we need to consider two distinct cases: $\lambda \geq 1-c_1$ and $\lambda < 1-c_1$ (see Table \ref{Table14}).

\begin{table}[!ht]
\setlength{\extrarowheight}{.25cm}
\begin{tabular}{l!{\vrule}ccccccc}
%
\multicolumn{1}{F}{Configuration} &
%
\multicolumn{1}{P}{$  c_2 \leq c_2' < \lambda \leq c_1 \mbox{ and } \lambda \geq 1-c_1$} &
\multicolumn{1}{P}{$c_2 \leq c_2'  < \lambda \leq c_1   \mbox{ and } \lambda < 1-c_1 $} &
\\ \hline

\multicolumn{1}{F}{Curves} &
\multicolumn{1}{A}{$F, D_{4\ell-1},D_{4\ell+4}$}   &
\multicolumn{1}{A}{$D_{4\ell+1},D_{4\ell-1},D_{4\ell+4}$} &
\\ \hline

\multicolumn{1}{F}{$\omega_{abe}$} &
\multicolumn{1}{P}{$\frac{\omucc+ a d_{4\ell+4}+bd_{4\ell-1}+ef}{1+a+b}$}   &
\multicolumn{1}{P}{$\frac{\omucc+ a d_{4\ell+4}+bd_{4\ell-1}+ed_{4\ell+1}}{1+a+b+e}$} &
\\ \hline

\multicolumn{1}{F}{$a$} &
\multicolumn{1}{A}{$\frac{(1-c_1)b}{c_1}$}   &
\multicolumn{1}{A}{$\frac{\lambda b}{c_1}$} &
\\ \hline

\multicolumn{1}{F}{$e$} &
\multicolumn{1}{P}{$\frac{(\lambda -1+c_1)b}{c_1}$}   &
\multicolumn{1}{P}{$\frac{(1-\lambda-c_1)b}{c_1}$} &
\\ \hline

\multicolumn{1}{F}{$\om_b(E_2)$} &
\multicolumn{1}{A}{$\frac{c_2c_1+(1-c_1)b}{c_1+b}$}   &
\multicolumn{1}{A}{$\frac{c_2c_1+ \lambda b}{c_1+b}$} &

\end{tabular}
\bigskip

\caption{Inflation process to show that $\Aaucc \subset \Aa_{\mu,c_1,c_2'}$ when $c_2 \leq c_2'$.}
\label{Table14}
\end{table}

And finally for configurations (17) and (18) we also need to consider two distinct cases: $\lambda -c_1 \leq c_1 \Leftrightarrow \lambda \leq 2c_1$ and $c_1 < \lambda -c_1 \Leftrightarrow \lambda > 2c_1$ (see Table \ref{Table15}).

\begin{table}[!ht]
\setlength{\extrarowheight}{.25cm}
\begin{tabular}{l!{\vrule}ccccccc}
%
\multicolumn{1}{F}{Configuration} &
%
\multicolumn{1}{P}{$  c_2+c_1 < \lambda \mbox{ and } \lambda -c_1 \leq c_1$} &
\multicolumn{1}{P}{$c_2+c_1 < \lambda  \mbox{ and } \lambda -c_1 > c_1 $} &
\\ \hline

\multicolumn{1}{F}{Curves} &
\multicolumn{1}{A}{$F_1=F-E_1, D_{4\ell+1},D_{4\ell+6}$}   &
\multicolumn{1}{A}{$F, D_{4\ell+1},D_{4\ell+6}$} &
\\ \hline

\multicolumn{1}{F}{$\omega_{abe}$} &
\multicolumn{1}{P}{$\frac{\omucc+ a d_{4\ell+6}+bd_{4\ell+1}+ef_1}{1+a+b}$}   &
\multicolumn{1}{P}{$\frac{\omucc+ a d_{4\ell+6}+bd_{4\ell+1}+ef}{1+a+b}$} &
\\ \hline

\multicolumn{1}{F}{$a$} &
\multicolumn{1}{A}{$\frac{(\lambda -c_1)b}{1- \lambda +c_1}$}   &
\multicolumn{1}{A}{$\frac{c_1 b}{1-c_1}$} &
\\ \hline

\multicolumn{1}{F}{$e$} &
\multicolumn{1}{P}{$\frac{(2c_1-\lambda)b}{1- \lambda +c_1}$}   &
\multicolumn{1}{P}{$\frac{(\lambda-2c_1)b}{1-c_1}$} &
\\ \hline

\multicolumn{1}{F}{$\om_b(E_2)$} &
\multicolumn{1}{A}{$\frac{c_2(1-\lambda+c_1)+ (\lambda-c_1)b }{1- \lambda+c_1+b}$}   &
\multicolumn{1}{A}{$\frac{c_2(1-c_1)+c_1b }{1-c_1+b}$} &

\end{tabular}
\bigskip

\caption{Inflation process to show that $\Aaucc \subset \Aa_{\mu,c_1,c_2'}$ when $c_2 \leq c_2'$.}
\label{Table15}
\end{table}

\section{Compatible Integrable Complex Structures}

The purpose of this Appendix is to briefly explain some differential and topological results relative to the stratification of the space $\jj_{\om}:=\tjjucc$ of compatible almost complex structures on $\tMucc$, and to the stratification it induces on the subspace $\jj_{\om}^{int}\subset \jj_{\om}$ of compatible, integrable, complex structures. 

Recall from Section 2.2 that the space $\jj_{\om}$ is the disjoint union of finitely many strata, $\jj_{\om}=\sqcup U_{i}$, $i=0,\ldots,k$. Each stratum is characterized by the existence of a unique chain of holomorphic embedded spheres of one of the types (1) -- (18) and containing a curve in class $D_{-m}$, for some $m\geq 0$ (see Figures 2 to 6 and the definition of the classes $D_i$ in Section 2.2).
\begin{prop}
Suppose the stratum $U_{\Aa}\subset\jj_{\om}$ is characterized by the existence of a configuration of $J$-holomorphic embedded spheres $C_{1}\cup C_{2}\cup\cdots\cup C_{N}$ representing a given set of distinct homology classes $A_{1},\ldots, A_{N}$ of negative self-intersection. Then $U_{\Aa}$ is a cooriented Fréchet submanifold of $\jj_{\om}$ of (real) codimension $\codim_{\R}(U_{\Aa})= 2N-2c_{1}(A_{1}+\cdots+A_{N})$.
\end{prop}
\begin{proof}[Sketch of proof]

Given a symplectic $4$-manifold $(M,\om)$ and a set of $N$ distinct spherical homology classes 
\[
\Aa=\{A_{1},\ldots,A_{N}\}\subset H_{2}(M;\Z)
\]
let denote by $\Mm(\Aa,\jj_{\om})$ the space of $(N+1)$-tuples
\[(u,J):=(u_{1},\ldots,u_{N},J)\in C^{\infty}\left((\CP^{1})^{N},M\right)\times \jj_{\om}\] 
such that $u_{i}:\CP^{1}\to M$ is a somewhere injective $J$-holomorphic map whose image represents the homology class $A_{i}$. By Proposition~6.2.7 in~\cite{MS:J-HolomorphicCurves2}, this space is always a smooth infinite dimensional Fréchet manifold. Its tangent space at $(u,J)$ can be identified with the vector space of $(N+1)$-tuples
\[ 
(\xi,Y):=(\xi_{1},\ldots,\xi_{N},Y)\in \Omega^{0}(u_{1}^{*}(TM))\times\cdots\times\Omega^{0}(u_{N}^{*}(TM))\times S\Omega_{J}^{0,1}(TM)\]
such that
\[
\delbar_{J}\xi + u_{i}^{*}Y = 0,~\forall i
\]
where $S\Omega_{J}^{0,1}(TM)$ is the image of the complex symmetric $(0,2)$-tensors $S^{0,2}_{J}(M)\subset T^{(0,2)}_{J}(M)$ under the natural identification $T^{(0,2)}_{J}(M)\simeq\Omega^{(0,1)}_{J}(TM)$ defined by the Hermitian metric $h_{J}(v,w):=\om(v,Jw)-i\om(v,w)$, and where $\delbar_{J}:\Omega^{0}(u^{*}(TM))\to\Omega^{0,1}_{J}(u^{*}(TM))$ is a $0$-order perturbation of a genuine Dolbeault operator. The image under the projection $\pi:\Mm(A,\jj_{\om})\to\jj_{\om}$ is the set $U_{\Aa}$ of all $J$ such that each class in $\Aa$ is represented by an irreducible $J$-holomorphic sphere. The differential is given by the projection
\begin{align*}
\pi_{*}: T_{(u,J)}\Mm(\Aa,\jj_{\om}) &\to T_{J}\jj_{\om}\\
(\xi,Y) &\mapsto Y
\end{align*}
and is Fredholm of index $2c_{1}(A_{1}+\cdots+A_{N})+4N$. Its image at $(u,J)$ is
\begin{align*}
\image(\pi_{*}) &= \left\{ Y~|~u_{i}^{*}(Y)\in\image\left(\delbar_{J}:\Omega^{0}(u_{i}^{*}(TM))\to \Omega^{0,1}_{J}(u_{i}^{*}(TM)\right),~\forall i\right\}\\
&\simeq \image\left(\oplus\delbar_{J}:\bigoplus \Omega^{0}(u_{i}^{*}(TM))\to \bigoplus\Omega^{0,1}_{J}(u_{i}^{*}(TM))\right)
\end{align*}
while its kernel is isomorphic to $\ker(\oplus\delbar_{J})$. Moreover, because the classes $A_{i}$ have negative self-intersection, each element in the preimage $\pi^{-1}(J)$ is obtained from a single $J$-holomorphic map $u:(\CP^{1})^{N}\to M$ by reparametrization of each component under the free action of $\PSL(2,\C)$. It follows that $\ker(\pi_{*})\simeq\ker(\oplus\delbar_{J})$ is isomorphic to $(\psl(2,\C))^{N}$. Consequently, $\coker(\pi_{*})$ has constant dimension $6N - \left(2c_{1}(A_{1}+\cdots A_{N})+4N\right)$ and the projection $\pi$ factors through a smooth embedding
\[
\Mm(\Aa, \jj_{\om}) / (\PSL(2,\C))^{N} \into \jj_{\om}
\]
To show that the bundle $\coker(\pi_{*})\to U_{i}$ is oriented, we first note that $\pi_{*}$ lies in a component of Fredholm operators that contains a complex linear operator (namely the projection $\pi_{*}$ at an integrable $J$). Therefore, the determinant bundle $\Lambda^{\max}\ker(\pi_{*})\otimes\Lambda^{\max}\coker(\pi_{*})$ is oriented. Since $\ker(\pi_{*})\simeq\psl(2,\C)$ is also oriented, this defines an orientation on the cokernel bundle. For more details, see Appendix A in~\cite{Ab}, Proposition~2.8 in~\cite{AGK}, and Proposition 6.2.7 in~\cite{MS:J-HolomorphicCurves2}.
\end{proof}
We now turn our attention to the subspace $\jj_{\om}^{int}\subset\jj_{\om}$ of compatible, integrable, complex structures. Given a stratum $U_{i}\subset\jj_{\om}$, let write $V_{i}=U_{i}\cap\jj_{\om}^{int}$.

\begin{lemma}\label{lemma:Hirzebruch}
Given any $J\in\jj_{\om}^{int}$, the complex surface $(\X_{3}, J)$ is the twofold blow-up of a Hirzebruch surface $(\F_{m},J_{m})$.
\end{lemma}
\begin{proof}
This follows from the fact that given any $J\in\jj_{\om}^{int}$, the class $E_{2}$ is always represented by an exceptional complex curve that we can blow-down, and that on the resulting surface $(\X_{2},J)$, the exceptional classes $E_{1}$ and $F-E_{1}$ are both represented. The blow-down along the class $E_{1}$ yields an even Hirzebruch surface $\F_{2k}$ while the blow-down  along the class $F-E_{1}$ yields an odd Hirzebruch surface $\F_{2k-1}$.
\end{proof}

\noindent From a complex analytic (or algebraic) point of view, the 18 types of compatible complex structures on $\X_{3}$ can be constructed as follows:
\begin{description}\setlength{\itemsep}{1mm}
\item [Type (1)] Twofold blow-up of $\F_{0}$ at two generic points (not lying on the same fiber $F$ nor on the same section $B$).
\item [Type (2)] Twofold blow-up of $\F_{0}$ at two distinct points on the same fiber $F$. 
\item [Type (3)] Twofold blow-up of $\F_{0}$ at two distinct points on the same section $B$.
\item [Type (4)] Twofold blow-up of $\F_{0}$ at two ``infinitely near'' points on a fiber, that is, the blow-up of $\F_{0}$ at $p$ followed by the blow-up of $\widetilde{\F}_{0}$ at the line $\ell_{p}=T_{p}F\subset T_{p}\F_{0}$ on the exceptional divisor.
\item [Type (5)] Twofold blow-up of $\F_{0}$ at two ``infinitely near'' points on a flat section $B$, that is, at $(p, \ell_{p}=T_{p}B)$.
\item [Type (6)] Twofold blow-up of $\F_{0}$ at two ``infinitely near'' points $(p, \ell_{p})$ with the direction $\ell_{p}$ transverse to $T_{p}F$ and $T_{p}B$.
\item [Type (7)] Twofold blow-up of $\F_{2m-1}$ at two distinct points on the same fiber, one of which lying on the section $s_{0}$ of self-intersection $-2m+1$. 
\item [Type (8)] Twofold blow-up of $\F_{2m-1}$ at two ``infinitely near'' points $(p,\ell_{p})$, where $p\in s_{0}$ and $\ell_{p}$ is tranverse to $s_{0}$ and to the fiber $F_{p}$. 
\item [Type (9)] Twofold blow-up of $\F_{2m-1}$ at two ``infinitely near'' points $(p,\ell_{p})$, where $p\in s_{0}$ and $\ell_{p}=T_{p}F$.
\item [Type (10)] Twofold blow-up of $\F_{2m}$ at two generic points, that is, at two points on two different fibers and away from $s_{0}$.
\item [Type (11)] Twofold blow-up of $\F_{2m}$ at two distinct points on the same fiber, one of which belonging to $s_{0}$.
\item [Type (12)] Twofold blow-up of $\F_{2m}$ at two points on different fibers, one points belonging to $s_{0}$.
\item [Type (13)] Twofold blow-up of $\F_{2m}$ at two distinct points on the same fiber, one of which belonging to $s_{0}$ (note that only the order of the two blow-ups distinguishes (13) from (11)).
\item [Type (14)] Twofold blow-up of $\F_{2m}$ at two ``infinitely near'' points $(p,\ell_{p})$, where $p\in s_{0}$ and $\ell_{p}$ is tranverse to $s_{0}$ and to the fiber $F_{p}$.
\item [Type (15)] Twofold blow-up of $\F_{2m}$ at two ``infinitely near'' points $(p,\ell_{p})$, where $p\in s_{0}$ and $\ell_{p}=T_{p}F$.
\item [Type (16)] Twofold blow-up of $\F_{2m}$ at two points on different fibers, one points belonging to $s_{0}$ (note that only the order of the two blow-ups distinguishes (16) from (12)).
\item [Type (17)] Twofold blow-up of $\F_{2m}$ at two ``infinitely near'' points $(p,\ell_{p})$, where $p\in s_{0}$ and $\ell_{p}=T_{p}s_{0}$.
\item [Type (18)] Twofold blow-up of $\F_{2m}$ at two distinct points on $s_{0}$.
\end{description}

\noindent Let $\Diff_{h}$ denote the group of diffeomorphisms of $\X_{3}$ acting trivially on homology.

\begin{prop}\label{prop:Teichmuller}
Given any two compatible complex structures $J_{1}, J_{2}$ in $V_{i}$, there exists a diffeomorphism $\phi$ acting trivially on homology such that $J_{2}=\phi_{*}J_{1}$. Consequently, given any $J_{i}\in V_{i}$, we have 
\[V_{i} = U_{i}\cap\jj_{\om}^{int}=(\Diff_{h}\cdot J_{i})\cap \jj_{\om}\]
\end{prop}
\begin{proof}
First recall that the complex automorphism group of $\F_{0}$ is isomorphic to 
$$\Aut(\F_{0})\simeq(\PSL(2,\C)\times \PSL(2,\C))\ltimes \Z_{2}$$
while the automorphism group of $\F_{m}$, $m\geq 1$, is isomorphic to the semi-direct product 
$$\Aut(\F_{m}) \simeq GL(2,\C)/\mu_{m}\ltimes H^{0}(\CP^{1}, (\Lambda^{m})^{*})$$
where $\mu_{m}$ is the group of $m$\textsuperscript{th} roots of unity.  The action of $GL(2,\C)/\mu_{m}$ lifts the action of $\PSL(2,\C)$ on $\CP^{1}$ and thus acts triply transitively on the set of fibers of $\F_{m}$. This action preserves the zero section $s_{0}$ (of self-intersection $-m$) and the section at infinity $s_{\infty}$ (of self-intersection $+m$) and is transitive on their complement $\F_{m}\setminus\{s_{0},s_{\infty}\}$. The group of sections of the dual line bundle $(\Lambda^{m})^{*}$ is isomorphic to the space $\Sym^{m}(\C^{2})$ of symmetric $m$-linear forms on $\C^{2}$ -- that is, homogenous polynomials of degree $m$ in two variables -- and acts fiberwise by 
$$\alpha\cdot[z:u^{\otimes m}] =  [z+\alpha(u^{\otimes m}):u^{\otimes m}]$$
Once restricted to a single fiber $F_{p}\simeq \CP^{2}$, the action of $\Sym^{m}(\C^{2})$ fixes the point $s_{0}\cap F_{p}=[1:0]$ while it is simply transitive on $F_{p}\setminus[1:0]$. We note that for $m\geq1$, these automorphism groups are all connected and, consequently, they act trivially on homology.

Now, the proof of the statement reduces to showing that the appropriate automorphism group acts transitively on the pairs of points $(p_{1},p_{2})$ or $(p,\ell_{p})$ that define the compatible complex structures of a given type (1) -- (18) as blow-ups of Hirzebruch structures. \end{proof}

In their paper~\cite{AGK}, Abreu-Granja-Kitchloo proved that, under some cohomological conditions, $\jj^{int}$ is a genuine Fr\'echet submanifold of $\jj$ whose tangent bundle may be described using standard deformation theory. In order to state their result, let write $H^{0,q}_{J}(M)$ for the $q^{\text{th}}$ Dolbeault cohomology group with coefficients in the sheaf of germs of holomorphic functions, and $H^{0,q}_{J}(TM)$ for the $q^{\text{th}}$ Dolbeault cohomology group with coefficients in the sheaf of germs of holomorphic vector fields. Then, 
\begin{thm}[\cite{AGK} Theorem 2.3]\label{thmA}
If $(M,\om)$ is a symplectic $4$-manifold, $J\in\jj^{int}_{\om}$ is a compatible integrable complex structure, and the cohomology groups $H^{0,2}_{J}(M)$ and $H^{0,2}_{J}(TM)$ are zero, then $\jj^{int}_{\om}$ is a submanifold of $\jj_{\om}$ in the neighborhood of $J$. Moreover, the moduli space of infinitesimal compatible deformations of $J$ in $\jj^{int}_{\om}$ coincides with the moduli space of infinitesimal deformations of $J$ in the set of all integrable structures, that is, it is given by $H^{0,1}_{J}(TM)$. Finally, the tangent space of $\jj^{int}_{\om}$ at $J$ is naturally identified with   

$T_{J}((\Diff\cdot J)\cap \jj^{int}_{\om})\oplus H^{0,1}_{J}(TM)$.
\end{thm}

\begin{lemma}\label{le:Vanishing}
For any $J\in\jj_{\om}^{int}$, the cohomology groups $H^{0,2}_{J}(\X_{3})$ and $H^{0,2}_{J}(T\X_{3})$ are zero. 
\end{lemma}
\begin{proof}
(See also~\cite{Ko} \S 5.2(a)(iv) p.220.) By Lemma~\ref{lemma:Hirzebruch}, we know that any compatible complex structure $J\in\jj_{\om}^{int}$ is obtained by blowing-up a Hirzebruch structure on $\STS$. 
Now, for any complex surface $(X,J)$, the geometric genus $p_{g}:=\rk H^{0,2}_{J}(X)$ is a   birational invariant. Consequently, the first assertion follows from the classical fact that $p_{g}(\F_{m})=0$, for all $m\geq 0$.

As for $H^{0,2}_{J}(T\X_{3}):=\check{H}^{2}_{J}(T\X_{3})$, Serre duality implies that $\check{H}^{2}_{J}(T\X_{3})^{\vee}\simeq \check{H}^{0}(\mathcal{K}_{J}\otimes\Omega^{1}_{J})$. Now, by Lemma~\ref{lm:E2}, $(\X_{3},J)$ contains a two-dimensional family of embedded rational curves of zero self-intersection (the fibers in class $[F]$) which cover a dense open set, and the restriction of the rank $2$ bundle $\mathcal{K}_{J}\otimes\Omega^{1}_{J}$ to any of those curves is isomorphic to $\Oo(-4)\oplus\Oo(-2)$. Hence, $\mathcal{K}_{J}\otimes\Omega^{1}_{J}$ cannot have any nontrivial holomorphic sections. 
\end{proof}

Let $\Aut_{h}(J)\subset\Diff_{h}$ for the subgroup of complex automorphism of $(M,J)$. Let $\Iso_{h}(\om,J)\subset \Aut_{h}(J)$ denote the K\"ahler isometry group of $(M,\om,J)$. The next result shows that in some cases the part of the $\Diff_{h}$-orbit of $J$ which lies in $\jj^{int}_{\om}$ may be identified with the $\Symp_{h}(M,\om)$-orbit:

\begin{thm}[\cite{AGK} Corollary 2.6]\label{thmB}
If $J\in\jj^{int}_{\om}$ is such that the inclusion $\Iso_{h}(\om,J)\into\Aut_{h}(J)$ is a weak homotopy equivalence, then the inclusion of the $\Symp_{h}(M,\om)$-orbit of $J$ in $(\Diff_{h}\cdot J)\cap\jj^{int}_{\om}$ 
\[
\Symp_{h}(M,\om)/\Iso_{h}(\om,J)\into (\Diff_{h}\cdot J)\cap\jj^{int}_{\om}
\] 
is also a weak homotopy equivalence.
\end{thm}
\begin{remark}
The actual statement of \cite{AGK}~Corollary 2.6 gives a condition for the full $\Symp(M,\om)$ orbit of $J\in\jj_{\om}^{int}$ to be homotopy equivalent to its $\Diff_{[\om]}$ orbit, where $\Diff_{[\om]}$ is the group of diffeomorphisms preserving the class $[\om]$. However, the same arguments apply to the orbits of $J$ under $\Symp_{h}$ and $\Diff_{h}$.
\end{remark}

\begin{lemma}\label{le:isometries}
Given any $J\in\jj^{int}(\tMucc)$, the inclusion $\Iso_{h}(\omucc,J)\into\Hol_{h}(J)$ is a weak homotopy equivalence.
\end{lemma}
\begin{proof}
First observe that if $(\tM,\tJ)$ is the blow-up of $(M,J)$ at a point $p$, then the automorphism group $\Hol_{h}(\tJ)$ is isomorphic to the stabilizer subgroup of $p$ in $\Hol_{h}(J)$. Looking at the 18 possible types of compatible complex structures defined on $\X_{3}\simeq \tMucc$, it is easy, but tedious, to check that the complex automorphism groups $\Hol_{h}(\X_{3},J)$ are homotopy equivalent to
\[
\Aut_{h}(\X_{3},J)\simeq
\begin{cases}
S^{1} & \text{if $J$ is of type 6, 8, or 14},\\
T^{2} & \text{otherwise.}
\end{cases}
\]
On the other hand, the isometry groups of the Hirzebruch surfaces (with respect to any K\"ahler form $\om$) are the maximal compact Lie subgroups of their complex automorphism groups, namely
\[
\Iso(\F_{m},\om) \simeq 
\begin{cases}
(\SO(3)\times \SO(3))\ltimes \Z_{2} & \text{if~} m=0,\\
S^{1}\times \SO(3) & \text{if $m$ is even},\\
\U(2) & \text{if $m$ is odd}.
\end{cases}
\]
In particular, they are deformation retracts of the automorphism groups $\Aut(\F_{m})$. After blow-up, they induce isometry groups $\Iso_{h}(\X_{3},\om,J)$ isomorphic to
\[
\Iso_{h}(\X_{3},\om,J)\simeq
\begin{cases}
S^{1} & \text{if $J$ is of type 6, 8, or 14},\\
T^{2} & \text{otherwise.}
\end{cases}
\] 
\end{proof}
Together with Proposition~\ref{prop:Teichmuller}, this gives
\begin{cor}
Given $J\in V_{i}\subset \jj_{\om}^{int}$, there is a weak homotopy equivalence
\[
\Symp_{h}(\tMucc,\om)/\Iso_{h}(\om,J) \simeq V_{i}
\]
\end{cor}
\noindent Let $\Aa=\{A_{1},\ldots,A_{N}\}$ be a set of distinct spherical homology classes of negative self-intersections. Let $U_{\Aa}$ be the stratum it defines in $\jj_{\om}$. The next proposition gives conditions ensuring that $U_{\Aa}$ is tranversal to $\jj_{\om}^{int}$ and that its normal bundle at $J\in\jj_{\om}^{int}$ may be described in terms of deformation theory.

\begin{thm}[\cite{AGK} Theorem 2.9]\label{thmC}
Let $(M,\om,J)$ be a K\"ahler $4$-manifold such that $J\in V_{\Aa}:=U_{\Aa}\cap\jj_{\om}^{int}$ and that $V_{\Aa}=(\Diff_{[\om]}\cdot J)\cap\jj_{\om}$. Suppose that the cohomology groups $H^{0,2}_{J}(M)$ and $H^{0,2}_{J}(TM)$ are zero. Suppose also that $(u,J):=(u_{1},\ldots,u_{N},J)\in\Mm(\Aa,\jj_{\om})$ is such that $u^{*}:H^{0,1}_{J}(TM)\to \bigoplus H^{0,1}(u_{i}^{*}(TM))$ is surjective. Then the projection $\pi:\Mm(\Aa,\jj_{\om})\to\jj_{\om}$ is tranversal at $(u,J)$ to $\jj_{\om}^{int}\subset\jj_{\om}$ and the infinitesimal complement to the image $U_{\Aa}$ of $\pi$ in a neighborhood of $J$ can be identified with the moduli space of infinitesimal deformations $H^{0,1}_{J}(TM)$.
\end{thm}

\begin{lemma}\label{le:Transversality}
Given $J\in V_{\Aa}\subset\jj^{int}(\X_{3},\om)$, let denote by $C$ the unique $J$-holomorphic configuration of type $\Aa$, and let $u=(u_{1},\ldots,u_{N})$ be some $J$-holomorphic parametrization of $C$. Then the induced map $u^{*}:H^{0,1}(T\X_{3})\to \bigoplus H^{0,1}(u^{*}(T\X_{3}))$ is surjective.
\end{lemma}
\begin{proof}
Consider the exact sequence of sheaves associated to the inclusion $f:C\into\X_{3}$ of $C$ viewed as a nodal curve:
\[
0\to \Oo_{\X_{3}}(-C)\to \Oo_{\X_{3}}\to \Oo_{C}\to 0
\]
Tensoring with $T\X_{3}$ we get the short exact sequence
\[
0\to \Oo_{\X_{3}}(-C)\otimes T\X_{3}\to T\X_{3}\to \Oo_{C}\otimes T\X_{3}\simeq f_{*}f^{*}T\X_{3}\to 0
\]
whose associated cohomology sequence is
\[
\cdots\to H^{1}(\X_{3};TX_{3})\stackrel{f^{*}}{\to} H^{1}(C;f^{*}T\X_{3})\to H^{2}(\X_{3};\Oo_{\X_{3}}(-C)\otimes T\X_{3})\to\cdots
\]
The sheaf $\Oo_{\X_{3}}(-C)\otimes T\X_{3}$ being locally free,  implies, by Serre duality, that $H^{2}(\X_{3};\Oo_{\X_{3}}(-C)\otimes T\X_{3})\simeq H^{0}(\X_{3};\Oo_{\X_{3}}(C)\otimes T\X_{3}^{\vee}\otimes K_{\X_{3}}.$ But, since $C\cdot F = 1$, the restriction of $\Oo_{\X_{3}}(C)\otimes T\X_{3}^{\vee}\otimes K_{\X_{3}}$ to any fiber $F$ is isomorphic to $\Oo(1)\otimes(\Oo(-4)\oplus\Oo(-2))$. Since the fibers cover an open dense subset of $\X_{3}$, it follows that $\Oo_{\X_{3}}(C)\otimes T\X_{3}^{\vee}\otimes K_{\X_{3}}$ has no nontrivial sections and, by duality, that $H^{2}(\X_{3};\Oo_{\X_{3}}(-C)\otimes T\X_{3})=0$.

Let $\nu:\widetilde{C}\to C$ be the normalization defined by $u=f\circ \nu$. We claim that $\nu^{*}: H^{1}(C;f^{*}T\X_{3})\to H^{1}(\widetilde{C};\nu^{*}f^{*}T\X_{3})$ is surjective. To see this, we observe that $H^{1}(\widetilde{C};\nu^{*}f^{*}T\X_{3})\simeq H^{1}(C;\nu_{*}\nu^{*}f^{*}T\X_{3})$ since $\nu$ is a finite-to-one proper map. Now consider the short exact sequence 
\[
0\to f^{*}T\X_{3}\to\nu_{*}\nu^{*}f^{*}T\X_{3}\to S \to 0
\]
Since the cokernel $S$ is supported on a finite number of points of $C$, it follows that $H^{1}(C,S)=0$, so that $H^{1}(C,f^{*}T\X_{3})\to H^{1}(C,\nu_{*}\nu^{*}f^{*}T\X_{3})$ is surjective.
\end{proof}
\begin{cor}
The action of $\Symp(\tMucc)$ on $\jj_{\om}$ is homotopy equivalent to its restriction to $\jj_{\om}^{int}$.
\end{cor} 
\begin{cor}\label{cor:Contractibility}
The space $\jj_{\om}^{int}$ of compatible integrable complex structures on $\tMucc$ is contractible.
\end{cor}


\end{document}